\documentclass[a4paper,12pt]{amsart}
\usepackage{amssymb,amsmath,a4wide,graphicx,stmaryrd,fullpage,setspace,microtype}
\usepackage[pagebackref=true]{hyperref}
\usepackage[all]{xy}
\usepackage{etoolbox}
\apptocmd{\thebibliography}{\raggedright}{}{}

\title{Combinatorics of past-similarity in higher dimensional transition systems}

\author[P. Gaucher]{Philippe Gaucher}

\address{Institut de Recherche en Informatique Fondamentale\\
  CNRS et Universit\'e Paris Diderot\\
 Case 7014\\ 75205 PARIS Cedex 13 \\ France}

\urladdr{http://www.irif.fr/{\~{}}gaucher} 

\subjclass{18C35,55U35,18G55,68Q85}

\keywords{left determined model category, combinatorial model category, discrete model structure, higher dimensional transition system, causal structure, bisimulation}

\swapnumbers

\newcommand{\C}{\mathcal{C}}

\newcommand{\K}{\mathcal{K}}

\newcommand{\I}{\mathcal{I}}

\newcommand{\de}{\partial}
\newcommand{\p}\times

\newtheorem{thm}{Theorem}[section]
\newtheorem{prop}[thm]{Proposition}
\newtheorem{lem}[thm]{Lemma}
\newtheorem{cor}[thm]{Corollary}

\newtheorem{defn}[thm]{Definition}
\newtheorem{nota}[thm]{Notation}
\newtheorem{rem}[thm]{Remark}
\newcommand{\bd}{\begin{defn}}
\newcommand{\ed}{\end{defn}}
\newcommand{\bp}{\begin{prop}}
\newcommand{\ep}{\end{prop}}
\newcommand{\bth}{\begin{thm}}
\renewcommand{\eth}{\end{thm}}
\newcommand{\bpf}{\begin{proof}}
\newcommand{\epf}{\end{proof}}
\newcommand{\bc}{\begin{cor}}
\newcommand{\ec}{\end{cor}}

\newtheorem{innercustomthm}{Theorem}
\newenvironment{customthm}[1]
  {\renewcommand\theinnercustomthm{#1}\innercustomthm}
  {\endinnercustomthm}

\newcommand{\fL}[1]{\ar@{->}[ll]_-{#1}}
\newcommand{\fR}[1]{\ar@{->}[rr]^-{#1}}
\newcommand{\fRr}[1]{\ar@{->}[rrr]^-{#1}}
\newcommand{\fD}[1]{\ar@{->}[dd]_-{#1}}
\newcommand{\fU}[1]{\ar@{->}[uu]^-{#1}}
\newcommand{\f}[2]{\ar@{->}[#1]|{#2}}
\newcommand{\ff}[2]{\ar@2{->}[#1]|{#2}}
\newcommand{\frr}[1]{\ar@{->}[rrrr]^-{#1}}

\newcommand{\fl}[1]{\ar@{->}[l]_-{#1}}
\newcommand{\fr}[1]{\ar@{->}[r]^-{#1}}
\newcommand{\fd}[1]{\ar@{->}[d]_-{#1}}
\newcommand{\fu}[1]{\ar@{->}[u]^-{#1}}

\newcommand{\iso}{\cong}

\renewcommand{\leq}{\leqslant}
\renewcommand{\geq}{\geqslant}
\newcommand{\dd}[1]{\uparrow\!\!{#1}\!\!\uparrow}

\def\cartesien{%
  \ar@{-}[]+R+<6pt,-2pt>;[]+RD+<6pt,-6pt>%
  \ar@{-}[]+D+<2pt,-6pt>;[]+RD+<6pt,-6pt>%
}
\def\cocartesien{%
  \ar@{-}[]+L+<-6pt,+2pt>;[]+LU+<-6pt,+6pt>%
  \ar@{-}[]+U+<-2pt,+6pt>;[]+LU+<-6pt,+6pt>%
}
\def\hocartesien{%
  \ar@{-}[]+R+<6pt,-2pt>;[]+RD+<6pt,-6pt>_{h}%
  \ar@{-}[]+D+<2pt,-6pt>;[]+RD+<6pt,-6pt>%
}
\def\hococartesien{%
  \ar@{-}[]+L+<-6pt,+2pt>;[]+LU+<-6pt,+6pt>_{h}%
  \ar@{-}[]+U+<-2pt,+6pt>;[]+LU+<-6pt,+6pt>%
}

\newcommand{\brm}[1]{\rm{\mathbf{#1}}}

\newcommand{\set}{{\brm{Set}}}

\DeclareMathOperator{\id}{Id}

\DeclareMathOperator{\Mor}{Mor}

\newcommand{\liminj}{\varinjlim}

\newcommand{\wts}{\mathcal{W\!T\!S}}
\newcommand{\cts}{\mathcal{C\!T\!S}}

\newcommand{\rts}{\mathcal{R\!T\!S}}
\newcommand{\csts}{\mathcal{C\!S\!T\!S}}

\DeclareMathOperator{\Int}{\underline{I}}
\DeclareMathOperator{\Ext}{\underline{E}}
\DeclareMathOperator{\St}{\underline{S}}
\DeclareMathOperator{\Ac}{\underline{L}}
\DeclareMathOperator{\Tr}{\underline{T}}

\makeatletter
\def\varholim@#1#2{%
  \vtop{\m@th\ialign{##\cr
    \hfil$#1\operator@font holim$\hfil\cr
    \noalign{\nointerlineskip\kern1.5\ex@}#2\cr
    \noalign{\nointerlineskip\kern-\ex@}\cr}}%
}
\def\holimproj{%
  \mathop{\mathpalette\varholim@{\leftarrowfill@\textstyle}}\nmlimits@
}
\def\holiminj{%
  \mathop{\mathpalette\varholim@{\rightarrowfill@\textstyle}}\nmlimits@
}
\makeatother

\DeclareMathOperator{\cell}{{\brm{cell}}}
\DeclareMathOperator{\cof}{{\brm{cof}}}

\DeclareMathOperator{\inj}{{\brm{inj}}}

\DeclareMathOperator{\cyl}{{Cyl}}
\DeclareMathOperator{\cocyl}{{Path}}
\DeclareMathOperator{\pscocyl}{{PseudoPath}}
\DeclareMathOperator{\CSA}{CSA}

\begin{document}

\begin{abstract} 
  The key notion to understand the left determined Olschok model
  category of star-shaped Cattani-Sassone transition systems is
  past-similarity.  Two states are past-similar if they have homotopic
  pasts. An object is fibrant if and only if the set of transitions is
  closed under past-similarity. A map is a weak equivalence if and
  only if it induces an isomorphism after the identification of all
  past-similar states. The last part of this paper is a discussion
  about the link between causality and homotopy.
\end{abstract}

\maketitle
\tableofcontents

\section{Introduction}

\subsection{Presentation}

This work belongs to our series of papers devoted to \emph{higher
  dimensional transition systems} \cite{hdts} \cite{cubicalhdts}
\cite{erratum_cubicalhdts} \cite{homotopyprecubical} \cite{biscsts1}
\cite{csts}. One of the goal of this series of papers is to explore
the link between causality and homotopy in this setting.

\begin{figure}
\[
\xymatrix{
& \beta \ar@{->}[rd]^-{\{v\}}&\\
\alpha\ar@{->}[ru]^-{\{u\}}\ar@{->}[rd]_-{\{v\}} & \{u,v\} & \delta\\
&\gamma\ar@{->}[ru]_-{\{u\}}&}
\]
\caption{$u|| v$: Concurrent execution of $u$ and $v$}
\label{concab}
\end{figure}
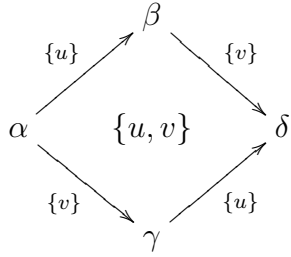

The notion of \emph{higher dimensional transition system} is a higher
dimensional analogue of the computer-scientific notion of labelled
transition system. The purpose is to model the concurrent execution of
$n$ actions by a multiset of actions, i.e. a set with a possible
repetition of some elements (e.g. $\{u,u,v,w,w,w\}$).  In the language
of Cattani and Sassone \cite{MR1461821}, the higher dimensional
transition system $u||v$ modeling the concurrent execution of the two
actions $u$ and $v$, depicted by Figure~\ref{concab}, consists of the
transitions $(\alpha,\{u\},\beta)$, $(\beta,\{v\},\delta)$,
$(\alpha,\{v\},\gamma)$, $(\gamma,\{u\},\delta)$ and
$(\alpha,\{u,v\},\delta)$, where the middle term is a multiset, not a
set. The labelling map is in this case the identity map.

This notion is reformulated in \cite{hdts} to make easier a
categorical and homotopical treatment. A higher dimensional system
consists of a set of states $S$, a set of actions $L$ together with a
labelling map $\mu:L\to \Sigma$ where $\Sigma$ is a set of labels, and
a set of tuples of $\bigcup_{n\geq 1} S \p L^n \p S$ satisfying at
least the following two axioms, to obtain the ``minimal'' notion of
\emph{weak transition system}:
\begin{itemize}
\item \textbf{Multiset axiom.} For every permutation $\sigma$ of
  $\{1,\dots,n\}$ with $n\geq 2$, if the tuple
  $(\alpha,u_1,\dots,u_n,\beta)$ is a transition, then the tuple
  $(\alpha,u_{\sigma(1)}, \dots, u_{\sigma(n)}, \beta)$ is a
  transition as well.
\item \textbf{Patching axiom.} For every $(n+2)$-tuple
  $(\alpha,u_1,\dots,u_n,\beta)$ with $n\geq 3$, for every $p,q\geq 1$
  with $p+q<n$, if the five tuples $(\alpha,u_1, \dots, u_n, \beta)$,
  $(\alpha,u_1, \dots, u_p, \nu_1)$, $(\nu_1, u_{p+1}, \dots, u_n,
  \beta)$, $(\alpha, u_1, \dots, u_{p+q}, \nu_2)$ and $(\nu_2,
  u_{p+q+1}, \dots, u_n, \beta)$ are transitions, then the
  $(q+2)$-tuple $(\nu_1, u_{p+1}, \dots, u_{p+q}, \nu_2)$ is a
  transition as well.
\end{itemize}
The multiset axiom avoids the use of multisets. The patching axiom
enables us to see, amongst other things, the $n$-cube as a free object
generated by a $n$-transition. The patching axiom looks like a $5$-ary
composition because it generates a new transition (the patch) from
five transitions satisfying a particular condition. These two axioms
are mathematically designed so that the forgetful functor forgetting
the set of transitions is topological \cite[Theorem~3.4]{hdts}. This
topological structure turns out to be a very powerful tool to deal
with these objects.  Figure~\ref{concab} has in this new formulation
the transitions $(\alpha,u,\beta)$, $(\beta,v,\delta)$,
$(\alpha,v,\gamma)$, $(\gamma,u,\delta)$, $(\alpha,u,v,\delta)$ and
$(\alpha,v,u,\delta)$.

This paper is the direct continuation of \cite{biscsts1}. However, its
reading is not required to read this one. In \cite{biscsts1}, we prove
the existence of a left determined Olschok model structure of weak
transition systems which restricts to left determined Olschok model
structures on various full subcategories without using the map
$R:\{0,1\}\to \{0\}$ in the set of generating cofibrations (unlike
what is done in \cite{cubicalhdts} and \cite{csts}). Then we prove
that the behavior of these homotopy theories break the causal
structure and that the solution to overcome this problem is to work
with \emph{star-shaped transition systems}. A star-shaped transition
system is by definition a \emph{pointed transition system} $(X,*)$,
that means a transition system $X$ together with a distinguished state
$*$ called the \emph{base state}, such that all other states of $X$
are reachable from $*$ by a \emph{path}, i.e. a finite sequence of
$1$-transitions. This path may also be called a \emph{past} of the
state.

In this paper, we want to make precise the above observations which
are only sketched in \cite{biscsts1}.  We work in a reflective
subcategory of all subcategories of higher dimensional transition
systems introduced in \cite{biscsts1}. It is called the subcategory
$\csts$ of \emph{Cattani-Sassone transition systems}. The axiom we add
is CSA1: \emph{if $\alpha \stackrel{u}\to \beta$ and
  $\alpha \stackrel{v}\to \beta$ are two $1$-transitions with
  $\mu(u)=\mu(v)$, then $u=v$}. This axiom is introduced in
\cite{MR1461821} in a more general formulation: we only use its
$1$-dimensional version. It appears also in
\cite[Proposition~72]{MR1365754} (it is called the condition-extension
condition) and in the notion of extensional asynchronous transition
system \cite[page~140]{MR1365754}. It is used in \cite{cubicalhdts}
for a different purpose. The main feature of this axiom is to simplify
the calculations of the cylinder and path functors
(cf. Table~\ref{csts-cyl-path}, Table~\ref{csts-pointed-cyl-path} and
Table~\ref{csts-star-shaped-cyl-path}) while keeping all examples
coming from process algebras \cite{ccsprecub} \cite{symcub}
\cite{hdts}.  The structure of the left determined model category
$\csts$ obtained by restricting the construction of \cite{biscsts1} is
unravelled in the following theorem:

\begin{customthm}{\ref{description-csts}} The left determined Olschok
  model category $\csts$ is Quillen equivalent to the full subcategory
  of Cattani-Sassone transition systems having at most one state
  equipped with the discrete model structure.
\end{customthm}

This theorem means that the homotopy category of Cattani-Sassone
transition systems destroys the causal structure in a very spectacular
way. Thus, localizing or colocalizing this model category will never
give anything interesting from a computer-scientific point of view
because it already contains too many weak equivalences.

The formalism of Cattani-Sassone transition systems is interesting
because, unlike any formalism of labelled precubical sets
\cite{ccsprecub} \cite{symcub} \cite{homotopyprecubical}, it only
contains objects satisfying the higher dimensional automata paradigm
\cite[Definition~7.1]{hdts}: mathematically, this paradigm states that
the boundary of a labelled $n$-cube, with $n\geq 2$, can be filled by
at most one $n$-cube; from a computer scientific point of view, this
paradigm means that the concurrent execution of $n$ actions is
modelled by exactly one $n$-cube, and not two or more having the same
$(n-1)$-dimensional boundary.  The drawback of the formalism of higher
dimensional transition system is that colimits are difficult to
compute because of the patching axiom which freely adds new
transitions in the colimits (cf. \cite[Proposition~A.1]{csts}). Thanks
to the axiom CSA2 (recalled in Section~\ref{reminder} and which plays
the role of face operators in this setting), the category of
Cattani-Sassone transition systems is better behaved with respect to
colimits and the following fact can be considered as an important
result of the paper:

\begin{customthm}{\ref{colim-csts}} The set of transitions of a
  colimit of Cattani-Sassone transition systems is the union of the
  sets of transitions of the components.
\end{customthm}

The calculations made to prove Theorem~\ref{description-csts} enable
us to study the model category $\csts_\bullet$ of star-shaped
(Cattani-Sassone) transition systems. The notion of
\emph{past-similarity} plays a key role in this study. Two states
$\alpha$ and $\beta$ of a star-shaped (Cattani-Sassone) transition
system $(X,*)$ are \emph{past-similar} if there exist two paths from
the base state $*$ to $\alpha$ and $\beta$ respectively which are
homotopic. A star-shaped transition system is \emph{reduced} if two
states are past-similar if and only if they are equal:
Figure~\ref{past-similar} page~\pageref{past-similar} gives an example
of a non-reduced transition system. Using these two new notions
(past-similarity and reduced transition system), the structure of the
model category of star-shaped transition systems is unravelled:

\begin{customthm}{\ref{fib-carac}}
  The fibrant objects of $\csts_\bullet$ are the Cattani-Sassone
  transition systems such that the set of transitions is closed under
  past-similarity. In particular, all reduced transition systems are
  fibrant.
\end{customthm}

\begin{customthm}{\ref{carac-bullet}}
  The left determined Olschok model category $\csts_\bullet$ is
  Quillen equivalent to the full subcategory of reduced objects
  equipped with the discrete model structure. In particular, a map of
  star-shaped Cattani-Sassone transition system is a weak equivalence
  if and only if it is an isomorphism after the identification of all
  past-similar states.
\end{customthm}

The interpretation of Theorem~\ref{carac-bullet} is postponed to the
discussion of Section~\ref{futur} which speculates about the link between 
causality and homotopy in this setting.

\subsection{Outline of the paper}

The paper is structured as follows. Section~\ref{reminder} is a
reminder about weak, cubical, regular and Cattani-Sassone transition
systems. It avoids the reader to have to read the previous papers of
this series. All these notions are necessary because calculations of
limits and colimits often require to start from the topological
structure of weak transition systems, and then to restrict to the
coreflective subcategory of cubical transition systems, and then to
restrict twice to the reflective subcategories of regular and
Cattani-Sassone transition systems. It also happens that some proofs
can be written only by working with cubical transition systems
(e.g. the proof of Proposition~\ref{onto-by-step}).
Section~\ref{section-colim} is a technical section about the
calculation of colimits in the category of regular transition systems.
Theorem~\ref{colim-rts} must be considered as a vast generalization of
\cite[Theorem~4.7]{hdts} and \cite[Proposition~A.3]{csts}. Roughly
speaking, it says that the set of transitions of a colimit of regular
transition systems is the union of the set of transitions.
Section~\ref{def-csts} expounds some basic properties of
Cattani-Sassone transition systems. It also extends
Theorem~\ref{colim-rts} to this new setting: in the category of
Cattani-Sassone transition systems as well, the set of transitions of
a colimit is also the union of the set of transitions.
Section~\ref{modelcsts} uses the toolbox \cite{leftdet} and the
results of \cite{biscsts1} to construct the left determined Olschok
model structure of Cattani-Sassone transition systems. The cylinder
functor is calculated in detail.  Section~\ref{section-path-calcul}
gives a very explicit formulation of the path functor of the model
category constructed in Section~\ref{modelcsts}.
Section~\ref{description-csts-section} proves the first main result of
the paper (Theorem~\ref{description-csts}).  The formulation is chosen
to highlight the destruction of the causal structure.  The notions of
pointed and star-shaped transition systems are recalled in
Section~\ref{model-bullet}. Then the toolbox \cite{leftdet} is used to
prove the existence of the left determined model structures on pointed
and star-shaped Cattani-Sassone transition systems. Meanwhile, we give
precise formulations of the cylinder and path functors of a
star-shaped transition system.  The notion of past-similar states is
introduced and succinctly studied in
Section~\ref{past-similar-sec}. Section~\ref{fibrant-section}
characterizes the fibrant object of the left determined model
structure of star-shaped transition systems (Theorem~\ref{fib-carac}).
Section~\ref{desc-bullet} introduces a particular case of fibrant
objects: the reduced star-shaped transition systems. By definition, a
star-shaped transition system is reduced if past-similarity and
equality coincide. Then the last part of the second main result of the
paper (Theorem~\ref{carac-bullet}) is established. Section~\ref{futur}
is a discussion about an interpretation of Theorem~\ref{carac-bullet}
and about possible future works. In particular, Theorem~\ref{big-pb}
rules out a lot of candidates of model categories. The appendix is an
erratum of the paper \cite{biscsts1}.

\subsection{Prerequisites and notations}

All categories are locally small. The set of maps in a category $\K$
from $X$ to $Y$ is denoted by $\K(X,Y)$. The cardinal of a set $S$ is
denoted by $\# S$.  The class of morphisms of a category $\K$ is
denoted by $\Mor(\K)$. The composite of two maps is denoted by $fg$
instead of $f \circ g$. The initial (final resp.)  object, if it
exists, is always denoted by $\varnothing$ ($\mathbf{1}$ resp.). The
identity of an object $X$ is denoted by $\id_X$.  A subcategory is
always isomorphism-closed (i.e. replete).  A reflective or
coreflective subcategory is always full. By convention,
$A \p B \sqcup C \p D$ means $(A \p B) \sqcup (C \p D)$ where $\p$
denotes the binary product and $\sqcup$ the binary coproduct.  Let $f$
and $g$ be two maps of a locally presentable category $\K$. Denote by
$f\square g$ when $f$ satisfies the \emph{left lifting property} (LLP)
with respect to $g$, or equivalently when $g$ satisfies the
\emph{right lifting property} (RLP) with respect to $f$. Let us
introduce the notations
$\inj_\K(\C) = \{g \in \K, \forall f \in \C, f\square g\}$ called the
class of \emph{$\C$-injective} maps and
$\cof_\K(\C) = \{f \in \K, \forall g\in \inj_\K(\C), f\square g\}$
where $\C$ is a class of maps of $\K$. The class of morphisms of $\K$
that are transfinite compositions of pushouts of elements of $\C$ is
denoted by $\cell_\K(\C)$. There is the inclusion
$\cell_\K(\C)\subset \cof_\K(\C)$. Moreover, every morphism of
$\cof_\K(\C)$ is a retract of a morphism of $\cell_\K(\C)$ as soon as
the domains of $K$ are small relative to $\cell_\K(\C)$
\cite[Corollary~2.1.15]{MR99h:55031}, e.g. when $\K$ is locally
presentable. For every map $f:X \to Y$ and every natural
transformation $\alpha : F \to F'$ between two endofunctors of $\K$,
the map $f\star \alpha$ is defined by the diagram:
\[
\xymatrix@C=1em@R=1em{
FX \fD{\alpha_X}\fR{Ff} && FY \fD{}\ar@/^15pt/@{->}[dddr]_-{\alpha_Y} &\\
&& &&\\
F'X \ar@/_15pt/@{->}[rrrd]^-{F'f}\fR{} &&  \bullet \cocartesien\ar@{->}[rd]^-{f\star \alpha} & \\
&& & F'Y.
}
\]
For a set of morphisms $\mathcal{A}$, let
$\mathcal{A} \star \alpha = \{f\star \alpha, f\in \mathcal{A}\}$. A
cylinder functor $C:\K \to \K$ is equipped with two natural
transformations
$\gamma = \gamma^0 \sqcup \gamma^1:\id \sqcup \id \Rightarrow C$ and
$\sigma:C\Rightarrow \id$ such that
$\sigma\gamma:\id\sqcup \id \Rightarrow \id$ is the codiagonal. A path
functor $P:\K\to \K$ is equipped with two natural transformations
$\tau:\id \Rightarrow P$ and
$\pi=(\pi^0,\pi^1):P \Rightarrow \id \p \id$ such that
$\pi\tau:\id \Rightarrow \id\p \id$ is the diagonal. Sometimes, a
Greek letter denotes a state: the context always enables the reader to
avoid any confusion.  In a model category $\mathcal{M}$, the homotopy
class (left homotopy class, right homotopy class) of maps from $X$ to
$Y$ is denoted by $\pi_\mathcal{M}(X,Y)$ ($\pi^l_\mathcal{M}(X,Y)$,
$\pi^r_\mathcal{M}(X,Y)$ resp.). A cofibrant replacement functor will
be denoted by $(-)^{cof}$ and a fibrant replacement functor by
$(-)^{fib}$. The \emph{discrete model structure} is the model
structure such that all maps are cofibrations and fibrations and such
that the weak equivalences are the isomorphisms \cite{discrete}. A map
of model categories, i.e. a left Quillen functor
$L:\mathcal{M} \to \mathcal{N}$ is \emph{homotopically surjective}
\cite[Definition~3.1]{MR1870516} if for every fibrant object $Y$ of
$\mathcal{N}$ and every cofibrant replacement
$(RY)^{cof} \stackrel{\sim}\to RY$, where $R$ is a right adjoint of
$L$, the induced map $L((RY)^{cof}) \to Y$ is a weak equivalence of
$\mathcal{N}$. A homotopically surjective map of model categories
$L:\mathcal{M} \to \mathcal{N}$ is a Quillen equivalence if and only
if for every cofibrant object $X$ of $\mathcal{M}$ and every fibrant
replacement $LX \stackrel{\sim}\to (LX)^{fib}$, the map
$X \to R((LX)^{fib})$ is a weak equivalence of $\mathcal{M}$.

We refer to \cite{MR95j:18001} for locally presentable categories, to
\cite{MR2506258} for combinatorial model categories, and to
\cite{topologicalcat} for topological categories, i.e. categories
equipped with a topological functor towards a power of the category of
sets.  We refer to \cite{MR99h:55031} and to \cite{ref_model2} for
model categories. For general facts about weak factorization systems,
see also \cite{ideeloc}. The reading of the first part of
\cite{MOPHD}, published in \cite{MO}, is recommended for any reference
about good, cartesian, and very good cylinders.  We use the paper
\cite{leftdet} as a toolbox for constructing the model structures. To
keep this paper short, we refer to \cite{leftdet} for all notions
related to Olschok model categories.

\section{Higher dimensional transition systems}
\label{reminder}

This section is a reminder about weak, cubical, regular and
Cattani-Sassone transition systems.

An infinite set of \emph{labels} $\Sigma$ is fixed. A \emph{transition
  presystem} consists of a triple $X=(\St(X),\mu:\Ac(X)\rightarrow
\Sigma,\Tr(X)=\bigcup_{n\geq 1}\Tr_n(X))$ where $\St(X)$ is a set of
\emph{states}, where $\Ac(X)$ is a set of \emph{actions}, where
$\mu:\Ac(X)\rightarrow \Sigma$ is a set map called the \emph{labelling
  map}, and finally where $\Tr_n(X)\subset \St(X)\p \Ac(X)^n\p \St(X)$
for $n \geq 1$ is a set of \emph{$n$-transitions} or
\emph{$n$-dimensional transitions}. A $n$-transition
$(\alpha,u_1,\dots,u_n,\beta)$ is also called a {\rm transition from
  $\alpha$ to $\beta$}: $\alpha$ is the initial state and $\beta$ the
final state of the transition. It can be denoted by
$\xymatrix@1@C=3em{\alpha \fr{u_1,\dots,u_n} & \beta}$.  This set of
data satisfies one or several of the following axioms (note that the
Intermediate state axiom is a consequence of CSA2):
\begin{itemize}
\item \textbf{Multiset axiom.} For every permutation $\sigma$ of
  $\{1,\dots,n\}$ with $n\geq 2$, if the tuple
  $(\alpha,u_1,\dots,u_n,\beta)$ is a transition, then the tuple
  $(\alpha,u_{\sigma(1)}, \dots, u_{\sigma(n)}, \beta)$ is a
  transition as well.
\item \textbf{Patching axiom.} For every $(n+2)$-tuple
  $(\alpha,u_1,\dots,u_n,\beta)$ with $n\geq 3$, for every $p,q\geq 1$
  with $p+q<n$, if the five tuples $(\alpha,u_1, \dots, u_n, \beta)$,
  $(\alpha,u_1, \dots, u_p, \nu_1)$, $(\nu_1, u_{p+1}, \dots, u_n,
  \beta)$, $(\alpha, u_1, \dots, u_{p+q}, \nu_2)$ and $(\nu_2,
  u_{p+q+1}, \dots, u_n, \beta)$ are transitions, then the
  $(q+2)$-tuple $(\nu_1, u_{p+1}, \dots, u_{p+q}, \nu_2)$ is a
  transition as well.
\item \textbf{All actions are used.} For every $u\in L$, there is a
  $1$-transition $(\alpha,u,\beta)$.
\item \textbf{Intermediate state axiom.} For every $n\geq 2$, every $p$ with
  $1\leq p<n$ and every transition $(\alpha,u_1,\dots,u_n,\beta)$ of
  $X$, there exists a state $\nu$ (not necessarily unique) such that
  both $(\alpha,u_1,\dots,u_p,\nu)$ and
  $(\nu,u_{p+1},\dots,u_n,\beta)$ are transitions.
\item \textbf{CSA2 or Unique intermediate state axiom.} For every $n\geq 2$,
  every $p$ with $1\leq p<n$ and every transition
  $(\alpha,u_1,\dots,u_n,\beta)$ of $X$, there exists a unique state
  $\nu$ such that both $(\alpha,u_1,\dots,u_p,\nu)$ and
  $(\nu,u_{p+1},\dots,u_n,\beta)$ are transitions.
\item \textbf{CSA1.}  If $(\alpha,u,\beta)$ and $(\alpha,v,\beta)$ are
  two transitions such that $\mu(u)=\mu(v)$, then $u=v$.
\end{itemize}
A map of transition presystems consists of two set maps, one between
the sets of states, the other one between the set of actions
preserving the labelling map, such that any transition of the domain
is mapped to a transition of the codomain. For a map $f:X\to Y$ of
transition presystems, the image by $f$ of a transition
$(\alpha,u_1,\dots,u_n,\beta)$ should be noted
$f((\alpha,u_1,\dots,u_n,\beta))$. The notation
$f(\alpha,u_1,\dots,u_n,\beta)$ will be used instead to not overload
the calculations. This convention is already implicitly used in our
previous papers. The mapping $X\to \St(X)$ ( $X\to \Ac(X)$ resp.)
induces a functor from the category of transition presystems to the
category of sets $\set$.

Table~\ref{alldefs} lists the definitions of the categories $\wts$ of
\emph{weak transition systems} \cite[Definition~3.2]{hdts}, $\cts$ of
\emph{cubical transition systems} \cite[Proposition~6.7]{cubicalhdts},
$\rts$ of \emph{regular transition systems}
\cite[Definition~2.2]{csts} and $\csts$ of \emph{Cattani-Sassone
  transition systems} \cite[Table~1]{csts}. All examples coming from
process algebras belong to $\csts$ \cite{ccsprecub} \cite{symcub}
\cite{hdts}.

\begin{table}[ht]
\begin{tabular}{lcccc}
& Cattani-Sassone & Regular & Cubical & Weak  \\
\hline
Multiset axiom & yes & yes & yes & yes \\
Patching axiom & yes & yes & yes & yes \\
All actions used &  yes & yes & yes & no \\
Intermediate state axiom & yes & yes & yes & no\\
Unique intermediate state axiom  & yes & yes & no & no \\
CSA1 & yes & no & no & no \\
\hline
\\
\end{tabular}
\caption{Cattani-Sassone, regular, cubical and weak transition systems.}
\label{alldefs}
\end{table}

The category $\csts$ is a reflective subcategory of $\cts$ by
\cite[Proposition~7.2]{cubicalhdts}, but also of $\rts$ by
Proposition~\ref{csa1-reflective}. We will come back on the category
$\csts$ in Section~\ref{def-csts}. The category $\rts$ is locally
finitely presentable by \cite[Proposition~4.5]{csts}. It is a
reflective subcategory of $\cts$ by \cite[Proposition~4.4]{csts}. The
reflection is the functor $\CSA_2:\cts\to \cts$ which forces CSA2 to
hold. This functor is extensively studied in
\cite[Section~4]{csts}. The category $\cts$ is locally finitely
presentable by \cite[Corollary~3.15]{cubicalhdts}. By
\cite[Corollary~3.15]{cubicalhdts}, the category $\cts$ is a
coreflective subcategory of $\wts$. The latter is locally finitely
presentable by \cite[Theorem~3.4]{hdts}. The forgetful functor
$\omega:\wts \to \set^{\{s\}\cup \Sigma}$ taking the weak higher
dimensional transition system $X$ to the $(\{s\}\cup \Sigma)$-tuple of
sets $(\St(X),(\mu^{-1}(x))_{x\in \Sigma}) \in \set^{\{s\}\cup
  \Sigma}$ is topological by \cite[Theorem~3.4]{hdts}.  There is the
chain of functors
\[\csts \subset_{\hbox{\tiny reflective}} \rts \subset_{\hbox{\tiny reflective}} \cts \subset_{\hbox{\tiny coreflective}} \wts \stackrel{\omega}\longrightarrow_{\hbox{\tiny topological}} \set^{\{s\}\cup \Sigma}.\]
We give now some important examples of regular transition systems.

\begin{enumerate}
\item Every set $X$ may be identified with the cubical transition
  system having the set of states $X$, with no actions and no
  transitions.
\item For every $x\in \Sigma$, let us denote by $\dd{x}$ the cubical
  transition system with four states $\{1,2,3,4\}$, one action $x$ and
  two transitions $(1,x,2)$ and $(3,x,4)$. The cubical transition
  system $\dd{x}$ is called the \emph{double transition (labelled by
    $x$)} where $x\in \Sigma$.
\end{enumerate}

\begin{nota} For $n\geq 1$, let $0_n = (0,\dots,0)$ ($n$-times) and
  $1_n = (1,\dots,1)$ ($n$-times). By convention, let
  $0_0=1_0=()$. \end{nota}

Let us introduce the weak transition system corresponding to the
labelled $n$-cube.

\bp \label{cas_cube} \cite[Proposition~5.2]{hdts} Let $n\geq 0$ and
$x_1,\dots,x_n\in \Sigma$. Let $T_d\subset \{0,1\}^n \p
\{(x_1,1),\dots,(x_n,n)\}^d \p \{0,1\}^n$ (with $d\geq 1$) be the
subset of $(d+2)$-tuples
\[((\epsilon_1,\dots,\epsilon_n), (x_{i_1},i_1),\dots,(x_{i_d},i_d),
(\epsilon'_1,\dots,\epsilon'_n))\] such that
\begin{itemize}
\item $i_m = i_n$ implies $m = n$, i.e. there are no repetitions in the
  list $(x_{i_1},i_1),\dots,(x_{i_d},i_d)$
\item for all $i$, $\epsilon_i\leq \epsilon'_i$
\item $\epsilon_i\neq \epsilon'_i$ if and only if
  $i\in\{i_1,\dots,i_d\}$. 
\end{itemize}
Let $\mu : \{(x_1,1),\dots,(x_n,n)\} \rightarrow \Sigma$ be the set
map defined by $\mu(x_i,i) = x_i$. Then \[C_n[x_1,\dots,x_n] =
(\{0,1\}^n,\mu : \{(x_1,1),\dots,(x_n,n)\}\rightarrow
\Sigma,(T_d)_{d\geq 1})\] is a well-defined weak transition system
called the {\rm $n$-cube}. \ep

The $n$-cubes $C_n[x_1,\dots,x_n]$ for all $n\geq 0$ and all
$x_1,\dots,x_n\in \Sigma$ are regular by \cite[Proposition~4.6]{hdts}
and \cite[Proposition~5.2]{hdts}.  For $n = 0$, $C_0[]$, also denoted
by $C_0$, is nothing else but the set $\{()\}$.  

Here are two important families of weak transition systems which are
not cubical, and therefore not regular:
\begin{enumerate}
\item The weak transition system $\underline{x} = (\varnothing, \{x\}
  \subset \Sigma, \varnothing)$ for $x \in \Sigma$ is not cubical
  because the action $x$ is not used. 
\item Let $n \geq 0$. Let $x_1,\dots,x_n \in \Sigma$. The \emph{pure
    $n$-transition} $C^{ext}_n[x_1,\dots,x_n]$ is the weak transition
  system with the set of states $\{0_n,1_n\}$, with the set of actions
  \[\{(x_1,1), \dots, (x_n,n)\}\] and with the transitions all
  $(n+2)$-tuples $(0_n,(x_{\sigma(1)},\sigma(1)), \dots,
  (x_{\sigma(n)},\sigma(n)),1_n)$ for $\sigma$ running over the set of
  permutations of the set $\{1,\dots ,n\}$. It is not cubical for $n>
  1$ because it does not contain any $1$-transition. Intuitively, the
  pure transition is a cube without faces of lower dimension.
\end{enumerate}
The main use of the family of pure transitions is summarized in the
following two facts:
\begin{enumerate}
\item For all weak transition systems $X$, the set
  $\wts(C^{ext}_n[x_1,\dots,x_n],X)$ is the set of transitions
  $(\alpha,u_1,\dots,u_n,\beta)$ of $X$ such that for all
  $1\leq i\leq n$, $\mu(u_i)=x_i$ and
  \[\bigsqcup_{x_1,\dots,x_n\in \Sigma}\wts(C^{ext}_n[x_1,\dots,x_n],X)\]
  is the set of transitions of $X$.
\item Every map of weak transition systems
$f:C^{ext}_n[x_1,\dots,x_n] \to X$ where $X$ satisfies CSA2 factors
uniquely as a composite $f:C^{ext}_n[x_1,\dots,x_n] \to
C_n[x_1,\dots,x_n] \to X$ by \cite[Theorem~5.6]{hdts}.
\end{enumerate}

We conclude this section by recalling some important facts:

\bp \label{121-121-trans} Let $f:A\to B$ be a map of weak transition
systems which is one-to-one on states and on actions. Then it is
one-to-one on transitions.  \ep

\bpf It is mutatis mutandis the proof of
\cite[Proposition~.4.4]{homotopyprecubical}. \epf

\bth \cite[Theorem~3.3]{csts} \label{cubical-lift} Let $(f_i:\omega
(A_i)\to W)_{i\in I}$ be a cocone of $\set^{\{s\}\cup \Sigma}$ such
that the weak transition systems $A_i$ are cubical for all $i\in I$
and such that every action $u$ of $W$ is the image of an action of
$A_{i_u}$ for some $i_u\in I$. Then the $\omega$-final lift
$\overline{W}$ is cubical.  \eth

\bp \cite[Proposition~4.1]{csts} \label{mor-cub-reg} Let $X$ be a
cubical transition system. Let $Y$ be a weak transition system
satisfying CSA2. Let $f:X\to Y$ be a map of weak transition systems
which is one-to-one on states.  Then $X$ is regular, and in particular
satisfies CSA2.  \ep

\section{Colimit of regular transition systems}
\label{section-colim}

Theorem~\ref{colim-rts} states that the set of transitions of a
colimit of regular transition systems is the union of the set of
transitions. Its proof is similar to the proofs of
\cite[Theorem~4.7]{hdts} and \cite[Proposition~A.3]{csts}. 

We need to recall two lemmas about colimits of weak transition systems
and cubical transition systems coming from \cite{biscsts1}.

\begin{lem} \label{colim-preserv0} The forgetful functor mapping a
  weak transition system to its set of states is colimit-preserving.
  The forgetful functor mapping a weak transition system to its set of
  actions is colimit-preserving.
\end{lem}

\bpf The lemma is a consequence of the fact that the functor
$\omega : \wts \longrightarrow \set^{\{s\}\cup \Sigma}$ taking the
weak higher dimensional transition system
$(S,\mu : L \rightarrow \Sigma,(T_n)_{n\geq 1})$ to the
$(\{s\}\cup \Sigma)$-tuple of sets
$(S,(\mu^{-1}(x))_{x\in \Sigma}) \in \set^{\{s\}\cup \Sigma}$ is
topological.  \epf

\begin{lem} \label{colim-preserv} The forgetful functor mapping a
  cubical transition system to its set of states is
  colimit-preserving.  The forgetful functor mapping a cubical
  transition system to its set of actions is colimit-preserving.
\end{lem}

\bpf Since the category of cubical transition systems is a
coreflective subcategory of the category of weak transition systems by
\cite[Corollary~3.15]{cubicalhdts}, this lemma is a consequence of
Lemma~\ref{colim-preserv0}.  \epf

\bth \label{colim-rts} Let $(i\mapsto X_i)$ be a small diagram of
$\rts$. The set of states $\St(\liminj X_i)$ is a quotient of the set
$\liminj \St(X_i)$, the set of actions $\Ac(\liminj X_i)$ is equal to
$\liminj \Ac(X_i)$ and the set of transitions of $\liminj X_i$ is
equal to $\bigcup_i \phi_i(\Tr(X_i))$ where $\phi_i:X_i \to \liminj
X_i$ is the canonical map. In particular, the regular transition
system $\liminj X_i$ is equipped with the $\omega$-final structure.
\eth

\bpf The proof is divided in two parts. The first one is easy. The
second one requires to be more careful and shows how the patching
axiom can be used.

\underline{The case of states and actions}. The colimit $\liminj X_i$
in $\rts$ is equal to $\CSA_2(\liminj\nolimits^\cts X_i)$ where
$\liminj^\cts X_i$ is the colimit calculated in $\cts$. Since the
functors $\St:\cts \to \set$ and $\Ac:\cts \to \set$ are
colimit-preserving by Lemma~\ref{colim-preserv}, we have the bijection
of sets $\liminj \St(X_i) \iso \St(\liminj\nolimits^\cts X_i)$ and
$\liminj \Ac(X_i) \iso \Ac(\liminj\nolimits^\cts X_i)$. The unit map
$\liminj\nolimits^\cts X_i \to \CSA_2(\liminj\nolimits^\cts X_i)$ is
onto on states and bijective on actions by
\cite[Proposition~4.2]{csts}: the functor $\CSA_2:\cts\to \rts$ forces
CSA2 to hold by making identifications of states.  We deduce that the
set of states of $\liminj X_i$ is a quotient of $\liminj \St(X_i)$ and
that the set of actions of $\liminj X_i$ is exactly
$\liminj \Ac(X_i)$.

\underline{The case of transitions}. Let $T=\bigcup_i
\phi_i(\Tr(X_i))$. Let
$(f_i(\alpha),f_i(u_1),\dots,f_i(u_n),f_i(\beta))$ be a tuple of $T$.
Then the tuple $(\alpha,u_{\sigma(1)},\dots,u_{\sigma(n)},\beta)$ is a
transition of $X_i$ for all permutations $\sigma$ of
$\{1,\dots,n\}$. So the tuple
$(f(\alpha),f(u_{\sigma(1)}),\dots,f(u_{\sigma(n)}),f(\beta))$ belongs
to $T$. This means that the set of tuples $T$ satisfies the multiset
axiom. Let $n\geq 3$ and $p,q \geq 1$ with $p+q < n$. Let
  \begin{multline*}(\alpha,u_1, \dots, u_n, \beta),(\alpha,u_1, \dots,
u_p, \mu),(\mu, u_{p+1}, \dots, u_n, \beta),\\(\alpha, u_1, \dots,
u_{p+q}, \nu),(\nu, u_{p+q+1}, \dots, u_n, \beta)
\end{multline*} be five tuples of $T$. Let $(\alpha,u_1, \dots, u_n,
\beta) = (f_i(\gamma),f_i(v_1),\dots,f_i(v_n),f_i(\delta))$ where the
tuple $(\gamma,v_1,\dots,v_n,\delta)$ is a transition of $X_i$.  There
exist two states $\epsilon$ and $\eta$ of $X_i$ such that the five
tuples\begin{multline*} (\gamma,v_1,\dots,v_p,\epsilon),
  (\gamma,v_1,\dots,v_{p+q},\eta),\\
  (\epsilon,v_{p+1},\dots,v_n,\delta),
  (\eta,v_{p+q+1},\dots,v_n,\delta),
  (\epsilon,v_{p+1},\dots,v_{p+q},\eta) \end{multline*} are
transitions of $X_i$ since $X_i$ is cubical and by using the patching
axiom in $X_i$.  Therefore, the five tuples
\begin{multline*}
  (\alpha,u_1,\dots,u_p,f_i(\epsilon)),
  (\alpha,u_1,\dots,u_{p+q},f_i(\eta)),\\
  (f_i(\epsilon),u_{p+1},\dots,u_n,\beta),
  (f_i(\eta),u_{p+q+1},\dots,u_n,\beta),\\
  (f_i(\epsilon),u_{p+1},\dots,u_{p+q},f_i(\eta))
\end{multline*}
are transitions of $\liminj X_i$. The point is that $\liminj X_i$
satisfies CSA2. We deduce that $f_i(\epsilon)=\mu$ and
$f_i(\eta)=\nu$.  We obtain that
\[(\mu,u_{p+1},\dots,u_{p+q},\nu)=(f_i(\epsilon),f_i(v_{p+1}),\dots,f_i(v_{p+q}),f_i(\eta))
\in T.\] This means that the set of tuples $T$ satisfies the patching
axiom. Write $Z$ for the weak transition system having the set of
transitions $T$.  We have obtained a morphism of cocones of \emph{weak
  transition systems}
\[
\xymatrix@C=3em@R=3em
{
(X_i) \ar@{=}[r] \fd{} & (X_i) \fd{} \\
Z \fr{\subset} & \liminj X_i.
}
\]
The map $Z\to \liminj X_i$ is bijective on states and on actions, and
therefore one-to-one on transitions by
Proposition~\ref{121-121-trans}. We have also proved that
the weak transition system $Z$ is the $\omega$-final lift of the
cocone $(\omega(X_i) \to \omega(Z))$ of $\set^{\{s\}\cup \Sigma}$.  By
Theorem~\ref{cubical-lift}, the weak transition system $Z$ is cubical.
Since the map $Z\to \liminj X_i$ is one-to-one on states, the cubical
transition system $Z$ satisfies CSA2 by Proposition~\ref{mor-cub-reg}.
We have proved that $Z= \liminj X_i$.  \epf

\begin{cor} \label{onto-trans-rts} Consider a pushout diagram of $\rts$
\[
\xymatrix@C=3em@R=3em
{
A \fr{\phi} \fd{p} & X \fd{f} \\
B \fr{\psi} & \cocartesien Y
}
\]
such that the left vertical map $A\to B$ is onto on states, on actions
and on transitions. Then the right vertical map is onto on states, on
actions and on transitions.
\end{cor}

\bpf Since the map $A\to B$ is onto on states and on actions, the map
$X\to Y$ is onto on states and on actions as well by
Theorem~\ref{colim-rts}.  By Theorem~\ref{colim-rts}, we have the
equality $\Tr(Y)= f(\Tr(X)) \cup \psi(\Tr(B))$. Since $A\to B$ is onto
on transitions, we also have the equalities $\psi(\Tr(B)) =
\psi(p(\Tr(A)) = f(\phi(\Tr(A))$. We deduce the inclusion
$\psi(\Tr(B)) \subset f(\Tr(X))$. We obtain $\Tr(Y)= f(\Tr(X))$.  \epf

In Corollary~\ref{onto-trans-rts}, $\rts$ cannot be replaced by $\cts$
(cf. \cite[Proposition~A.1]{csts}). This ``good'' behavior of colimits
in $\rts$ is due to CSA2.

\section{Cattani-Sassone transition system}
\label{def-csts}

All sets, all double transitions and all cubes are Cattani-Sassone
transition systems. All weak transition systems coming from a process
algebra are Cattani-Sassone transition systems \cite{ccsprecub}
\cite{symcub}.

\begin{nota}
  Let 
$\mathcal{U} =
  \{C_1[x]\sqcup_{\{0_1,1_1\}} C_1[x] \to C_1[x] \mid x \in \Sigma\}$.
\end{nota}

\bp \label{CSA1-carac} Let $X$ be a weak transition system. The
following conditions are equivalent:
\begin{enumerate}
\item $X$ is $\mathcal{U}$-injective.
\item $X$ is $\mathcal{U}$-orthogonal.
\item $X$ satisfies CSA1.
\end{enumerate}
\ep

\bpf Every map of $\mathcal{U}$ is bijective on states and onto
on actions.  Thus, every map of $\mathcal{U}$ is epic. We deduce
the equivalence $(1) \Leftrightarrow (2)$. The equivalence 
$(1) \Leftrightarrow (3)$ is clear. \epf

\bp \label{epic_map} Every map of $\cell_\rts(\mathcal{U})$ is onto on
states, on actions and on transitions. In particular, every map of
$\cell_\rts(\mathcal{U})$ is epic.  \ep

\bpf Every pushout in $\rts$ of a map of $\mathcal{U}$ is onto on
states, on actions and on transitions by
Corollary~\ref{onto-trans-rts}. Every transfinite composition in
$\rts$ of maps onto on states, on actions and on transitions is onto
on states, on actions and on transitions by
Theorem~\ref{colim-rts}. Therefore, every map of
$\cell_\rts(\mathcal{U})$ is onto on states, on actions and on
transitions. In particular, every map of $\cell_\rts(\mathcal{U})$ is
epic.  \epf

Since every map of $\cell_\rts(\mathcal{U})$ is epic by
Proposition~\ref{epic_map}, the canonical map $X\to \mathbf{1}$ from
$X$ to the terminal object of $\rts$ factors functorially and
uniquely, up to isomorphism, as a composite
\[\xymatrix@C=4em{X \fr{\in \cell_\rts(\mathcal{U})}& \CSA_1(X)
  \fr{\in \inj_\rts(\mathcal{U})}& \mathbf{1}}\]
in $\rts$ where the left-hand map belongs to $\cell_\rts(\mathcal{U})$
and the right-hand map belongs to $\inj_\rts(\mathcal{U})$ (see
\cite[Proposition~A.1]{biscsts1} for a proof of uniqueness). By
Proposition~\ref{CSA1-carac}, the regular transition system
$\CSA_1(X)$ satisfies CSA1. We have obtained a well-defined functor
$\CSA_1:\rts \to \rts$.

\bp \label{csa1-reflective} The category $\csts$ is a reflective
locally finitely presentable subcategory of $\rts$. The left adjoint
of the inclusion $\csts \subset \rts$ is precisely the functor
$\CSA_1:\rts \to \rts$. \ep

\bpf By Proposition~\ref{CSA1-carac}, $\csts$ is a small-orthogonality
class of $\rts$. Hence it is a reflective subcategory. Thus, $\csts$
is complete and cocomplete.  The set of cubes and double transitions
is a dense and hence strong generator of finitely presentable objects
of $\rts$ by \cite[Theorem~3.11]{cubicalhdts}, and therefore of
$\csts$. Consequently, the category $\csts$ is locally finitely
presentable by \cite[Theorem~1.20]{MR95j:18001}. Let $f:X\to Z$ be a
map of regular transition systems such that $Z$ satisfies CSA1. Then
we have the commutative diagram of $\rts$
\[
\xymatrix@C=3em@R=3em
{
X \fr{} \fd{} & Z \fd{} \\
\CSA_1(X) \ar@{-->}[ru]^-{\ell} \fr{} & \mathbf{1}.
}
\]
By construction, the left vertical map belongs to
$\cell_\rts(\mathcal{U})$. Since $Z$ satisfies CSA1, the right
vertical map belongs to $\inj_\rts(\mathcal{U})$ by
Proposition~\ref{CSA1-carac}. Therefore the lift $\ell$ exists and it
is unique since the left vertical map is epic by
Proposition~\ref{epic_map}.  Thus, the map $X\to Z$ factors uniquely
as a composite $X\to \CSA_1(X)\to Z$.  \epf

\bth \label{colim-csts} Let $(i\mapsto X_i)$ be a small diagram of
$\csts$. The set of states $\St(\liminj X_i)$ is a quotient of the set
$\liminj \St(X_i)$, the set of actions $\Ac(\liminj X_i)$ is a
quotient of the set $\liminj \Ac(X_i)$ and the set of transitions of
$\liminj X_i$ is equal to $\bigcup_i \phi_i(\Tr(X_i))$ where
$\phi_i:X_i \to \liminj X_i$ is the canonical map. In particular, the
Cattani-Sassone transition system $\liminj X_i$ is equipped with the
$\omega$-final structure.  \eth

\bpf We have $\liminj X_i \iso \CSA_1(\liminj^\rts X_i)$ where
$\liminj^\rts$ is the colimit calculated in $\rts$. The unit map
$\liminj^\rts X_i \to \CSA_1(\liminj^\rts X_i)$ belongs to
$\cell_\rts(\mathcal{U})$.  Therefore it is onto on states, on actions
and on transitions by Proposition~\ref{epic_map}. The proof is
complete thanks to Theorem~\ref{colim-rts}.  \epf

\begin{cor} \label{onto-trans-csts} Consider a pushout diagram of $\csts$
\[
\xymatrix@C=3em@R=3em
{
A \fr{\phi} \fd{p} & X \fd{f} \\
B \fr{\psi} & \cocartesien Y
}
\]
such that the left vertical map $A\to B$ is onto on states, on actions
and on transitions. Then the right vertical map is onto on states, on
actions and on transitions. 
\end{cor}

\bpf The proof is mutatis mutandis the proof of
Corollary~\ref{onto-trans-rts}.  \epf

\section{Homotopy theory of non star-shaped objects}
\label{modelcsts}

Let $X$ be a weak transition system. By
\cite[Proposition~2.10]{biscsts1}, there exists a weak transition
system $\cyl(X)$ with the set of states $\St(X)\p \{0,1\}$, the set of
actions $\Ac(X)\p \{0,1\}$, the labelling map the composite map
$\mu:\Ac(X)\p \{0,1\} \to \Ac(X) \to \Sigma$, and such that a tuple
$((\alpha,\epsilon_0),(u_1,\epsilon_1),\dots,(u_n,\epsilon_n),(\beta,\epsilon_{n+1}))$
is a transition of $\cyl(X)$ if and only if the tuple
$(\alpha,u_1,\dots,u_n,\beta)$ is a transition of $X$.

\begin{nota} Let $Z\subset \St(X)$. The $\omega$-final lift of the map
\[\omega(\cyl(X)) \longrightarrow (Z\p \{0\} \sqcup (\St(X)\backslash
Z)\p \{0,1\},\Ac(X) \p \{0,1\})\] of $\set^{\{s\}\cup \Sigma}$ taking
$(s,\epsilon)\in Z\p \{0,1\}$ to $(s,0)$ and which is the identity on
$(\St(X)\backslash Z)\p \{0,1\}$ and on $\Ac(X)\p \{0,1\}$ is denoted
by $\cyl(X)//Z$. \end{nota}

The set of transitions of $\cyl(X)//Z$ is the set of tuples
$((\alpha,\epsilon_0),(u_1,\epsilon_1),\dots,(u_n,\epsilon_n),\linebreak[4](\beta,\epsilon_{n+1}))$
such that $(\alpha,u_1,\dots,u_n,\beta)$ is a transition of
$X$. Indeed, the $\omega$-final structure contains this set of tuples
and this set of tuples obviously satisfies the multiset axiom and the
patching axiom; consequently, it is the $\omega$-final structure.

By \cite[Theorem~3.3]{biscsts1}, if the weak transition system $X$ is
cubical, then the weak transition system $\cyl(X)$ introduced above is
cubical as well.  The map
\[\omega(\cyl(X)) \longrightarrow (Z\p \{0\} \sqcup (\St(X)\backslash
Z)\p \{0,1\},\Ac(X) \p \{0,1\})\]
of $\set^{\{s\}\cup \Sigma}$ induces the identity of
$\Ac(X)\p \{0,1\}$ on actions, and thus, is onto on actions. By
Theorem~\ref{cubical-lift}, we deduce that if the weak transition
system $X$ is cubical, then the weak transition system $\cyl(X)//Z$ is
cubical as well.

\begin{nota} Let $Z\subset \St(X)$. The $\omega$-final lift of the map
  of $\set^{\{s\}\cup \Sigma}$ \[\omega(\cyl(X)) \longrightarrow (Z\p
  \{0\} \sqcup (\St(X)\backslash Z)\p \{0,1\},\Ac(X))\] taking
  $(s,\epsilon)\in Z\p \{0,1\}$ to $(s,0)$ and taking the action
  $(u,\epsilon) \in \Ac(X)\p \{0,1\}$ to $u$ is denoted by
  $\cyl(X)///Z$. \end{nota}

\bp \label{precalcul0} Let $X$ be a cubical transition system. Let
$Z\subset \St(X)$ be a subset of the set of states of $X$.  The weak
transition system $\cyl(X)///Z$ is cubical. A tuple
$((\alpha,\epsilon_0),u_1,\dots,u_n,(\beta,\epsilon_{n+1}))$ is a
transition of $\cyl(X)///Z$ if and only if the tuple
$(\alpha,u_1,\dots,u_n,\beta)$ is a transition of $X$. Moreover, if
$X$ satisfies CSA1, then $\cyl(X)///Z$ satisfies CSA1 as well. \ep

\bpf By Theorem~\ref{cubical-lift}, the weak transition system
$\cyl(X)///Z$ is cubical since the projection map $\Ac(X)\p \{0,1\}
\to \Ac(X)$ is onto.  The set of transitions of $\cyl(X)///Z$ is the
$\omega$-final structure generated by the set of tuples \[\mathcal{T}
= \{((\alpha,\epsilon_0),u_1,\dots,u_n,(\beta,\epsilon_{n+1})) \mid
(\alpha,u_1,\dots,u_n,\beta) \in \Tr(X) \}.\] The set $T$ of
transitions satisfies the multiset axiom and the patching
axiom. Therefore it is the $\omega$-final structure:
$\mathcal{T}=\Tr(\cyl(X)///Z)$. Assume that $X$ satisfies CSA1. Let
$((\alpha,\epsilon_0),u,(\beta,\epsilon_1))$ and
$((\alpha,\epsilon_0),v,(\beta,\epsilon_1))$ be two transitions of the
cubical transition system $\cyl(X)///Z$ with $\mu(u)=\mu(v)$.  Then
the tuples $(\alpha,u,\beta)$ and $(\alpha,v,\beta)$ are two
transitions of $X$. Since $X$ satisfies CSA1 by hypothesis, we obtain
$u=v$. Consequently, $\cyl(X)///Z$ satisfies CSA1 if $X$ does.  \epf

\bd \cite{biscsts1} Let $X$ be a weak transition system. A state
$\alpha$ of $X$ is {\rm internal} if there exists three transitions
$(\gamma,u_1,\dots,u_n,\delta)$, $(\gamma,u_1,\dots,u_p,\alpha)$ and
$(\alpha,u_{p+1},\dots,u_n,\delta)$ with $n\geq 2$ and $p\geq 1$. A
state $\alpha$ is {\rm external} if it is not internal. \ed

\begin{nota} The set of internal states of a weak transition system
  $X$ is denoted by $\Int(X)$. The complement is denoted by
  $\Ext(X)=\St(X) \backslash \Int(X)$. \end{nota}

For any map $f:X\to Y$ of weak transition systems, we have
$f(\Int(X)) \subset \Int(Y)$: any internal state of $X$ is mapped to
an internal state of $Y$. In general, we have
$f(\Ext(X)) \nsubseteq \Ext(Y)$. In the example of the inclusion
$00\stackrel{u}\longrightarrow 01 \stackrel{v}\longrightarrow 11
\subset C_2[\mu(u),\mu(v)]$,
all states of the domains are external, and the middle state $01$ is
mapped to an internal state of $C_2[\mu(u),\mu(v)]$. We actually have
the following characterization:

\bp \label{int_charac} A state $\alpha$ of a cubical transition
systems $X$ is internal if and only if there exist three transitions
$(\gamma,u_1,u_2,\delta)$, $(\gamma,u_1,\alpha)$ and
$(\alpha,u_2,\delta)$. \ep

\bpf The ``if'' part is a consequence of the definition. The ``only
if'' part is a consequence of the fact that $X$ is cubical.
\epf

\bp \label{precalcul} Let $X$ be a Cattani-Sassone transition
system. Then the cubical transition system $\cyl(X)/// \Int(X)$ is a
Cattani-Sassone transition system.  \ep

\bpf 
Consider the five transitions
\begin{align*}
&((\alpha,\epsilon_0),u_1,\dots,u_n,(\beta,\epsilon_{n+1}))\\
&((\alpha,\epsilon_0),u_1,\dots,u_p,(\gamma,\epsilon)),((\gamma,\epsilon),u_{p+1},\dots,u_n,(\beta,\epsilon_{n+1}))\\
&((\alpha,\epsilon_0),u_1,\dots,u_p,(\gamma',\epsilon')),((\gamma',\epsilon'),u_{p+1},\dots,u_n,(\beta,\epsilon_{n+1}))
\end{align*} 
of $\cyl(X)/// \Int(X)$. Then the five tuples 
\begin{align*}
&(\alpha,u_1,\dots,u_n,\beta)\\
&(\alpha,u_1,\dots,u_p,\gamma),(\gamma,u_{p+1},\dots,u_n,\beta)\\
&(\alpha,u_1,\dots,u_p,\gamma'),(\gamma',u_{p+1},\dots,u_n,\beta)
\end{align*}
are transitions of $X$. Since $X$ is regular, we have
$\gamma=\gamma'$. And since $\gamma=\gamma'$ is an internal state of
$X$, we have $(\gamma,\epsilon)=(\gamma,0)$ and
$(\gamma',\epsilon)=(\gamma',0)$ in $\cyl(X) /// \Int(X)$. We deduce
that $\cyl(X) /// \Int(X)$ is regular. The proof is complete using
Proposition~\ref{precalcul0}.  \epf

\bth \label{precalcul2} Let $X$ be a Cattani-Sassone transition
system. Then there is the natural isomorphism $\CSA_1\CSA_2 \cyl(X)
\iso \cyl(X)/// \Int(X)$.  \eth

\bpf By \cite[Lemma~3.14]{biscsts1}, we have the isomorphism of weak
transition systems \[\CSA_2 \cyl(X) \iso \cyl(X)// \Int(X).\] In
particular, this means that $\cyl(X)// \Int(X)$ is regular. Let $u\in
\Ac(X)$. Since $X$ is cubical, there exists a transition
$(\alpha_u,u,\beta_u)$ of $X$. We deduce that
$((\alpha_u,0),(u,0),(\beta_u,0))$ and
$((\alpha_u,0),(u,1),(\beta_u,0))$ are two transitions of $\cyl(X)$,
and therefore two transitions of $\cyl(X)//\Int(X)$. Consider the map
$\phi_u:C_1[\mu(u)] \sqcup_{0_1,1_1} C_1[\mu(u)] \to \cyl(X)//\Int(X)$
which sends the two transitions of the domain to the transitions
$((\alpha_u,0),(u,0),(\beta_u,0))$ and
$((\alpha_u,0),(u,1),(\beta_u,0))$ respectively.  Consider the pushout
diagram of $\rts$
\[
\xymatrix@C=5em@R=3em
{
\bigsqcup_{u\in \Ac(X)} C_1[\mu(u)] \sqcup_{0_1,1_1} C_1[\mu(u)] \fd{} \fr{\bigsqcup_{u\in \Ac(X)} \phi_u} & \cyl(X)//\Int(X) \fd{} \\
\bigsqcup_{u\in \Ac(X)} C_1[\mu(u)] \fr{} & \cocartesien Z.
}
\]
By Theorem~\ref{colim-rts}, the set of states $\St(Z)$ is a quotient
of the set of states $\St(\cyl(X)//\Int(X))$, the set of actions
$\Ac(Z)$ of $Z$ is equal to $\Ac(X)$ and the set of transitions of $Z$
is given by the $\omega$-final structure.  The map of $\set^{\{s\}\cup
  \Sigma}$ \[\omega(\cyl(X))
\longrightarrow (\Int(X)\p \{0\} \sqcup \Ext(X)\p \{0,1\},\Ac(X))\]  factors as a composite
\begin{multline*}\omega(\cyl(X))\longrightarrow (\Int(X)\p \{0\} \sqcup \Ext(X)\p \{0,1\},\Ac(X)\p \{0,1\})
  \\\longrightarrow (\Int(X)\p \{0\} \sqcup \Ext(X)\p \{0,1\},\Ac(X))\end{multline*} where the left hand map is
the identity on $\Ac(X)\p \{0,1\}$.  Consequently, the map $\cyl(X)
\to \cyl(X)/// \Int(X)$ factors uniquely as a composite \[\cyl(X) \to
\cyl(X)//\Int(X) \to \cyl(X)/// \Int(X).\] The right-hand map
\[\cyl(X)//\Int(X)\to \cyl(X)/// \Int(X)\] is bijective on states and
$\cyl(X)/// \Int(X)$ is obtained from $\cyl(X)//\Int(X)$ by making the
identifications $(u,0)=(u,1)$ for all actions $u$ of
$X$. Consequently, by the universal property of the pushout, the map
$\cyl(X)//\Int(X)\to \cyl(X)/// \Int(X)$ factors uniquely as a
composite \[\cyl(X)//\Int(X)\longrightarrow Z \longrightarrow
\cyl(X)/// \Int(X).\] The latter composite set map yields the
factorization of $\id_{\St(\cyl(X)//\Int(X))}$ as the composite
\[\St(\cyl(X)//\Int(X)) \twoheadrightarrow \St(Z) \longrightarrow \St(\cyl(X)///\Int(X)).\] We
deduce that the left-hand map is one-to-one, and then bijective. We
have obtained the isomorphism $Z\iso \cyl(X)/// \Int(X)$ because the
two Cattani-Sassone transition systems are the $\omega$-final lift of
the same map of $\set^{\{s\}\cup \Sigma}$. We obtain a factorization
of the canonical map $\cyl(X)// \Int(X) \to \mathbf{1}$ as a composite
$\cyl(X)// \Int(X) \to \cyl(X)/// \Int(X)\to \mathbf{1}$ where the
left-hand map belongs to $\cell_\rts(\mathcal{U})$ and where, by
Proposition~\ref{precalcul0} and Proposition~\ref{CSA1-carac}, the
right-hand map belongs to $\inj_\rts(\mathcal{U})$. We deduce the
isomorphism $\CSA_1(\cyl(X)// \Int(X)) \iso \cyl(X)/// \Int(X)$.  \epf

\begin{nota} For every $X\in \csts$, let $\cyl^\csts(X):= \CSA_1
  \CSA_2\cyl(X)$. \end{nota}

Let $X$ be an object of $\csts$. The cylinder $\cyl^\csts(X)$ of $X$
in $\csts$ is then the Cattani-Sassone transition systems with the set
of states $\Int(X)\p\{0\} \sqcup \Ext(X) \p \{0,1\}$, the set of
actions $\Ac(X)$ with the same labelling map as $X$, and such that a
tuple \[((\alpha,\epsilon_0),u_1,\dots,u_n,(\beta,\epsilon_{n+1}))\]
is a transition of $\cyl^\csts(X)$ if and only if the tuple
$(\alpha,u_1,\dots,u_n,\beta)$ is a transition of $X$.  The canonical
map $\gamma^\epsilon:X \to \cyl^\csts(X)$ with $\epsilon \in \{0,1\}$
is induced by the mapping $\alpha \mapsto (\alpha,\epsilon)$ for
$\alpha \in \St(X)\backslash \Int(X)$ or $\epsilon=0$, by the mapping
$\alpha \mapsto (\alpha,0)$ for $\alpha \in \Int(X)$ and $\epsilon=1$
and by the identity of $\Ac(X)$ on actions. The canonical map
$\cyl^\csts(X)\to X$ is induced by the mapping
$(\alpha,\epsilon)\mapsto \alpha$ on states and by the identity of
$\Ac(X)$ on actions.

\bp \label{another-section} For every Cattani-Sassone
  transition system $X$, the unit map \[\CSA_2 \cyl(X) \to \CSA_1
  \CSA_2\cyl(X)\] is split epic.
\ep

\bpf For every Cattani-Sassone transition system $X$, the map \[\CSA_2
\cyl(X) \to \CSA_1 \CSA_2\cyl(X)\] is the identity on states and takes
the action $(u,\epsilon)$ of $\CSA_2 \cyl(X)$ to the action $u$ of
$\CSA_1 \CSA_2\cyl(X)$. The inclusion map \[\CSA_1 \CSA_2\cyl(X)
\subset \CSA_2 \cyl(X)\] induced by the identity on states and the
mapping $u\mapsto (u,0)$ on actions yields a section of the map
$\CSA_2 \cyl(X) \to \CSA_1 \CSA_2\cyl(X)$.  \epf

\bd \label{boundary-def} Let $n\geq 1$ and $x_1,\dots,x_n \in
\Sigma$. Let $\de C_n[x_1,\dots,x_n]$ be the regular transition system
defined by removing from the $n$-cube $C_n[x_1,\dots,x_n]$ all its
$n$-transitions. It is called the {\rm boundary} of
$C_n[x_1,\dots,x_n]$. \ed

\begin{nota} \label{fix}
  Let \begin{multline*}
\I^{\cts} = \{C:\varnothing \to \{0\}\} \cup \{\de
  C_n[x_1,\dots,x_n] \subset C_n[x_1,\dots,x_n]\mid \hbox{$n\geq 1$
    and $x_1,\dots,x_n \in \Sigma$}\} \\\cup \{C_0 \sqcup C_0 \sqcup
  C_1[x] \to \dd{x}\mid x\in \Sigma\}
\end{multline*}
where the maps $C_0 \sqcup C_0 \sqcup C_1[x] \to \dd{x}$ for $x$
running over $\Sigma$ are defined to be bijective on states.
\end{nota}

The definition of $\I^{\cts}$ is not the same as the one of
\cite{biscsts1}. The maps $C_1[x]\to \dd{x}$ for $x$ running over
$\Sigma$ are replaced by the maps
$C_0 \sqcup C_0 \sqcup C_1[x] \to \dd{x}$ for $x$ running over
$\Sigma$ which are defined to be bijective on states.  The class of
maps $\cell_\cts(\I^\cts)$ with $\I^\cts$ defined as above is exactly
the class of monomorphisms of weak transition systems (i.e. one-to-one
on states and actions by \cite[Proposition~3.1]{cubicalhdts}) between
cubical transition systems by \cite[Theorem~3.7]{biscsts1} (which is
fixed in Appendix~\ref{erratum}).  This choice is more convenient
because the maps $C_0 \sqcup C_0 \sqcup C_1[x]\to \dd{x}$ for
$x\in \Sigma$ are bijective on states and on actions. Note that the
maps $C_0 \sqcup C_0 \sqcup C_1[x]\to \dd{x}$ for $x$ running over
$\Sigma$ are already used e.g. in the proof of
\cite[Theorem~4.6]{homotopyprecubical}.

Let $X$ be a weak transition system. By
\cite[Proposition~2.14]{biscsts1}, there exists a well-defined weak
transition system $\cocyl(X)$ with the set of states $\St(X)\p
\St(X)$, the set of actions is the set $\Ac(X)\p_\Sigma \Ac(X)$ and
the labelling map is the composite map $\Ac(X)\p_\Sigma \Ac(X) \to
\Ac(X) \to \Sigma$ and such that a tuple
$((\alpha^0,\alpha^1),(u_1^0,u_1^1),\dots,(u_n^0,u_n^1),(\beta^0,\beta^1))$
is a transition if and only if for any
$\epsilon_0,\dots,\epsilon_{n+1}\in \{0,1\}$, the tuple
$(\alpha^{\epsilon_0},u_1^{\epsilon_1},\dots ,u_n^{\epsilon_n},
\beta^{\epsilon_{n+1}})$ is a transition of $X$. The functor
$\cocyl:\wts \to \wts$ is a right adjoint of the functor $\cyl:\wts\to
\wts$ by \cite[Proposition~2.15]{biscsts1}. The right adjoint
$\cocyl^{\cts}:\cts \to \cts$ of the restriction of $\cyl$ to the full
subcategory of cubical transition systems is the composite
map \[\cocyl^{\cts}:\cts \subset \wts \stackrel{\cocyl}
\longrightarrow \wts \longrightarrow \cts\] where the right-hand map
is the coreflection.

\bp \label{subpath} Let $X$ be a cubical transition system. The
cubical transition system $\cocyl^{\cts}(X)$ can be identified to a
subobject of $\cocyl(X)$ having the same set of states. \ep

\bpf By \cite[Proposition~3.4]{biscsts1}, the counit map
$\cocyl^{\cts}(X)\to \cocyl(X)$ is bijective on states and one-to-one
on actions and transitions.  \epf

\bp \label{csa1-path} For every Cattani-Sassone transition system $X$,
the cubical transition system $\cocyl^\cts(X)$ is a Cattani-Sassone
regular transition system. \ep

\bpf Since $X$ is regular, the cubical transition system
$\cocyl^\cts(X)$ is regular as well by
\cite[Proposition~3.10]{biscsts1}.  Let
\[
((\alpha^-,\alpha^+),(u^-,u^+),(\beta^-,\beta^+)),((\alpha^-,\alpha^+),(v^-,v^+),(\beta^-,\beta^+))\]
be two transitions of $\cocyl^\cts(X)$. By definition of
$\cocyl^\cts(X)$, the tuples $ (\alpha^-,u^\pm,\beta^-)$ and
$(\alpha^-,v^\pm,\beta^-)$ are transitions of $X$. Since $X$ satisfies
CSA1, we deduce that $u^-=u^+=v^+=v^-$. We obtain
$(u^-,u^+)=(v^-,v^+)$, which means that $\cocyl^\cts(X)$ satisfies
CSA1.  \epf

\bth \label{csts-model} There exists a left determined model structure
on $\csts$ with the set of generating cofibrations $\I^{\cts}$. This
model category is an Olschok model category with the very good
cartesian cylinder $\cyl^\csts:\csts \to \csts$. All objects are
cofibrant. \eth

\bpf By \cite[Theorem~3.16]{biscsts1}, there exists a unique left
determined model structure on $\rts$ such that the set of generating
cofibrations is $\CSA_2(\I^{\cts})=\I^{\cts}$. This model structure is
an Olschok model structure with the very good cylinder $\CSA_2 \cyl$
and, by \cite[Proposition~3.10]{biscsts1}, with the cocylinder
$\cocyl^\cts$.  We want to restrict this model structure to the
subcategory $\csts$.  1) It turns out that $\csts$ is a reflective
subcategory of $\rts$ by Proposition~\ref{csa1-reflective}. 2) We have
$\CSA_1(\I^\cts)=\I^\cts$ since every domain and every codomain of a
map of $\I^\cts$ satisfies CSA1 and since $\csts$ is a full
subcategory of $\rts$. 3) We have $\cocyl^\cts(\csts)\subset \csts$ by
Proposition~\ref{csa1-path}. 4) The map
$\CSA_2 \cyl(X) \to \CSA_1 \CSA_2\cyl(X)$ is split by
Proposition~\ref{another-section}. Therefore, we can apply
\cite[Theorem~3.1]{leftdet} and the proof is complete.  \epf

The following proposition will be used later: 

\bp \label{121-cof} Let $f:X\to Y$ be a map of $\csts$. If $f$ is
one-to-one on states and on actions, then it belongs to
$\cell_\csts(\I^{\cts})$. The converse is false: there exist maps of
$\cell_\csts(\I^{\cts})$ which are not one-to-one on states and on
actions.  \ep

\bpf Since $f$ is one-to-one on states and on actions, it belongs to
$\cell_\cts(\I^{\cts})$. The inclusion functor $\csts \subset \cts$
has a left adjoint $\CSA_1\CSA_2:\cts \to \csts$. Thus, the functor
$\CSA_1\CSA_2$ is colimit-preserving and
$\CSA_1\CSA_2(f):\CSA_1\CSA_2(X) \to \CSA_1\CSA_2(Y)$ belongs to
$\cell_\csts(\CSA_1\CSA_2(\I^{\cts}))$. Since $\CSA_1\CSA_2(f)=f$ and
$\CSA_1\CSA_2(\I^{\cts}) = \I^{\cts}$, we deduce that $f$ belongs to
$\cell_\csts(\I^{\cts})$.  Let $x\in \Sigma$. The map
$\gamma_{C_2[x,x]} : C_2[x,x] \sqcup C_2[x,x] \to
\cyl^\csts(C_2[x,x])$
is a cofibration of $\csts$. The state $01$ is an internal state of
$C_2[x,x]$. Thus, the states $01$ of the two copies of $C_2[x,x]$ are
identified to the same state of $\cyl^\csts(C_2[x,x])$. This means
that $\gamma_{C_2[x,x]}$ is not one-to-one on states. The actions
$(x,1)$ of the two copies of $C_2[x,x]$ are identified to the same
action $(x,1)$ of $\cyl^\csts(C_2[x,x])$. This implies that
$\gamma_{C_2[x,x]}$ is not one-to-one on actions.  \epf

\section{Computation of the path functor}
\label{section-path-calcul}

\bp \label{ps-construct} Let $X$ be a regular transition
system. Consider the family of sets $(\mathcal{T}_n)_{n\geq 1}$
constructed by induction on $n$ as follows:
\begin{itemize}
\item $\mathcal{T}_1$ is the set of triples
  $((\alpha^-,\alpha^+),u_1,(\beta^-,\beta^+))$ such that the
  triples $(\alpha^\pm,u_1,\beta^\pm)$ are transitions of $X$.
\item For $n\geq 2$, $\mathcal{T}_n$ is the set of tuples
  $((\alpha^-,\alpha^+),u_1,\dots,u_n,(\beta^-,\beta^+))$ such that
  the tuples $(\alpha^\pm,u_1,\dots,u_n,\beta^\pm)$ are transitions of
  $X$ and such that for all permutations $\sigma$ of $\{1,\dots,n\}$
  and all $1\leq p <n$, there exists a pair of states
  $(\gamma_\sigma^-,\gamma_\sigma^+)$ of $X$ such that the tuple
  $((\alpha^-,\alpha^+),u_{\sigma(1)},\dots,u_{\sigma(p)},(\gamma_\sigma^-,\gamma_\sigma^+))$
  belongs to $\mathcal{T}_p$ and such that the tuple
  $((\gamma_\sigma^-,\gamma_\sigma^+),u_{\sigma(p+1)},\dots,u_{\sigma(n)},(\beta^-,\beta^+))$
  belongs to $\mathcal{T}_{n-p}$.
\end{itemize}
Let $\mathcal{T}=\bigcup_{n\geq 1} \mathcal{T}_n$. Then we have:
\begin{enumerate}
\item[a)] For every transition $(\alpha,u_1,\dots,u_n,\beta)$ of $X$, the
  tuple $((\alpha,\alpha),u_1,\dots,u_n,(\beta,\beta))$ belongs to
  $\mathcal{T}$.
\item[b)] The set of states $\St(X)\p \St(X)$, the labelling map $\mu:\Ac(X)\to
  \Sigma$ and the set of tuples $\mathcal{T}$ assemble to a regular transition
  system denoted by $\pscocyl(X)$.
\end{enumerate} \ep

\bpf {a)} Consider the statement $\mathcal{P}_n$ defined for $n\geq 1$ by:
\begin{center}
\emph{
  for every transition $(\alpha,u_1,\dots,u_m,\beta)$ of $X$ with
  $m\leq n$, the tuple $((\alpha,\alpha),u_1,\dots,u_m,(\beta,\beta))$
  belongs to $\mathcal{T}$.}
\end{center} Let $(\alpha,u_1,\beta)$ be a transition of $X$. The
tuple $((\alpha,\alpha),u_1,(\beta,\beta))$ belongs to $\mathcal{T}_1$ by
definition of $\mathcal{T}_1$. We have proved $\mathcal{P}_1$. Let $n\geq
2$. Assume $\mathcal{P}_{n-1}$. We want to prove $\mathcal{P}_n$. Let
$(\alpha,u_1,\dots,u_n,\beta)$ be a transition of $X$. Let $1\leq p
<n$. Let $\sigma$ be a permutation of $\{1,\dots,n\}$. Since $X$ is
cubical, there exists a state $\gamma$ of $X$ such that the tuples
$(\alpha,u_{\sigma(1)},\dots,u_{\sigma(p)},\gamma)$ and
$(\gamma,u_{\sigma(p+1)},\dots,u_{\sigma(n)},\beta)$ are transitions
of $X$. By the induction hypothesis, the tuples
$((\alpha,\alpha),u_{\sigma(1)},\dots,u_{\sigma(p)},(\gamma,\gamma))$
and
$((\gamma,\gamma),u_{\sigma(p+1)},\dots,u_{\sigma(n)},(\beta,\beta))$
belong to $\mathcal{T}$ since $p\leq n-1$ and $n-p \leq n-1$. Hence we have
proved $\mathcal{P}_n$ from $\mathcal{P}_{n-1}$.

{b)} The set of tuples $\mathcal{T}$ satisfies the multiset axiom
because of the internal symmetry of its definition. Let
$((\alpha^-,\alpha^+),u_1,\dots,u_n,(\beta^-,\beta^+))$ be an element
of $\mathcal{T}$ with $n\geq 3$. Let $p,q\geq 1$ with $p+q<n$ such that the
tuples
\begin{align*}
 &((\alpha^-,\alpha^+),u_1, \dots, u_n, (\beta^-,\beta^+)), \\ &((\alpha^-,\alpha^+),u_1, \dots, u_p,
  (\nu_1^-,\nu_1^+)), ((\nu_1^-,\nu_1^+), u_{p+1}, \dots, u_n, (\beta^-,\beta^+)),\\ &((\alpha^-,\alpha^+), u_1,
  \dots, u_{p+q}, (\nu_2^-,\nu_2^+)), ((\nu_2^-,\nu_2^+), u_{p+q+1}, \dots, u_n, (\beta^-,\beta^+))
\end{align*}
belong to $\mathcal{T}$. Then by definition of $\mathcal{T}$, the tuples 
\begin{align*}
 &(\alpha^\pm,u_1, \dots, u_n, \beta^\pm), \\ &(\alpha^\pm,u_1, \dots, u_p,
  \nu_1^\pm), (\nu_1^\pm, u_{p+1}, \dots, u_n, \beta^\pm),\\ &(\alpha^\pm, u_1,
  \dots, u_{p+q}, \nu_2^\pm), (\nu_2^\pm, u_{p+q+1}, \dots, u_n, \beta^\pm)
\end{align*}
are transitions of $X$. Since $X$ satisfies CSA2, we obtain
$\nu_1^-=\nu_1^+=\nu_1$ and $\nu_2^-=\nu_2^+=\nu_2$.  By the patching
axiom, the tuple $(\nu_1, u_{p+1}, \dots, u_{p+q}, \nu_2)$ is a
transition of $X$. By a), we deduce that the $((\nu_1^-,\nu_1^+),
u_{p+1}, \dots, u_{p+q}, (\nu_2^-,\nu_2^+))$ belongs to
$\mathcal{T}_{q}$. We have proved that $\mathcal{T}$ satisfies the
patching axiom and that the set of states $\St(X)\p \St(X)$, the
labelling map $\mu:\Ac(X)\to \Sigma$ and the set of tuples
$\mathcal{T}$ assemble to a weak transition system $\pscocyl(X)$. Let
$u$ be an action of $X$. Since $X$ is cubical, there exists a
transition $(\alpha,u,\beta)$ of $X$. Thus, the triple
$((\alpha,\alpha),u,(\beta,\beta))$ is a transition of
$\pscocyl(X)$. Hence, all actions of $\pscocyl(X)$ are used.  By
definition of $\mathcal{T}$, $\pscocyl(X)$ satisfies the Intermediate
state axiom. This means that $\pscocyl(X)$ is a cubical transition
system. Since $X$ is regular, it satisfies CSA2.  This means that in
the definition of $\mathcal{T}$, the equality
$\gamma_\sigma^-=\gamma_\sigma^+$ always holds and that this state is
unique. Hence the cubical transition system $\pscocyl(X)$ satisfies
CSA2 as well.  \epf

The natural mapping $X\mapsto \pscocyl(X)$ yields a well-defined
functor \[\pscocyl:\rts \to \rts.\]

\bth \label{path-calcul} Let $X$ be a Cattani-Sassone transition
system. There exists a natural isomorphism $\pscocyl(X) \iso
\cocyl^\cts(X)$.  \eth

\bpf By definition of $\cocyl:\wts \to \wts$, the identity of
$\St(X)\p \St(X)$ and the diagonal map $\Ac(X) \to \Ac(X)\p_\Sigma
\Ac(X)$ induces a map of weak transition systems \[\pscocyl(X) \to
\cocyl(X).\] This map is one-to-one on states and on actions, and
therefore one-to-one on transitions by
Proposition~\ref{121-121-trans}.  By Proposition~\ref{ps-construct},
the weak transition system $\pscocyl(X)$ is cubical. Since $\cts$ is a
coreflective subcategory of $\wts$, the map $\pscocyl(X) \to
\cocyl(X)$ then factors uniquely as a composite \[\pscocyl(X) \to
\cocyl^\cts(X)\to \cocyl(X).\] By Proposition~\ref{subpath}, the
cubical transition system $\cocyl^{\cts}(X)$ has the set of states
$\St(X)\p \St(X)$ and the set of actions (of transitions resp.)  of
$\cocyl^{\cts}(X)$ is a subset of the set of actions (of transitions
resp.) of $\cocyl(X)$, i.e. $\Ac(X)\p_\Sigma \Ac(X)$.  We deduce that
the map $\pscocyl(X) \to \cocyl^\cts(X)$ induces the identity of
$\St(X)\p \St(X)$ on states and the diagonal of $\Ac(X)$ on actions.
Therefore, the set of actions of $\cocyl^{\cts}(X)$ contains the
diagonal of $\Ac(X)$. Let $(u^-,u^+) \in \Ac(X)\p_\Sigma \Ac(X)$ be an
action of $\cocyl^{\cts}(X)$. Since $\cocyl^{\cts}(X)$ is cubical,
there exists a transition
$((\alpha^-,\alpha^+),(u^-,u^+),(\beta^-,\beta^+))$ of
$\cocyl^{\cts}(X)$. Since the latter tuple is also a transition of
$\cocyl(X)$ by Proposition~\ref{subpath}, the triples
$(\alpha^\pm,u^\pm,\beta^\pm)$ are transitions of $X$ by definition of
$\cocyl:\wts \to \wts$. Since $X$ satisfies CSA1 by hypothesis, we
deduce that $u^-=u^+$. Thus, the set of actions of $\cocyl^{\cts}(X)$
is a subset of the diagonal of $\Ac(X)$. This implies that the set of
actions of $\cocyl^{\cts}(X)$ is exactly the diagonal of $\Ac(X)$.  We
deduce that the map $\pscocyl(X) \to \cocyl^\cts(X)$ induces the
bijection $\Ac(X) \iso \{(u,u)\mid u\in \Ac(X)\}$ on actions. For the
sequel, the diagonal of $\Ac(X)$ can be identified with the set
$\Ac(X)$.

We have proved so far that the map $\pscocyl(X) \to \cocyl^\cts(X)$ is
bijective on states (it is the identity of $\St(X)\p \St(X)$) and
bijective on actions (it is the bijection $\Ac(X) \iso \{(u,u)\mid
u\in \Ac(X)\}$). It is then one-to-one on transitions by
Proposition~\ref{121-121-trans}.  Let
$((\alpha^-,\alpha^+),(u_1,u_1),\dots,(u_n,u_n),(\beta^-,\beta^+))$ be
a transition of $\cocyl^\cts(X)$. It is a transition of $\cocyl(X)$ by
Proposition~\ref{subpath}. By definition of $\cocyl:\wts \to \wts$,
the tuples $(\alpha^\pm,u_1,\dots,u_n,\beta^\pm)$ are transitions of
$X$. If $n=1$, then the tuple
$((\alpha^-,\alpha^+),u_1,(\beta^-,\beta^+))$ is a transition of
$\pscocyl(X)$ by definition of $\pscocyl(X)$. Assume that $n\geq 2$.
Let $\sigma$ be a permutation of $\{1,\dots,n\}$. Let $1\leq p <n$.
Since $\cocyl^\cts(X)$ is cubical, there exists a pair of states
$(\gamma_\sigma^-,\gamma_\sigma^+)$ of $X$ such that the tuples
\[
((\alpha^-,\alpha^+),(u_{\sigma(1)},u_{\sigma(1)}),\dots,(u_{\sigma(p)},u_{\sigma(p)}),(\gamma_\sigma^-,\gamma_\sigma^+))\]
and 
\[
((\gamma_\sigma^-,\gamma_\sigma^+),(u_{\sigma(p+1)},u_{\sigma(p+1)}),\dots,(u_{\sigma(n)},u_{\sigma(n)}),(\beta^-,\beta^+))\]
are transitions of $\cocyl^\cts(X)$. By definition of $\cocyl:\wts\to
\wts$, this implies that the tuples
$(\alpha^\pm,u_{\sigma(1)},\dots,u_{\sigma(p)},\gamma_\sigma^\pm)$ and
$(\gamma_\sigma^\pm,u_{\sigma(p+1)},\dots,u_{\sigma(n)},\beta^\pm)$
are transitions of $X$.  Therefore the tuple
$((\alpha^-,\alpha^+),(u_1,u_1),\dots,(u_n,u_n),(\beta^-,\beta^+))$ is
a transition of $\pscocyl(X)$ by Proposition~\ref{ps-construct}. We
have proved that the map $\pscocyl(X) \to \cocyl^\cts(X)$ is onto on
transitions.  \epf

The path object of $X$ in $\csts$ is then the Cattani-Sassone
transition system denoted by $\cocyl^\csts(X)$ with the set of states
$\St(X)\p \St(X)$, the set of actions $\Ac(X)$ and such that a tuple
\[((\alpha^-,\alpha^+),u_1,\dots,u_n,(\beta^-,\beta^+))\] is a
transition of $\cocyl^\csts(X)$ if and only if the tuples
$(\alpha^\pm,u_1,\dots,u_n,\beta^\pm)$ are transitions of $X$. The
canonical map $\tau_X:X \to \cocyl^\csts(X)$ is induced by the
diagonal of $\St(X)$ on states and the identity on actions. The
canonical map $\pi_X^\epsilon: \cocyl^\csts(X) \to X$ with $\epsilon
\in \{0,1\}$ is induced by the projection on the $(\epsilon+1)$-th
component $\St(X) \p \St(X) \to \St(X)$ on states and by the identity
on actions.

Table~\ref{csts-cyl-path} summarizes the computation of the cylinder
functor and of the path functor of $\csts$.

\begin{table}
\begin{tabular}{|c|c|c|}
\hline
&&\\
 $X$ & $\cyl^\csts(X)$ & $\cocyl^\csts(X)$ \\&&\\
 $\St(X)$ & $\Int(X)\p \{0\} \sqcup \Ext(X)\p \{0,1\}$ & $\St(X) \p \St(X)$ \\&&\\
 $\Ac(X)$ & $\Ac(X)$ &  $\Ac(X)$ \\&&\\
 $\Tr(X)$ & $((\alpha,\epsilon_0),u_1,\dots,u_n,(\beta,\epsilon_{n+1}))$ & $((\alpha^-,\alpha^+),u_1,\dots,u_n,(\beta^-,\beta^+))$ \\
& such that $(\alpha,u_1,\dots,u_n,\beta) \in \Tr(X)$ & such that $(\alpha^\pm,u_1,\dots,u_n,\beta^\pm) \in \Tr(X)$\\
&&\\
\hline
\end{tabular}\newline\newline
\caption{Cylinder functor and path functor of $\csts$}
\label{csts-cyl-path}
\end{table}

\section{Characterization in the non star-shaped case}
\label{description-csts-section}

\bp \label{p2} Let $X$ be a fibrant object of $\csts$. Then the
labelling map $\mu:\Ac(X) \to \Sigma$ is one-to-one. \ep

\bpf Let $u$ and $v$ be two actions of $X$ with $\mu(u)=\mu(v)$. Since
$X$ is cubical, there exist two transitions $(\alpha_u,u,\beta_u)$ and
$(\alpha_v,v,\beta_v)$ of $X$. Consider the commutative diagram of
$\csts$
\[
\xymatrix@C=6em@R=3em
{
\{0_1,1_1\} = \de C_1[\mu(u)] \fr{\hbox{\tiny $\begin{array}{c}0_1\mapsto (\alpha_u,\alpha_v)\\ 1_1\mapsto (\beta_u,\beta_v)\end{array}$}} \fd{} & \cocyl^\csts(X) \fd{\pi^0} \\
C_1[\mu(u)] \ar@{-->}[ru]^-{\ell} \ar@{->}[r]_-{\hbox{\tiny $\begin{array}{c}0_1\mapsto \alpha_u\\ 1_1\mapsto \beta_u\\(\mu(u),1) \mapsto u\end{array}$}} & X.
}
\]
Since $X$ is fibrant by hypothesis, it is injective with respect to
the anodyne cofibration $(\de C_1[\mu(u)] \subset C_1[\mu(u)])\star
\gamma^0$. By adjunction, we deduce that the lift $\ell$ in the
diagram above exists. Therefore the triple $\ell(0_1,(\mu(u),1),1_1) =
((\alpha_u,\alpha_v),u,(\beta_u,\beta_v))$ is a transition of
$\cocyl^\csts(X)$. This implies that the triple $(\alpha_v,u,\beta_v)$
is a transition of $X$. Since $X$ satisfies CSA1, we deduce that
$u=v$.  \epf

\bp \label{p3} Let $X$ and $Y$ be two objects of $\csts$ with
  $Y$ fibrant. Then the set $\pi_\csts(X,Y)$ has at most one element:
  all maps from $X$ to $Y$ are homotopy equivalent.
\ep

\bpf Let $f,g:X\rightrightarrows Y$ be two maps of $\csts$.  We have
$\mu(f(u))=\mu(g(u))=\mu(u)$ since maps of $\csts$ preserve
labels. Since $Y$ is fibrant, $\mu$ is one-to-one by
Proposition~\ref{p2}. Thus, we deduce that $f$ and $g$ coincide on
actions.  Since $X$ is cofibrant by Theorem~\ref{csts-model} and $Y$
fibrant by hypothesis, it suffices to prove that $f$ and $g$ are right
homotopic to conclude the proof.  The only possible definition of this
right homotopy is $H(\alpha)=(f(\alpha),g(\alpha))$ for all states
$\alpha$ of $X$ and $H(u)=f(u)=g(u)=\mu(u)$ for all actions $u$ of
$X$.  We have to prove that $H$ yields a well-defined map of
transition systems, i.e. that for all transitions
$(\alpha,u_1,\dots,u_n,\beta)$ of $X$, the tuple
$((f(\alpha),g(\alpha)),H(u_1),\dots,H(u_n),(f(\beta),g(\beta)))$ is a
transition of $\cocyl^\csts(Y)$.

The $n$-transition $(f(\alpha),f(u_1),\dots,f(u_n),f(\beta))$
  gives rise to a (unique) map of weak transition systems
  $C^{ext}_n[\mu(u_1),\dots,\mu(u_n)] \to Y$ (the definition of
  $C^{ext}_n[\mu(u_1),\dots,\mu(u_n)]$ is recalled in
  Section~\ref{reminder}). Since $Y$ satisfies CSA2, it factors
  (uniquely) as a composite
\[C^{ext}_n[\mu(u_1),\dots,\mu(u_n)] \to
C_n[\mu(u_1),\dots,\mu(u_n)]\to Y\]
by \cite[Theorem~5.6]{hdts}. We obtain the commutative diagram of
$\csts$
\[
\xymatrix@C=6em@R=3em
{
\{0_n,1_n\} \fr{\hbox{\tiny $\begin{array}{c}0_n\mapsto (f(\alpha),g(\alpha))\\ 1_n\mapsto (f(\beta),g(\beta))\end{array}$}} \fd{} & \cocyl^\csts(Y) \fd{\pi^0} \\
C_n[\mu(u_1),\dots,\mu(u_n)] \ar@{-->}[ru]^-{\ell} \ar@{->}[r]_-{} & Y.
}
\]
Since $Y$ is fibrant by hypothesis, it is injective with respect to
the anodyne cofibration
$(\{0_n,1_n\}\subset C_n[\mu(u_1),\dots,\mu(u_n)]) \star \gamma^0$.
By adjunction, we deduce that the lift $\ell$ in the diagram above
exists. By Table~\ref{csts-cyl-path}, the right vertical map $\pi^0$
is the identity on actions. We deduce that the tuple
\[\ell(0_n,(\mu(u_1),1),(\mu(u_n),n),1_n)=((f(\alpha),g(\alpha)),f(u_1),\dots,f(u_n),(f(\beta),g(\beta)))\]
is a transition of $\cocyl^\csts(Y)$.
\epf

\bp \label{p4} Let $w:X \to Y$ be a map of $\csts$ which is onto on
states, on actions and on transitions. Then $w$ is a weak equivalence.
\ep

\bpf Let $Z$ be a fibrant object. We have to prove that the set map
$\pi_\csts(Y,Z)\to \pi_\csts(X,Z)$ induced by the precomposition with
$w$ is bijective. By Proposition~\ref{p3}, it suffices to prove that
it is onto. Suppose first that $Z=\varnothing$. Let $f:X\to Z$. Then
$X=\varnothing$. Since $w:X \to Y$ is onto on states, $Y$ does not
contain any state. Consequently, $Y=\varnothing$ because all actions
of $Y$ are used. We deduce that $w=\id_\varnothing$, and that
$\pi_\csts(Y,Z) = \pi_\csts(X,Z) =
\pi_\csts(\varnothing,\varnothing)=\{\id_\varnothing\}$.
Suppose now that $Z\neq\varnothing$. By Proposition~\ref{p2}, we have
the inclusion $\Ac(Z)\subset \Sigma$. Let $f:X\to Z$ be a map of
$\csts$. Let $\xi$ be a state of $Z$. Let $g:Y\to Z$ defined on states
by $g(\alpha)=\xi$ and on actions by $g(u)=\mu(u)$. It suffices to
prove that $g$ yields a well-defined map of $\csts$ from $Y$ to $Z$ to
complete the proof. Let $(\alpha,u_1,\dots,u_n,\beta)$ be a transition
of $Y$. By hypothesis, there exists a transition
$(\overline{\alpha},\overline{u_1},\dots,\overline{u_n},\overline{\beta})$
of $X$ with
$w(\overline{\alpha},\overline{u_1},\dots,\overline{u_n},\overline{\beta})
= (\alpha,u_1,\dots,u_n,\beta)$.
We obtain the transition
$(f(\overline{\alpha}),\mu(u_1),\dots,\mu(u_n),f(\overline{\beta}))$
of $Z$ which yields a map of cubical transition systems
$C_n[\mu(u_1),\dots,\mu(u_n)] \to Z$. Since $Z$ is fibrant, it is
injective with respect to the anodyne cofibration
$(\{0_n,1_n\} \subset C_n[\mu(u_1),\dots,\mu(u_n)]) \star
\gamma^0$.
By adjunction, we deduce that the lift $\ell$ exists in the
commutative diagram of solid arrows of $\csts$
\[
\xymatrix@C=4em@R=3em
{
\{0_n,1_n\}  \fr{\hbox{\tiny $\begin{array}{c}0_n\mapsto (f(\overline{\alpha}),\xi)\\ 1_n\mapsto (f(\overline{\beta}),\xi)\end{array}$}} \fd{} & \cocyl^\csts(Z) \fd{\pi^0} \\
C_n[\mu(u_1),\dots,\mu(u_n)] \ar@{-->}[ru]^-{\ell} \ar@{->}[r]_-{\hbox{\tiny $\begin{array}{c}0_n\mapsto f(\overline{\alpha})\\ 1_n\mapsto f(\overline{\beta})\\(\mu(u_i),i) \mapsto \mu(u_i)\end{array}$}} & Z.
}
\]
Therefore the tuple
$((f(\overline{\alpha}),\xi),\mu(u_1),\dots,\mu(u_n),(f(\overline{\beta}),\xi))$
is a transition of $\cocyl^\csts(Z)$. We deduce that the tuple
$(\xi,\mu(u_1),\dots,\mu(u_n),\xi)$ is a transition of $Z$. We have
proved that $g(\alpha,u_1,\dots,u_n,\beta) =
(\xi,\mu(u_1),\dots,\mu(u_n),\xi)$ is a transition of $Z$ and that $g$
is a well-defined map of Cattani-Sassone transition systems.  \epf

\bth \label{cellR-weak} Let $R:\{0,1\} \to \{0\}$. Every map of
$\cell_\csts(\{R\})$ is onto on states, on actions and on
transitions. Every map of $\cell_\csts(\{R\})$ is a weak equivalence
of $\csts$. \eth

\bpf Every pushout of $R:\{0,1\} \to \{0\}$ is onto on states, on
actions and on transitions by Corollary~\ref{onto-trans-csts}. Every
transfinite composition of pushouts of $R:\{0,1\} \to \{0\}$ is also
onto on states, on actions and on transitions by
Theorem~\ref{colim-csts}. We deduce the first assertion of the
theorem. The second assertion of the theorem is then a corollary of
Proposition~\ref{p4}.  \epf

  Note that every map of $\CSA_1(\CSA_2(\I^\cts))=\I^\cts$ is a map
  between finitely presentable objects. Thus, by
  \cite[Proposition~4.1]{rankweak}, the class of weak equivalences of
  $\csts$ is closed under transfinite composition.

Let $X$ be a Cattani-Sassone transition system.  For all objects $X$
of $\csts$, the unique map $X\to \mathbf{1}$ from $X$ to the terminal
object of $\csts$ factors functorially as a
composite
\[\xymatrix@1@C=6em{X \fr{\in \cell_\csts(\{R\})}& R^\perp(X) \fr{\in
    \inj_\csts(\{R\})}& \mathbf{1}}\]
in $\csts$ where the left-hand map belongs to $\cell_\csts(\{R\})$ and
the right-hand map belongs to $\inj_\csts(\{R\})$. Since all maps of
$\cell_\csts(\{R\})$ are epic by Theorem~\ref{cellR-weak}, this
decomposition is unique up to isomorphism by
\cite[Proposition~A.1]{biscsts1}.  Since $R:\{0,1\}\to \{0\}$ is epic,
an object is $R$-injective if and only if it is $R$-orthogonal,
i.e. if and only if $X$ has at most one state. Let us denote by
$\overline{\csts}$ this small-orthogonality class. The functor
$R^\perp:\csts \to \overline{\csts}$ is a left adjoint of the
inclusion $\overline{\csts} \subset \csts$. By
\cite[Theorem~1.39]{MR95j:18001}, the category $\overline{\csts}$ is
locally presentable.

\bth \label{description-csts} The left adjoint $R^\perp:\csts \to
\inj_\csts(\{R\})$ induces a Quillen equivalence between the model
category $\csts$ and the category $\overline{\csts}$ equipped with the
discrete model category structure. \eth

\bpf For any $R$-orthogonal Cattani-Sassone transition system $X$, we
have by Table~\ref{csts-cyl-path} the equalities $\cocyl^\csts(X)=X$
and $R^\perp(\cyl^\csts(X))=X$. Using \cite[Theorem~3.1]{leftdet}, we
obtain a left determined Olschok model structure on $\overline{\csts}$
such that the cylinder and the path functors are the identity
functor. Therefore, we have the equalities
$\pi_{\overline{\csts}}(X,Y) = \pi^l_{\overline{\csts}}(X,Y) =
\pi^r_{\overline{\csts}}(X,Y) = \overline{\csts}(X,Y)$
for all $R$-orthogonal Cattani-Sassone transition systems $X$ and
$Y$. This means that the weak equivalences of $\overline{\csts}$ are
the isomorphisms.  For every $R$-orthogonal Cattani-Sassone transition
system $Y$, the canonical map $R^\perp(Y)\to Y$ is an isomorphism.  We
deduce that the left adjoint $R^\perp:\csts \to \overline{\csts}$ is a
homotopically surjective left Quillen adjoint from $\csts$ to
$\overline{\csts}$ equipped with the discrete model structure.  By
Theorem~\ref{cellR-weak}, the canonical map $X\to R^\perp(X)$ is a
weak equivalence of $\csts$ for all Cattani-Sassone transition systems
$X$. Thus, this left Quillen adjoint is a left Quillen equivalence.
\epf

\begin{cor} \label{w1} A map $f$ of $\csts$ is a weak equivalence if
  and only if $R^\perp(f)$ is an isomorphism. \end{cor}

\bpf Since all objects of $\csts$ are cofibrant by
Theorem~\ref{csts-model}, a weak equivalence $f$ is mapped to a weak
equivalence $R^\perp(f)$ of $\overline{\csts}$, i.e. an isomorphism.
Conversely, if $R^\perp(f)$ is an isomorphism, then by
Theorem~\ref{cellR-weak} and the two-out-of-three property, $f$ is a
weak equivalence of $\csts$.  \epf

\section{Homotopy theory of star-shaped objects}
\label{model-bullet}

A \emph{pointed (Cattani-Sassone) transition system} is a pair
$(X,*)$ where $X$ is a Cattani-Sassone transition system and
where $*$ is a state of $X$ called the \emph{base state}. A map
of pointed transition system is a map of $\csts$ preserving the base
state. The category of pointed (Cattani-Sassone) transition systems is
denoted by $\csts_*$. Let $\omega^*:\csts_* \to \csts$ be the
forgetful functor. It is a right adjoint which preserves connected
colimits. The left adjoint $\rho^*:\csts \rightarrow \csts_*$ is
defined on objects by $\rho^*(X) = (\{*\} \sqcup X,*)$ and on
morphisms by $\rho^*(f) = \id_{\{*\}} \sqcup f$.  By
\cite[Proposition~1.57]{MR95j:18001} and \cite[Theorem~2.7]{undercat},
there exists a structure of combinatorial model category on $\csts_*$
such that the cofibrations (the fibrations, the weak equivalences
resp.)  belong to the inverse image by the forgetful functor of the
class of cofibrations (fibrations, weak equivalences resp.) of
$\csts$. The set of generating cofibrations of the category $\csts_*$
is the set $\rho^*(\I^\cts)$.

Note that \emph{pointed weak, cubical or regular transition systems}
are defined exactly as pointed Cattani-Sassone transition
systems. Without further precision, a pointed transition system is
supposed to be a pointed Cattani-Sassone transition system.

\begin{rem} It is important to keep in mind that the functors
  $\omega:\wts \to \set^{\{s\}\cup \Sigma}$ and $\omega^*:\csts_* \to
  \csts$ are two different functors.
\end{rem}

\begin{nota} Let $\St(X,*)=\St(X)$, $\Ac(X,*)=\Ac(X)$ and
  $\Tr(X,*)=\Tr(X)$.
\end{nota}

Let $\cyl_*:\csts_* \to \csts_*$ be the functor defined by the
following natural pushout diagram of $\csts$ ($(X,*)$ being an
object of $\csts_*$):
\[
\xymatrix@C=3em {
  \{*\} \sqcup \{*\} \fr{} \fd{} & \{*\} \ar[dd]^-{\cyl_*(X)} \\
  X \sqcup X \fd{\gamma_X} & \\
  \cyl^\csts(X) \fr{} & \cocartesien
  \omega^*(\cyl_*(X,*)), }
\]
where the pushout is taken in $\csts$.

\bth \label{csts-pointed-model} The combinatorial model category
$\csts_*$ is a left determined Olschok model category with the very
good cylinder $\cyl_*:\csts_* \to \csts_*$.  \eth

\bpf By Proposition~\ref{121-cof}, every map with domain a singleton
is a cofibration of $\csts$.  The map $\{*\} \sqcup \{*\} \to
\cyl^\csts(\{*\})$ is an isomorphism. By \cite[Theorem~5.8]{leftdet},
we deduce that the model category $\csts_*$ is an Olschok model
category. Let $(X,*)$ be an object of $\csts_*$. We know that the set
of states of $\cyl^\csts(X)$ is $\Int(X) \p\{0\} \sqcup \Ext(X) \p
\{0,1\}$, that the set of actions of $\cyl^\csts(X)$ is $\Ac(X)$, and
that a tuple
$((\alpha,\epsilon_0),u_1,\dots,u_n,(\beta,\epsilon_{n+1}))$ is a
transition of $\cyl^\csts(X)$ if and only if the tuple
$(\alpha,u_1,\dots,u_n,\beta)$ is a transition of $X$.  Consider now
the pushout diagram of cubical transition systems
\[
\xymatrix@C=3em {
  \{*\} \sqcup \{*\} \fr{} \fd{} & \{*\} \fD{} \\
  X \sqcup X \fd{\gamma_X} & \\
  \cyl^\csts(X) \fr{} & \cocartesien
  U. }
\]
By Lemma~\ref{colim-preserv}, the forgetful functor mapping a cubical
transition system to its set of states (to its set of actions resp.)
is colimit-preserving. We deduce that the set of states of $U$ is
$(\Int(X) \cup \{*\}) \p\{0\} \sqcup (\St(X)\backslash (\Int(X) \cup
\{*\})) \p \{0,1\}$
and that the set of actions of $U$ is $\Ac(X)$. Therefore, the cubical
transition system $U$ is the $\omega$-final lift of the cocone of
$\set^{\{s\}\cup \Sigma}$ consisting of the unique
map
\[\omega(\cyl(X)) \longrightarrow ((\Int(X) \cup \{*\})\p \{0\} \sqcup
(\St(X)\backslash (\Int(X) \cup \{*\}))\p \{0,1\},\Ac(X)).\]
By Proposition~\ref{precalcul0}, we obtain
$U = \cyl(X)///(\Int(X) \cup \{*\})$. There is an inclusion
\[\cyl(X)///(\Int(X) \cup \{*\}) \subset \cyl(X)///\Int(X).\]
By Theorem~\ref{precalcul2}, we obtain the inclusion 
\[\cyl(X)///(\Int(X) \cup \{*\}) \subset \CSA_1\CSA_2 \cyl(X).\]
By Proposition~\ref{mor-cub-reg}, the cubical transition system 
$\cyl(X)///(\Int(X) \cup \{*\})$ is regular. And it satisfies CSA1 by 
Proposition~\ref{precalcul0}. Thus, the cubical transition system 
$\cyl(X)///(\Int(X) \cup \{*\})$ is a Cattani-Sassone transition system and we obtain 
the pushout diagram of $\csts$
\[
\xymatrix@C=3em {
  \{*\} \sqcup \{*\} \fr{} \fd{} & \{*\} \fD{} \\
  X \sqcup X \fd{\gamma_X} & \\
  \cyl^\csts(X) \fr{} & \cocartesien
  \cyl(X)///(\Int(X) \cup \{*\}). }
\] 
We deduce the isomorphism
\[\omega^*(\cyl_{*}(X,*)) \iso \cyl(X)///(\Int(X) \cup \{*\}).\]
The inclusion $\cyl(X)///(\Int(X) \cup \{*\}) \subset
\CSA_1\CSA_2\cyl(X)$ yields a section of the bottom horizontal
map. Thus, by \cite[Corollary~5.9]{leftdet}, the Olschok model
category $\csts_*$ is left determined and the functor
$\cyl_*:\csts_* \to \csts_*$ yields a very good cylinder. \epf

After the previous calculations, the following definition makes
sense:

\bd A state of a pointed transition system $(X,*)$ is {\rm internal}
if it is equal to the base state $*$ or it is internal in $X$. Let
$\Int(X,*)=\Int(X) \cup \{*\}$, and $\Ext(X,*)=\St(X) \backslash
\Int(X,*)$. A state which is not internal is {\rm external}.  Note
that for any Cattani-Sassone transition system $X$, there are the
equalities $\Int(\rho^*(X),*) = \Int(X) \sqcup \{*\}$ and
$\Ext(\rho^*(X),*) = \Ext(X)$.  \ed

Thus, we have for any pointed transition system $(X,*)$ the equalities
$\cyl^\csts(X)=\cyl(X)///\Int(X)$ and
$\omega^*(\cyl_*(X,*))=\cyl(X)///\Int(X,*)$. For any map $f:(X,*)\to
(Y,*)$ of pointed transition systems, we have $f(\Int(X,*)) \subset
\Int(Y,*)$: any internal state of $(X,*)$ is mapped to an internal
state of $(Y,*)$. In general, we have $f(\Ext(X,*)) \nsubseteq
\Ext(Y,*)$.

\bd \emph{The path} $(P(w),0)$ indexed by $w$ with $n\geq 0$,
$w=x_1\dots x_n \in \Sigma^n$ for $n\geq 0$ is by definition the
pointed transition system \[\xymatrix@C=3em{P(w)=0 \fr{(x_1,1)}& 1
  \fr{(x_2,2)}& 2 \dots n-1\fr{(x_n,n)}& n},\] which means that the
set of actions is $\{(x_i,i)\mid 1 \leq i \leq n\}$ with the labelling
map $\mu(x_i,i) = x_i$ for $1 \leq i \leq n$. Let $p^w:\rho^*(\{n\})
\subset (P(w),0)$ with $n \geq 0$ and $w\in \Sigma^n$ be the
inclusion. The state $0$ is also called the \emph{initial state} of
$(P(w),0)$ and the state $n$ the \emph{final state} of $(P(w),0)$. We
have the equality $(P(\varnothing),0)=\{0\}$. \ed

\bd Let $(X,*)$ be a pointed Cattani-Sassone transition system. A
state $\alpha$ of $X$ is {\rm reachable} if there exists $w\in
\Sigma^n$ with $n\geq 0$ and a map $(P(w),0) \to (X,*)$ taking $n$ to
$\alpha$.  A transition $(\alpha,u_1,\dots,u_n,\beta)$ of $(X,*)$ is
{\rm reachable} if its initial state $\alpha$ is reachable.  \ed

\bd A {\rm star-shaped (Cattani-Sassone) transition system} is a
pointed Cattani-Sassone transition system $(X,*)$ such that every
state of $X$ is reachable. The full subcategory of $\csts_*$ of
star-shaped transition systems is denoted by $\csts_{\bullet}$. \ed

Note that \emph{star-shaped weak, cubical or regular transition
  systems} are defined exactly as star-shaped Cattani-Sassone
transition systems. Without further precision, a star-shaped
transition system is supposed to be a star-shaped Cattani-Sassone
transition system.

The category $\csts_{\bullet}$ is a coreflective subcategory of
$\csts_*$. By \cite[Proposition~5.5]{biscsts1}, the coreflection
$\csts_* \to \csts_{\bullet}$ removes all states which are not
reachable, all actions which are not used by a reachable transition,
and all transitions which are not reachable.

\begin{nota} Let $(X,*)$ be a star-shaped transition system. Let
  $\alpha$ be a state of $(X,*)$. Let $\ell(\alpha) = \min \{n \in
  \mathbb{N} \mid \exists w\in \Sigma^n, \exists f:(P(w),0) \to (X,*),
  f(n)=\alpha\}$. The integer $\ell(\alpha)$ is well-defined since
  $(X,*)$ is star-shaped. We have $\ell(*)=0$.
\end{nota}

\bth \label{oslash} Consider a map $f\to g$ of $\Mor(\csts_*)$ (i.e. a
commutative square) where $f$ is a map of $\csts_*$ which is
one-to-one on states and on actions and where $g$ is a map of
$\csts_\bullet$. Then $f\to g$ factors as a composite $f\to
f\varoslash g \to g$ where $f\varoslash g$ is a map of $\csts_\bullet$
which is one-to-one on states and on actions. In particular,
$f\varoslash g$ is a cofibration of $\csts_\bullet$.  Finally, the
class of maps $\{f \varoslash g \mid g\in \Mor(\csts_\bullet)\}$ is a
set. \eth

Note that the factorization $f\to f\varoslash g \to g$ is not unique.

\bpf The proof of this theorem is similar to the proof of
\cite[Theorem~5.9]{biscsts1}.  This one is more explicit because we
need an explicit calculation of $f \varoslash g$ for the
sequel. Consider the commutative diagram of $\csts_*$
\[
\xymatrix@C=3em@R=3em
{
(A,*) \fd{f} \fr{\phi} & (X,*) \fd{g} \\
(B,*) \fr{\psi} & (Y,*)
}
\]
with $(X,*)$ and $(Y,*)$ star-shaped and $f$ one-to-one on states and
on actions. Since $(X,*)$ is star-shaped, for any state $\alpha$ of
$(A,*)$, the state $\phi(\alpha)$ is reachable from $*$ in $(X,*)$
using a path of length $\ell(\phi(\alpha))$ labelled with a finite
sequence $w_\alpha$ of length $\ell(\phi(\alpha))$ of $\Sigma$. Thus,
the composite map
\[\xymatrix@C=5em{\rho^*(\{\ell(\phi(\alpha))\}) \fr{\ell(\phi(\alpha)) \mapsto \alpha}& (A,*) \fr{\phi} & (X,*)}\] 
factors as a composite 
\[\xymatrix@C=5em{\rho^*(\{\ell(\phi(\alpha))\}) \fr{p^{w_\alpha}}& P(w_\alpha)\fr{} &(X,*)}.\]
We obtain the commutative diagram of $\csts_*$
\[
\xymatrix@C=3em@R=3em
{
(M,*)=\rho^*\left(\bigsqcup\limits_{\alpha\in \St(A)}\{\ell(\phi(\alpha))\}\right) \fd{\bigsqcup\limits_{\alpha\in \St(A)} p^{w_\alpha}} \fr{} & (A,*) \fd{\phi}  \\ 
(N,*)=\bigsqcup\limits_{\alpha\in \St(A)} P(w_\alpha) \fr{} & (X,*).
}
\]
For any state $\alpha$ of $\St(A)$, the state $g(\phi(\alpha))$ is
reachable from $*$ in $(Y,*)$ using a path of length
$\ell(\phi(\alpha))$ labelled with a finite sequence $w_\alpha$ of
length $\ell(\phi(\alpha))$ of $\Sigma$. Since $(Y,*)$ is star-shaped,
for any state $\beta$ of $\St(B)\backslash \St(A)$, the state
$\psi(\beta)$ is reachable from $*$ in $(Y,*)$ using a path of length
$\ell(\psi(\beta))$ labelled with a finite sequence $w_\beta$ of
length $\ell(\psi(\beta))$ of $\Sigma$. We obtain the commutative
diagram of $\csts_*$
\[
\xymatrix@R=5em
{
(P,*)=\rho^*\left(\bigsqcup\limits_{\alpha\in \St(A)}\{\ell(\phi(\alpha))\} \sqcup \bigsqcup\limits_{\beta\in \St(B)\backslash \St(A)}\{\ell(\psi(\beta))\}\right) \ar@{->}[d]|-{\bigsqcup\limits_{\alpha\in \St(A)} p^{w_\alpha} \sqcup\bigsqcup\limits_{\beta\in \St(B)\backslash \St(A)} p^{w_\beta}} \fr{} & (B,*) \fd{\psi}  \\ 
(Q,*)=\bigsqcup\limits_{\alpha\in \St(A)} P(w_\alpha) \sqcup \bigsqcup\limits_{\beta\in \St(B)\backslash \St(A)} P(w_\beta) \fr{} & (Y,*).
}
\]
We obtain the commutative diagram of $\csts_*$:
\[
\xymatrix@C=3em@R=3em
{
(M,*) \ar@{->}[rdd]\ar@{->}[d]_-{} \ar@{->}[rr]^-{} && (A,*) \ar@{->}[rdd]^/20pt/{f} \ar@{->}[r]^-{\phi} \ar@{->}[d]_-{} & X \ar@{->}[rrdd]^-{g}  \\ 
(N,*) \ar@{->}[rdd]_-{\subset} \ar@{->}[rr]^-{}|(0.25)\hole && \cocartesien
(\widehat{A},*)\ar@{-->}[rdd]_/10pt/{\exists !} \ar@{-->}[ru]^-{\exists !}|(0.32)\hole && \\
&(P,*) \ar@{->}[d]^-{} \ar@{->}[rr]^-{} && (B,*) \ar@{->}[rr]^-{\psi} \ar@{->}[d]_-{} && Y \\ 
&(Q,*) \ar@{->}[rr]^-{} && \cocartesien
(\widehat{B},*) \ar@{-->}[rru]^-{\exists !} &&
}
\]
The map $(\widehat{A},*) \to (\widehat{B},*)$ making the diagram commutative
exists by the universal property of the pushout and it is one-to-one
on states and on actions. By Proposition~\ref{121-cof}, the underlying
map is then a cofibration between Cattani-Sassone transition
system. Thus, since $(\widehat{A},*)$ and $(\widehat{B},*)$ are star-shaped by
construction, it is a cofibration of $\csts_\bullet$. We obtain the
factorization
\[
\xymatrix@C=5em@R=3em
{
(A,*) \ar@/^20pt/@{->}[rr]^-{\phi} \fr{\subset} \fd{f} & \widehat{A} \fd{f \varoslash g}\fr{} & X \fd{g} \\
(B,*) \ar@/_20pt/@{->}[rr]_-{\psi} \fr{\subset} & \widehat{B} \fr{} & Y.
}
\]
Finally, the class of maps $\{f \varoslash g \mid g\in
\Mor(\csts_\bullet)\}$ has at most $(\#\Sigma)^{(\#\St(B))\aleph_0}$
elements. Thus, it is a set.  \epf

\begin{cor} \label{generating-star-shaped} The set of generating cofibrations
$\rho^*(\I^\cts)$ of $\csts_*$ has a solution set of cofibrations with
respect to $\csts_\bullet$, i.e.  there exists a set $\I_\bullet$ of
cofibrations of $\csts_*$ between star-shaped objects such that every
map $i\to g$ from a generating cofibration $i$ of $\csts_*$ to a map
$g$ of $\csts_\bullet$ factors as a composite $i\to j\to g$ with $j\in
\I_\bullet$.  \end{cor}

\bpf The set $\I_\bullet = \{\rho^*(i) \varoslash g \mid i \in
\rho^*(\I^\cts) \hbox{ and }g \in \Mor(\csts_\bullet)\}$ is a
solution.  \epf

\bth \label{star-shaped-csts} There exists a left determined Olschok
model structure on the category $\csts_{\bullet}$ of star-shaped
transition systems with respect to the class of maps such that the
underlying map is a cofibration of $\csts$.  A very good cylinder is
given by the restriction of the functor $\cyl_*:\csts_* \to \csts_*$
to the coreflective subcategory $\csts_{\bullet}$ of $\csts_*$.  \eth

\bpf By Theorem~\ref{csts-pointed-model}, the combinatorial model
  category $\csts_*$ is a left determined Olschok model category.  We
  have seen that $\csts_\bullet$ is a coreflective subcategory of
  $\csts_*$.  By \cite[Lemma~A.3]{cubicalhdts},
  $\cof_{\csts_\bullet}(\I_\bullet) = \cof_{\csts_*}(\rho^*(\I^\cts))
  \cap \Mor(\csts_\bullet)$
  where $\I_\bullet$ is the set constructed in
  Corollary~\ref{generating-star-shaped}. The proof will be complete
  by applying \cite[Theorem~4.1]{leftdet}. It therefore remains to
  prove that $\cyl_*(\csts_\bullet)\subset \csts_\bullet$. Let $(X,*)$
  be a star-shaped transition system. We have the pushout diagram (in
  $\csts$):
\[
\xymatrix@C=3em {
  \{*\} \sqcup \{*\} \fr{} \fd{} & \{*\} \ar[dd]^-{\cyl_*(X)} \\
  X \sqcup X \fd{\gamma_X} & \\
  \cyl^\csts(X) \fr{} & \cocartesien
  \omega^*(\cyl_*(X,*)), }
\]
In $\cyl_*(X,*))$, every state $(\alpha,\epsilon)$ is reachable from
$(*,\epsilon)$ since $(X,*)$ is star-shaped by hypothesis.  Since
$(*,0)=(*,1)=*$ in $\omega^*(\cyl_*(X,*))$, every state of
$\cyl_*(X,*))$ is therefore reachable from $*$. Thus, the
Cattani-Sassone transition system $\cyl_*(X,*))$ is star-shaped.  \epf

We have $\cocyl^\cts(\{*\}) \iso \{(*,*)\}$.  By
\cite[Lemma~5.2]{leftdet}, the cylinder functor
$\cyl_*:\csts_* \to \csts_*$ has a right adjoint
$\cocyl_*:\csts_* \to \csts_*$ defined on objects by mapping
$(X,*)$ to the composite $\{*\}\stackrel{\iso}\to
\cocyl^\cts(\{*\})\to \cocyl^\cts(X)$ and on maps by mapping the
commutative triangle $\{*\}\to f$ to the composite
$\{*\}\stackrel{\iso}\to \cocyl^\cts(\{*\})\to \cocyl^\cts(f)$.  By
Theorem~\ref{star-shaped-csts}, the restriction of
$\cyl_*:\csts_* \to \csts_*$ to $\csts_{\bullet}$
gives rise to a very good cylinder $\cyl_\bullet:\csts_\bullet
\to \csts_\bullet$. A right adjoint is given by the composite functor
\[\cocyl_\bullet: \xymatrix@C=3em{\csts_\bullet \subset \csts_*
  \fr{\cocyl_*} & \csts_* \fr{} & \csts_{\bullet}}\] where the
right-hand map is the coreflection, i.e. the right adjoint of the
inclusion functor $\csts_{\bullet} \subset \csts_*$. In particular,
this means that the underlying Cattani-Sassone transition system of
$\cocyl_\bullet(X,*)$ is a subobject of the underlying Cattani-Sassone
transition system of $\cocyl_*(X,*)$, i.e. of $\cocyl^\cts(X)$
calculated in Table~\ref{csts-cyl-path}.

Table~\ref{csts-pointed-cyl-path} summarizes the computation of the
cylinder functor and of the path functor of $\csts_*$.

\begin{table}
\begin{tabular}{|c|c|c|}
\hline&&\\
 $(X,*)$ & $\cyl_*(X)$ & $\cocyl_*(X)$ \\&&\\
& $\Int(X,*) = \Int(X)\cup \{*\}$ & \\
& $\Ext(X,*) = \Ext(X)\backslash \{*\}$& \\&&\\
 $\St(X)$ & $\Int(X,*)\p \{0\} \sqcup \Ext(X,*)\p \{0,1\}$ & $\St(X) \p \St(X)$ \\&&\\
 $\Ac(X)$ & $\Ac(X)$ &  $\Ac(X)$ \\&&\\
 $\Tr(X)$ & $((\alpha,\epsilon_0),u_1,\dots,u_n,(\beta,\epsilon_{n+1}))$ & $((\alpha^-,\alpha^+),u_1,\dots,u_n,(\beta^-,\beta^+))$ \\
& such that $(\alpha,u_1,\dots,u_n,\beta) \in \Tr(X)$ & such that $(\alpha^\pm,u_1,\dots,u_n,\beta^\pm) \in \Tr(X)$\\&&\\
\hline
\end{tabular}\newline\newline
\caption{Cylinder functor and path functor of $\csts_*$}
\label{csts-pointed-cyl-path}
\end{table}

\section{Past-similarity of states}
\label{past-similar-sec}

\bd Let $(X,*)$ be a star-shaped transition system. Two states
$\alpha$ and $\beta$ of $X$ are {\rm past-similar} (denoted by $\alpha
\simeq_{past} \beta$) if there exists $w\in \Sigma^n$ with $n\geq 0$
and two right homotopic maps $(P(w),0)\rightrightarrows (X,*)$ sending
$n$ to $\alpha$ and $\beta$ respectively. \ed

\begin{figure}
\[
\xymatrix@C=3em@R=0.5em
{
& \bullet \fr{v}\ar@{->}[rdd]_-/-18pt/{v} & \alpha\\
{*} \ar@{->}[ru]^-{u} \ar@{->}[rd]_-{u}\\
& \bullet \fr{v}\ar@{->}[ruu]^-/-18pt/{v} & \beta \\
}
\]
\caption{$\alpha$ and $\beta$ are past-similar}
\label{past-similar}
\end{figure}
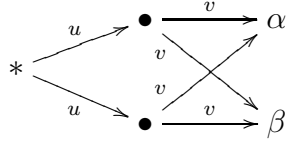

In Figure~\ref{past-similar}, the states $\alpha$ and $\beta$ are
past-similar because there is a homotopy between any path from $*$ to
$\alpha$ and any path from $*$ to $\beta$. In
Figure~\ref{not-past-similar}, the states $\alpha$ and $\beta$ are not
past-similar. However, they are past-similar in any fibrant
replacement.  In fact, the star-shaped transition system of
Figure~\ref{past-similar} is a fibrant replacement of the star-shaped
transition system of Figure~\ref{not-past-similar} because: 1) it is
fibrant by Theorem~\ref{fib-carac}, 2) the image by $R_\bullet^\perp$
of these two star-shaped transition systems is $*\stackrel{u}\to
\bullet \stackrel{v}\to \bullet\stackrel{w}\to \bullet$ (see
Section~\ref{desc-bullet}).

\begin{figure}
\[
\xymatrix@C=3em@R=0.5em
{
& \bullet \fr{v} & \alpha\\
{*} \ar@{->}[ru]^-{u} \ar@{->}[rd]_-{u}\\
& \bullet \fr{v} & \beta \\
}
\]
\caption{$\alpha$ and $\beta$ are not past-similar}
\label{not-past-similar}
\end{figure}
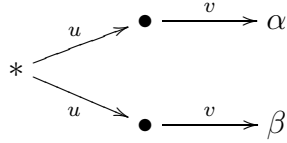

\bp \label{pastsimilar-pathsimilar} Let $(X,*)$ be a star-shaped
transition system. Let $\alpha$ and $\beta$ be two states of $X$.  We
have $\alpha \simeq_{past} \beta$ if and only if $(\alpha,\beta)$ is a
state of $\cocyl_\bullet(X,*)$. \ep

\bpf If $\alpha \simeq_{past} \beta$, then there exists $w\in
\Sigma^n$ for $n\geq 0$ and a map $(P(w),0) \to \cocyl_\bullet(X,*)$
taking $n$ to $(\alpha,\beta)$.  Thus, the pair $(\alpha,\beta)$ is a
state of $\cocyl_\bullet(X,*)$. Conversely, if $(\alpha,\beta)$ is a
state of $\cocyl_\bullet(X,*)$, then it is reachable. Consequently,
there exists a map $(P(w),0) \to \cocyl_\bullet(X,*)$ with $n\geq 0$
and $w\in \Sigma^n$ taking $n$ to $(\alpha,\beta)$.  \epf

For all actions $u$ of a Cattani-Sassone transition system $X$, there
exists a transition $(\alpha,u,\beta)$ of $X$ because $X$ is
cubical. Since $\alpha \simeq_{past} \alpha$ and $\beta \simeq_{past}
\beta$, the triple $((\alpha,\alpha),u,(\beta,\beta))$ is a transition
of $\cocyl_\bullet(X,*)$. Thus, all actions of $\Ac(X)$ are used in
$\cocyl_\bullet(X,*)$. Table~\ref{csts-star-shaped-cyl-path}
summarizes the computation of the cylinder functor and of the path
functor of $\csts_\bullet$.

\begin{table}
\begin{tabular}{|c|c|c|}
\hline&&\\
 $(X,*)$ & $\cyl_\bullet(X)$ & $\cocyl_\bullet(X)$ \\
&&\\
& $\Int(X,*) = \Int(X)\cup \{*\}$ & \\
& $\Ext(X,*) = \Ext(X)\backslash \{*\}$& \\
&&\\
 $\St(X)$ & $\Int(X,*)\p \{0\} \sqcup \Ext(X,*)\p \{0,1\}$ & $\simeq_{past}\subset \St(X) \p \St(X)$ \\
&&\\
 $\Ac(X)$ & $\Ac(X)$ &  $\Ac(X)$ \\
&&\\
 $\Tr(X)$ & $((\alpha,\epsilon_0),u_1,\dots,u_n,(\beta,\epsilon_{n+1}))$ & $((\alpha^-,\alpha^+),u_1,\dots,u_n,(\beta^-,\beta^+))$ \\
& such that $(\alpha,u_1,\dots,u_n,\beta) \in \Tr(X)$ & such that $(\alpha^\pm,u_1,\dots,u_n,\beta^\pm) \in \Tr(X)$\\
&&\\\hline
\end{tabular}\newline\newline
\caption{Cylinder functor and path functor of $\csts_\bullet$}
\label{csts-star-shaped-cyl-path}
\end{table}

\bp \label{eq-left-right} Let $f,g:(X,*) \rightrightarrows (Y,*)$ be
two maps of star-shaped transition systems. The following conditions
are equivalent:
\begin{enumerate}
\item $f$ and $g$ are left homotopic
\item $f$ and $g$ are right homotopic
\item $f$ and $g$ are homotopic.
\end{enumerate}
In particular, this means that 
\[\pi^l_{\csts_\bullet}((X,*),(Y,*)) =
\pi^r_{\csts_\bullet}((X,*),(Y,*)) = \pi_{\csts_\bullet}((X,*),(Y,*)).\]  \ep

Note that this proposition also holds in $\csts$ and in $\csts_*$.

\bpf We have $(3) \Rightarrow (1)$ and $(3) \Rightarrow (2)$ by
definition. It suffices to prove the equivalence $(1)\Leftrightarrow
(2)$. By adjunction, the existence of a left homotopy
$H_l:\cyl_\bullet(X,*) \to (Y,*)$ is equivalent to the existence of a
right homotopy $H_r:(X,*)\to \cocyl_\bullet(Y,*)$. The maps $f$ and
$g$ are left homotopic if and only if $H_l(\alpha,0)=f(\alpha)$,
$H_l(\alpha,1)=g(\alpha)$ for any state $\alpha$ of $X$ and
$H(u)=f(u)=g(u)$ for any action $u$ of $X$. The maps $f$ and $g$ are
right homotopic if and only if $H_r(\alpha)=(f(\alpha),g(\alpha))$ for
any state $\alpha$ of $X$ and $H_r(u)=u$ for any action $u$ of $X$.
\epf


\bd The set of transitions of a star-shaped transition system $(X,*)$
is {\rm closed under past-similarity} if for all $n \geq 1$, for all
transitions $(\alpha,u_1,\dots,u_n,\beta)$ of $X$, and for all states
$\gamma$ and $\delta$ of $X$, if $\alpha \simeq_{past} \gamma$ and
$\beta \simeq_{past} \delta$ then the tuple
$(\gamma,u_1,\dots,u_n,\delta)$ is a transition of $X$. \ed

\bp \label{fib-carac-demi} Let $(X,*)$ be a star-shaped transition
system. If $(X,*)$ is fibrant, then the set of transitions of $X$ is
closed under past-similarity. \ep

The converse is proved in Theorem~\ref{fib-carac}.

\bpf The proof is by induction on $n$.  

Suppose that $n=1$.  Let $(\alpha,u_1,\beta)$ be a transition of
$X$. Let $\gamma$ and $\delta$ be two states of $X$ such that $\alpha
\simeq_{past} \gamma$ and $\beta \simeq_{past} \delta$. Consider the
diagram of $\csts$
\[
\xymatrix@C=5em@R=4em
{
\{0_1,1_1\} \fd{\subset}  & \omega^*(\cocyl_\bullet(X,*)) \fd{\omega^*(\pi^0)} \\ 
C_1[\mu(u_1)] \ar@{-->}[ru]^-{\ell}\ar@{->}[r]_-{\text{\tiny {\(\begin{array}{c}0_1 \mapsto \alpha \\1_1 \mapsto \beta \\(\mu(u_1),1) \mapsto u_1\end{array}\)}}} & \omega^*((X,*)).
}
\]
By hypothesis, we have $\alpha \simeq_{past} \gamma$ and $\beta
\simeq_{past} \delta$. By Proposition~\ref{pastsimilar-pathsimilar},
$(\alpha,\gamma)$ and $(\beta,\delta)$ are two states of
$\cocyl_\bullet(X,*)$. We obtain the commutative diagram of solid
arrows of $\csts$
\[
\xymatrix@C=5em@R=4em
{
\{0_1,1_1\} \fd{\subset} \fr{\text{\tiny {\(\begin{array}{c}0_1 \mapsto (\alpha,\gamma) \\1_1 \mapsto (\beta,\delta)\end{array}\)}}} &\omega^*(\cocyl_\bullet(X,*)) \fd{\omega^*(\pi^0)} \\ 
C_1[\mu(u_1)] \ar@{-->}[ru]^-{\ell}\ar@{->}[r]_-{\text{\tiny {\(\begin{array}{c}0_1 \mapsto \alpha \\1_1 \mapsto \beta \\(\mu(u_1),1) \mapsto u_1\end{array}\)}}} & \omega^*((X,*)).
}
\]
The map $\rho^*(\{0_1,1_1\} \subset C_1[\mu(u_1)])\varoslash \pi^0$
is a cofibration by Theorem~\ref{oslash}.  Since $(X,*)$ is fibrant,
it is injective with respect to the trivial cofibration
\[(\rho^*(\{0_1,1_1\} \subset C_1[\mu(u_1)])\varoslash \pi^0)\star
\gamma^0.\] Thus, the right vertical map satisfies the RLP with
respect to $(\{0_1,1_1\} \subset C_1[\mu(u_1)])\varoslash \pi^0$, and
then with respect to the left vertical map $\{0_1,1_1\} \subset
C_1[\mu(u_1)]$. Since $\pi^0$ is the identity on actions, we obtain
that the triple $((\alpha,\gamma),u_1,(\beta,\delta))$ is a transition
of $\omega^*(\cocyl_\bullet(X,*))$. By
Table~\ref{csts-star-shaped-cyl-path}, we deduce that
$(\gamma,u_1,\delta)$ is a transition of $X$.  

Suppose the induction hypothesis proved for all $p \leq n$ with $n\geq
1$. Choose a transition $(\alpha,u_1,\dots,u_{n+1},\beta)$ of $X$,
which gives rise to a map $f:C_{n+1}[\mu(u_1),\dots,\mu(u_{n+1})] \to
\omega^*(X,*)$. Let $\gamma$ and $\delta$ be two states of $X$ such
that $\alpha \simeq_{past} \gamma$ and $\beta \simeq_{past} \delta$.
By induction hypothesis, there is a commutative diagram of $\csts$
\[
\xymatrix@C=15em@R=5em
{
\de C_{n+1}[\mu(u_1),\dots,\mu(u_{n+1})] \fd{\subset} \ar@{->}[r]^-{\text{\tiny {\(\begin{array}{c}0_{n+1} \mapsto (\alpha,\gamma) \\f:\zeta \in \{0,1\}^n\backslash \{0_{n+1},1_{n+1}) \mapsto (f(\zeta),f(\zeta))\\1_{n+1} \mapsto (\beta,\delta)\\ (\mu(u_i),i) \mapsto u_i \hbox{ for }1\leq i \leq n+1\end{array}\)}}} & \omega^*(\cocyl_\bullet(X,*)) \fd{\omega^*(\pi^0)} \\ 
C_{n+1}[\mu(u_1),\dots,\mu(u_{n+1})] \ar@{-->}[ru]^-{\ell}\ar@{->}[r]^-{f} & \omega^*(X,*).
}
\]
The map $\rho^*(\de C_{n+1}[\mu(u_1),\dots,\mu(u_{n+1})] \subset
C_{n+1}[\mu(u_1),\dots,\mu(u_{n+1})]) \varoslash \pi^0$ is a
cofibration by Theorem~\ref{oslash}. Thus, $(X,*)$ is injective with
respect to the trivial cofibration
\[\left(\rho^*(\de C_{n+1}[\mu(u_1),\dots,\mu(u_{n+1})] \subset C_{n+1}[\mu(u_1),\dots,\mu(u_{n+1})]) \varoslash \pi^0\right) \star \gamma^0.\]
Hence, the lift $\ell$ exists. Since $\pi^0$ is the identity on
actions, we obtain that the tuple
$((\alpha,\gamma),u_1,\dots,u_{n+1},(\beta,\delta))$ is a transition
of $\omega^*(\cocyl_\bullet(X,*))$. By
Table~\ref{csts-star-shaped-cyl-path}, we deduce that
$(\gamma,u_1,\dots,u_{n+1},\delta)$ is a transition of $X$.  \epf

\bp Let $(X,*)$ be a star-shaped transition system. Past-similarity is
a reflexive and symmetric relation. There exists a star-shaped
transition system $(X,*)$ such that the binary relation
$\simeq_{past}$ is not transitive. Let $f:(X,*)\to (Y,*)$ be a map of
star-shaped transition systems. If $\alpha$ and $\beta$ are two
past-similar states of $X$, then $f(\alpha)$ and $f(\beta)$ are two
past-similar states of $Y$. \ep

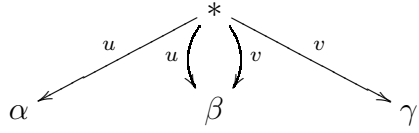
\begin{figure}
\[
\xymatrix@C=5em@R=2em
{
& {*} \ar@{->}[dl]_-{u}\ar@{->}[dr]^-{v}  \ar@/_10pt/@{->}[d]_-{u}\ar@/^10pt/@{->}[d]^-{v} &\\
\alpha & \beta & \gamma
}
\]
\caption{Example of non transitive past-similarity}
\label{not-trans}
\end{figure}

\bpf Reflexivity and symmetry are obvious.  Consider the following
example of Figure~\ref{not-trans} (with $\mu(u) \neq \mu(v)$).  We
have $\alpha \simeq_{past} \beta$ and $\beta \simeq_{past}
\gamma$. But $\alpha$ is not past-similar to $\gamma$.  If $\alpha$
and $\beta$ are two past-similar states of $X$, then there exists a
right homotopy $H:(P(w),0) \to \cocyl_\bullet(X,*)$ for $w\in
\Sigma^n$ sending $n$ to $(\alpha,\beta)$. The map $f$ induces a map
$\cocyl_\bullet(X,*) \to \cocyl_\bullet(Y,*)$ by functoriality, and
therefore a right homotopy $(P(w),0) \to \cocyl_\bullet(X,*) \to
\cocyl_\bullet(Y,*)$ sending $n$ to $(f(\alpha),f(\beta))$. Thus,
$f(\alpha)$ and $f(\beta)$ are past-similar.  \epf

\bp \label{pastsim-transitive} Let $(X,*)$ be a star-shaped transition
system such that the set of transitions is closed under
past-similarity. Then $\simeq_{past}$ is transitive. In particular, if
$(X,*)$ is fibrant, then $\simeq_{past}$ is transitive. \ep

\bpf Let $\alpha,\beta,\gamma$ be three states of $X$ with
$\alpha \simeq_{past} \beta$ and $\beta\simeq_{past} \gamma$. Then
there exists a map $f:\cyl_\bullet(P(w),0) \to (X,*)$ with
$w\in \Sigma^n$ and $n\geq 0$ such that $f(n,0)=\alpha$ and
$f(n,1)=\beta$. If $n=0$, then $\alpha=\beta$ and therefore
$\alpha \simeq_{past} \gamma$. Let us suppose $n\geq 1$.  Define
$g:\cyl_\bullet(P(w),0) \to (X,*)$ on states by $g(n,1)=\gamma$ and
$g=f$ otherwise, and $g=f$ on actions. The only transitions of
$\cyl_\bullet(P(w),0)$ involving the state $(n,1)$ are
$((n-1,0),(x_n,n),(n,1))$ and $((n-1,1),(x_n,n),(n,1))$ where
$w=x_1\dots x_n$. Since $f$ is a well-defined map of star-shaped
transition systems, the triples
$(f(n-1,0),f(x_n,n),f((n,1)))=(f(n-1,0),f(x_n,n),\beta)$ and
$(f(n-1,1),f(x_n,n),f((n,1)))=(f(n-1,1),f(x_n,n),\beta)$ are two
$1$-transitions of $X$.  Since $\beta\simeq_{past} \gamma$ and since
the set of transitions of $(X,*)$ is closed under past-similarity, the
triples $(g(n-1,0),g(x_n,n),g((n,1)))=(f(n-1,0),f(x_n,n),\gamma)$ and
$(g(n-1,1),g(x_n,n),g((n,1)))=(f(n-1,1),f(x_n,n),\gamma)$ are two
$1$-transitions of $(X,*)$. We deduce that the map
$g:\cyl_\bullet(P(w),0) \to (X,*)$ is a well-defined map of
star-shaped transition systems. Thus, $\alpha \simeq_{past}
\gamma$.
The last sentence is a consequence of
Proposition~\ref{fib-carac-demi}.  \epf

\bth \label{car_homeq-2} Let $(X,*)$ and $(Y,*)$ be two objects of
$\csts_\bullet$ with $(Y,*)$ fibrant. Two maps from $(X,*)$ to $(Y,*)$
are homotopy equivalent if and only if they coincide on actions and
send any state of $(X,*)$ to past-similar states of $(Y,*)$.  \eth

\bpf Let $f$ and $g$ be two homotopy equivalent maps from $(X,*)$ to
$(Y,*)$.  There exists a right homotopy $H:(X,*) \to
\cocyl_\bullet(Y,*)$ from $f$ to $g$.  Let $u$ be an action of
$X$. Since $X$ is cubical, there exists a transition
$(\alpha,u,\beta)$. We deduce that
\[H(\alpha,u,\beta) = ((f(\alpha),g(\alpha)),u',(f(\beta),g(\beta)))\]
with $(f(\alpha),u',f(\beta))=f(\alpha,u,\beta)$ and
$(g(\alpha),u',g(\beta))=g(\alpha,u,\beta)$. We obtain $u'=f(u)=g(u)$.
Let $\alpha$ be a state of $X$. Then $H(\alpha) =
(f(\alpha),g(\alpha))$. We deduce that $f(\alpha)\simeq_{past}
g(\alpha)$ by Proposition~\ref{pastsimilar-pathsimilar}.

Conversely, let $f$ and $g$ be two maps from $(X,*)$ to $(Y,*)$ which
coincide on actions and such that for any state $\alpha$ of $X$,
$f(\alpha)$ and $g(\alpha)$ are past-similar.  We have to construct a
right homotopy $H:(X,*) \to \cocyl_\bullet(Y,*)$ from $f$ to $g$. The
only possible definition is $H(\alpha)=(f(\alpha),g(\alpha))$ for all
states $\alpha$ of $(X,*)$ and $H(u)=f(u)=g(u)$ for all actions $u$ of
$(X,*)$. We have to prove that $H$ yields a well-defined map of
transition systems, i.e. that for all transitions
$(\alpha,u_1,\dots,u_n,\beta)$ of $(X,*)$, the tuple
$((f(\alpha),g(\alpha)),H(u_1),\dots,H(u_n),(f(\beta),g(\beta)))$ is a
transition of $\cocyl_\bullet(Y,*)$. We are going to prove by
induction on $n\geq 1$ that for all states $\alpha$ and $\beta$ of
$(X,*)$ and all actions $u_1,\dots,u_p$ of $X$ with $1\leq p \leq n$,
the tuple
$((f(\alpha),g(\alpha)),H(u_1),\dots,H(u_p),(f(\beta),g(\beta)))$ is a
transition of $\cocyl_\bullet(Y,*)$ if $(\alpha,u_1,\dots,u_p,\beta)$
is a transition of $(X,*)$.

\underline{The case $n=1$}. Consider the commutative diagram of solid
arrows of $\csts$
\[
\xymatrix@C=7em@R=3em
{
\{0_1,1_1\} \fd{\subset} \fr{\text{\tiny {\(\begin{array}{c}0_1 \mapsto (f(\alpha),g(\alpha)) \\1_1 \mapsto (f(\beta),g(\beta))\end{array}\)}}} & \omega^*(\cocyl_\bullet(Y,*)) \fd{\omega^*(\pi^0)} \\ 
C_1[\mu(u_1)] \ar@{-->}[ru]^-{\ell}\ar@{->}[r]_-{\text{\tiny {\(\begin{array}{c}0_1 \mapsto f(\alpha) \\1_1 \mapsto g(\alpha)\\(\mu(u_1),1) \mapsto f(u_1)\end{array}\)}}} & \omega^*((Y,*)).
}
\]
Since $(Y,*)$ is fibrant, the canonical map $(Y,*)\to \mathbf{1}$ satisfies the
RLP with respect to the trivial cofibration $(\rho^*(\de C_1[\mu(u_1)]
\subset C_1[\mu(u_1)])\varoslash \pi^0) \star \gamma^0$. By adjunction, this implies
that the lift $\ell$ in the diagram above exists. Thus, we have proved
that the triple
\[\ell(0_1,\mu(u_1),1_1))=((f(\alpha),g(\alpha)),f(u_1),(f(\beta),g(\beta)))\]
is a transition of $\cocyl_\bullet(Y,*)$.

\underline{From $n$ to $n+1$ with $n\geq 1$}. By induction hypothesis, we have the
commutative diagram of solid arrows of $\csts$:
\[
\xymatrix@C=12em@R=3em
{
\de C_{n+1}[\mu(u_1),\dots,\mu(u_{n+1})] \fd{\subset} \ar@{->}[r]^-{\text{\tiny {\(\begin{array}{c}0_{n+1} \mapsto (f(\alpha),g(\alpha)) \\1_{n+1} \mapsto (f(\beta),g(\beta))\\ (\mu(u_i),i) \mapsto f(u_i) \hbox{ for }1\leq i \leq n+1\end{array}\)}}} & \omega^*(\cocyl_\bullet(Y,*)) \fd{\omega^*(\pi^0)} \\ 
C_{n+1}[\mu(u_1),\dots,\mu(u_{n+1})] \ar@{-->}[ru]^-{\ell}\ar@{->}[r]_-{\text{\tiny {\(\begin{array}{c}0_{n+1} \mapsto f(\alpha) \\1_{n+1} \mapsto g(\alpha)\\(\mu(u_i),i) \mapsto f(u_i) \hbox{ for }1\leq i \leq n+1\end{array}\)}}} & \omega^*(Y,*).
}
\]
Since $(Y,*)$ is fibrant, the canonical map $(Y,*)\to \mathbf{1}$
satisfies the RLP with respect to the trivial
cofibration \[\left(\rho^*\left(\de
    C_{n+1}[\mu(u_1),\dots,\mu(u_{n+1})] \subset
    C_{n+1}[\mu(u_1),\dots,\mu(u_{n+1})]\right)\varoslash
  \pi^0\right)\star \gamma^0.\] By adjunction, this implies that the
lift $\ell$ in the diagram above exists.  Thus, we have proved that
the triple
\[\ell(0_{n+1},\mu(u_1),\dots,\mu(u_{n+1}),1_{n+1})=((f(\alpha),g(\alpha)),f(u_1),\dots,f(u_{n+1}),(f(\beta),g(\beta)))\]
is a transition of $\cocyl_\bullet(Y,*)$.  \epf

\begin{cor} \label{hom-eq11} Every weak equivalence between two
  fibrant star-shaped transition systems is bijective on actions.
\end{cor}

\section{Fibrant star-shaped transition systems}
\label{fibrant-section}

\begin{lem} \label{aide} Consider a commutative square of $\csts_\bullet$
\[
\xymatrix@C=2em@R=2em
{
(A,*) \fr{} \fd{} & (C,*) \fd{}\\
(B,*) \fr{} & (D,*).
}
\]
Then the set map $\St(B,*) \sqcup_{\St(A,*)} \St(C,*) \to \St(D,*)$
factors as a composite
\[\St(B,*) \sqcup_{\St(A,*)} \St(C,*) \twoheadrightarrow \St((B,*)
\sqcup_{(A,*)} (C,*)) \longrightarrow \St(D,*)\] where the left-hand
map is onto and the set map $\Ac(B,*) \sqcup_{\Ac(A,*)} \Ac(C,*) \to
\Ac(D,*)$ factors as a composite
\[\Ac(B,*) \sqcup_{\Ac(A,*)} \Ac(C,*) \twoheadrightarrow \Ac((B,*)
\sqcup_{(A,*)} (C,*)) \longrightarrow \Ac(D,*)\] where the left-hand map is onto.
\end{lem}

\bpf A pushout $(B,*)\sqcup_{(A,*)} (C,*)$ in $\csts_\bullet$ can be
calculated in $\csts_*$ since the category $\csts_\bullet$ is a
coreflective subcategory of $\csts_*$. Since a pushout is a connected
 colimit, we have the isomorphism 
 \[(B,*)\sqcup_{(A,*)} (C,*) \iso (B \sqcup_A C,*).\] A pushout $B
 \sqcup_A C$ in $\csts$ is calculated by taking the pushout $B
 \sqcup^\cts_A C$ in $\cts$, and then by applying the reflection
 $\CSA_1\CSA_2:\cts\to \csts$ which identifies states and actions to
 force CSA1 and CSA2 to hold.  We obtain the isomorphism of
 $\csts_\bullet$:
 \[(B \sqcup_A C,*) \iso (\CSA_1\CSA_2(B \sqcup^\cts_A C),*).\]
 By adjunction, the map of cubical transition systems
 \[B \sqcup^\cts_A C \to D\]
 factors uniquely as a composite
 \[B \sqcup^\cts_A C \longrightarrow B \sqcup_A C\longrightarrow D.\]
 By Lemma~\ref{colim-preserv}, we have the bijections of sets
 \[\St(B \sqcup^\cts_A C) \iso \St(B) \sqcup_{\St(A)} \St(C) \hbox{
   and }\Ac(B \sqcup^\cts_A C) \iso \Ac(B) \sqcup_{\Ac(A)} \Ac(C).\]
 We obtain the composite set map
 \[\St(B) \sqcup_{\St(A)} \St(C) \twoheadrightarrow \St((B,*)
 \sqcup_{(A,*)} (C,*)) \longrightarrow \St(D)\]
 where the left-hand map is onto and the composite set map
 \[\Ac(B) \sqcup_{\Ac(A)} \Ac(C) \twoheadrightarrow \Ac((B,*)
 \sqcup_{(A,*)} (C,*)) \longrightarrow \Ac(D)\]
 where the left-hand map is onto.  \epf

 \begin{nota} 
   Let \begin{enumerate}
\item $\Lambda_0(Z) = (Z\star \gamma^0) \cup (Z\star \gamma^1)$
\item 
     $\Lambda_{n+1}(Z)=\Lambda_{n}(Z) \star \gamma$ for $n\geq 0$
\item
     $\Lambda(Z) := \bigcup_{n\geq 0} \Lambda_n(Z)$
   \end{enumerate} where
   $\gamma^0,\gamma^1:\id \Rightarrow \cyl_\bullet$ are the natural
   transformations associated with the cylinder
   $\cyl_\bullet:\csts_\bullet \to \csts_\bullet$ and
   $\gamma = \gamma^0 \sqcup \gamma^1$ (see
   Table~\ref{csts-star-shaped-cyl-path} page
   \pageref{csts-star-shaped-cyl-path}) and where $Z$ is a set of maps
   of $\csts_\bullet$.
\end{nota}

\bp \label{ex0} Every map of $\Lambda(\I_\bullet)$ is bijective on
actions.  \ep

\bpf Let $f:(A,*) \to (B,*)$ be a map of star-shaped transition
systems.  Using Lemma~\ref{aide}, the map $f\star \gamma^\epsilon$
gives rise to the composite set map
\[\Ac(A,*) \sqcup_{\Ac(A,*)} \Ac(B,*) \twoheadrightarrow \Ac\left(\cyl_\bullet(A,*) \sqcup_{(A,*)} (B,*)\right)
\longrightarrow \Ac(B,*)\] which is bijective. Thus, the left-hand map
is injective and then bijective. This implies that the right-hand map
is bijective and that $f\star \gamma^\epsilon$ is bijective on
actions.  Suppose now that $f$ is bijective on actions. By
Lemma~\ref{aide}, the map $f\star \gamma$ gives rise to the composite
set map
\begin{multline*}\Ac(\cyl_\bullet(A,*)) \sqcup_{\Ac((A,*)\sqcup
    (A,*))} \Ac((B,*)\sqcup (B,*))\\\twoheadrightarrow
  \Ac\left(\cyl_\bullet(A,*) \sqcup_{(A,*)\sqcup (A,*)} ((B,*)\sqcup
    (B,*))\right) \longrightarrow
  \Ac(\cyl_\bullet(B,*))\end{multline*} which is bijective. Thus, the
left-hand map is injective and then bijective. This implies that the
right-hand map is bijective and that $f\star \gamma$ is bijective on
actions.  \epf

\bp \label{int-ext-cyl} Let $(X,*)$ be a star-shaped transition
system. Then we have the equalities
$\Int(\cyl_\bullet(X,*))=\Int(X,*)\p \{0\}$ and
$\Ext(\cyl_\bullet(X,*))=\Ext(X,*)\p \{0,1\}$. \ep

\bpf The natural map $\gamma^0_X: (X,*) \to \cyl_\bullet(X,*)$
induces a one-to-one set map \[\St(\gamma^0_X):\Int(X,*) \sqcup \Ext(X,*) \longrightarrow 
\Int(X,*) \p \{0\} \sqcup \Ext(X,*)\p \{0,1\}.\]  We deduce the set inclusion
\[\Int(\gamma^0_X): \Int(X,*) \iso \Int(X,*)\p \{0\} \subset \Int(\cyl_\bullet(X,*)).\] 
Let $(\alpha,\epsilon)\in \Int(\cyl_\bullet(X,*))$. By
Proposition~\ref{int_charac}, there exists a $2$-transition
\[((\mu_0,\epsilon_0),u_1,u_2,(\mu_1,\epsilon_1))\] of
$\cyl_\bullet(X,*)$ such that the tuples
$((\mu_0,\epsilon_0),u_1,(\alpha,\epsilon))$ and
$((\alpha,\epsilon),u_2,(\mu_1,\epsilon_1))$ are two transitions of
$\cyl_\bullet(X,*)$. Using Table~\ref{csts-star-shaped-cyl-path}, we
deduce that the tuples $(\mu_0,u_1,u_2,\mu_1)$, $(\mu_0,u_1,\alpha)$
and $(\alpha,u_2,\mu_1)$ are transitions of $(X,*)$.  This implies
that $\alpha$ is an internal state of $(X,*)$ and that
$\epsilon=0$. We deduce the set inclusion
\[\Int(\cyl_\bullet(X,*)) \subset
\Int(X,*)\p \{0\}.\] We obtain the equality \[\Int(\cyl_\bullet(X,*))
= \Int(X,*)\p \{0\} \] and, by taking the complement, the equality
 \[\Ext(\cyl_\bullet(X,*)) = \Ext(X,*)\p \{0,1\}.\]
\epf

\bd A map $f:(A,*)\to (B,*)$ of pointed transition systems is {\rm
  proper} if the set map $\Int(f):\Int(A,*) \to \Int(B,*)$ is
bijective and if $f$ induces a one-to-one set map $f(\Ext(A,*))\subset
\Ext(B,*)$ on external states. In particular, a proper map is
one-to-one on states. A map $f:(A,*)\to (B,*)$ of pointed transition
systems is {\rm strongly proper} if it is proper and bijective on
states (which means that not only it is bijective on internal states,
but also on external states).\ed

\bp \label{oslash-proper} Let $f:(A,*)\to (B,*)$ be a (strongly resp.)
proper map of pointed transition systems. Then for any $g\in
\Mor(\csts_\bullet)$, the map $f \varoslash g$ is (strongly resp.)
proper.  \ep

\bpf Let $f \varoslash g:(\widehat{A},*) \to (\widehat{B},*)$. The
star-shaped transition system $(\widehat{A},*)$ is obtained from
$(A,*)$ by doing for each state $\alpha$ of $(A,*)$ different from the
base state one of the following operations (cf. the proof of
Theorem~\ref{oslash}): 1) either identifying $\alpha$ with $*$, 2) or
attaching a finite sequence of $1$-transitions with initial state $*$
and final state $\alpha$. When $\alpha$ is internal in $(A,*)$, it
remains internal in $(\widehat{A},*)$. When it is external in $(A,*)$,
it becomes internal (case 1) and actually equal to $*$, or it remains
external (case 2). In case 2, all intermediate states of the finite
sequence of $1$-transitions are external. The star-shaped transition
system $(\widehat{B},*)$ is obtained from $(B,*)$ by doing for each
state $\alpha$ of $(B,*)$ different from the base state one of the
following operations (cf. the proof of Theorem~\ref{oslash}): 1)
nothing for a state of $f(\St(A))$, 2) one of the two operations above
on the states of $\St(B)\backslash f(\St(A))$. Thus, $f \varoslash g$
is proper.  \epf

\begin{nota} Let $\Sigma^+ = \bigcup_{n\geq 1} \Sigma^n$ denote the
  set of nonempty words over $\Sigma$. \end{nota}

\begin{cor} \label{ex1} Every map of $\I_\bullet$ is proper. Every map
  of $\I_\bullet \backslash \{(\varnothing \to (P(w),0)) \mid w\in
  \Sigma^+\}$ is strongly proper (where $\varnothing=(\{*\},*)$ is the
  initial object of $\csts_\bullet$).
\end{cor}

\bpf For any set $S$, we have $\Int(S)=\varnothing$, $\Ext(S)=S$. For
any $n$-cube $C_n[x_1,\dots,x_n]$ with $n\geq 1$ and $x_1,\dots,x_n\in
\Sigma$, we have $\Ext(C_n[x_1,\dots,x_n]) = \{0_n,1_n\}$ and
$\Int(C_n[x_1,\dots,x_n]) = \{0,1\}^n \backslash \{0_n,1_n\}$. Thus,
every map of $\rho^*(\I^\cts)$ is proper and every map of
$\rho^*(\I^\cts \backslash \{\varnothing \subset \{0\}\})$ is strongly
proper. The proof of the first statement is complete using
Proposition~\ref{oslash-proper}. We have \[\{\rho^*(\varnothing
\subset \{0\})\varoslash g \mid g\in \Mor(\csts_\bullet)\} =
\{(\varnothing \to (P(w),0)) \mid w\in \Sigma^+\}.\] Then the proof of
the last part is complete using Proposition~\ref{oslash-proper}.  \epf

\bp \label{ex2} Let $f:(A,*)\to (B,*)$ be a (strongly resp.) proper
map of star-shaped transition systems. Then the map $f \star
\gamma^\epsilon$ for $\epsilon=0,1$ is (strongly resp.) proper. 
\ep

\bpf We write the proof for $\epsilon = 0$.  We have
\begin{align*}
&\St(\cyl_\bullet(A,*)) \sqcup_{\St(A,*)} \St(B,*) \\
&\iso \left(\Int(A,*)\p\{0\} \sqcup \Ext(A,*)\p\{0,1\}\right) \sqcup_{\Int(A,*)\p \{0\} \sqcup \Ext(A,*)\p\{0\}} \left(\Int(B,*)\p\{0\}\sqcup \Ext(B,*)\p\{0\}\right) \\
& \iso \Int(B,*)\p\{0\} \sqcup \Ext(B,*)\p\{0\} \sqcup \Ext(A,*)\p\{1\}.
\end{align*}
The map $f\star \gamma^0$ gives rise by Lemma~\ref{aide} to the
composite set map
\begin{multline*}\Int(B,*)\p\{0\} \sqcup \Ext(B,*)\p\{0\} \sqcup \Ext(A,*)\p\{1\} \\\twoheadrightarrow \St\left(\cyl_\bullet(A,*)
    \sqcup_{(A,*)} (B,*)\right) \longrightarrow \Int(B,*)\p\{0\}
  \sqcup \Ext(B,*)\p\{0,1\}\end{multline*} which is one-to-one. Thus,
the left-hand map is one-to-one, and then bijective.  We deduce the
bijection of sets
\[\St\left(\cyl_\bullet(A,*)
  \sqcup_{(A,*)} (B,*)\right) \iso \Int(B,*)\p\{0\} \sqcup
\Ext(B,*)\p\{0\} \sqcup \Ext(A,*)\p\{1\}\] and therefore, that $f\star
\gamma^0$ is one-to-one on states. The map of star-shaped
transition systems \[\iota_2:(B,*) \to
\cyl_\bullet(A,*) \sqcup_{(A,*)} (B,*)\] induces on internal states 
the set map 
\[\Int(\iota_2):\Int(B,*) \longrightarrow 
\Int(\cyl_\bullet(A,*) \sqcup_{(A,*)} (B,*)).\] The latter yields the
inclusion of sets \[\Int(B,*) \p \{0\} \subset \Int(\cyl_\bullet(A,*)
\sqcup_{(A,*)} (B,*)).\] The map $f \star \gamma^0:\cyl_\bullet(A,*)
\sqcup_{(A,*)} (B,*) \to \cyl_\bullet(B,*)$ induces a set inclusion
between the sets of internal states
\[\Int\left(\cyl_\bullet(A,*) \sqcup_{(A,*)}
  (B,*)\right) \subset \Int(\cyl_\bullet(B,*)).\] By
Proposition~\ref{int-ext-cyl}, we have
$\Int(\cyl_\bullet(B,*))=\Int(B,*) \p\{0\}$.  We obtain the set
inclusion
\[\Int\left(\cyl_\bullet(A,*) \sqcup_{(A,*)}
  (B,*)\right) \subset \Int(B,*) \p\{0\}.\] Thus, we obtain the
equality
\[\Int(\cyl_\bullet(A,*) \sqcup_{(A,*)} (B,*)) = \Int(B,*) \p \{0\}\]
and, by taking the complement, the equality 
\[\Ext(\cyl_\bullet(A,*) \sqcup_{(A,*)} (B,*)) = \Ext(A,*)
\p \{1\} \sqcup \Ext(B,*)\p \{0\}.\] 
The proof is complete thanks to Proposition~\ref{int-ext-cyl}.
\epf 

\bp \label{ex3} Let $f:(A,*)\to (B,*)$ be a proper map. The map
$f\star \gamma$ is strongly proper. \ep

\bpf We have
\begin{align*}
&\St(\cyl_\bullet(A,*)) \sqcup_{\St((A,*)\sqcup(A,*))} \St((B,*)\sqcup (B,*)) \\
&\iso (\Int(A,*) \p \{0\} \sqcup \Ext(A,*) \p\{0,1\} )
\\ &\hspace{3cm}\sqcup_{\Int(A,*)\p \{0,1\} \sqcup \Ext(A,*)\p \{0,1\}} \left(\Int(B,*)\p \{0,1\} \sqcup \Ext(B,*)\p \{0,1\}\right)\\
&\iso \Int(A,*) \p \{0\}\sqcup \Ext(B,*)\p \{0,1\}.
\end{align*}
Thus the map $f\star \gamma$ gives rise by Lemma~\ref{aide} to the composite set map 
\begin{multline*}
\Int(A,*) \p \{0\}\sqcup \Ext(B,*)\p \{0,1\} \twoheadrightarrow \St(\cyl_\bullet(A,*)
\sqcup_{(A,*)\sqcup (A,*)} ((B,*)\sqcup (B,*))) \\\longrightarrow \Int(B,*) \p \{0\}\sqcup \Ext(B,*)\p \{0,1\}.
\end{multline*}
Since this composite is bijective, the left-hand map is injective and
then bijective. We obtain
\[\St(\cyl_\bullet(A,*)
\sqcup_{(A,*)\sqcup (A,*)} ((B,*)\sqcup (B,*))) = \Int(A,*)\p \{0\}
\sqcup \Ext(B,*)\p \{0,1\}\] and $f\star \gamma$ is bijective on
states. The states of $\Ext(B,*)\p \{0,1\}$ are external in
$\cyl_\bullet(B,*)$ by Proposition~\ref{int-ext-cyl}. Thus, the states
$\Ext(B,*)\p \{0,1\}$ cannot be internal states of $\cyl_\bullet(A,*)
\sqcup_{(A,*)\sqcup (A,*)} ((B,*)\sqcup (B,*))$. We have proved the
inclusion \[\Int(\cyl_\bullet(A,*) \sqcup_{(A,*)\sqcup (A,*)}
((B,*)\sqcup (B,*))) \subset \Int(B,*)\p \{0\}.\] The map
$\cyl_\bullet(A,*) \to \cyl_\bullet(A,*) \sqcup_{(A,*)\sqcup (A,*)}
((B,*)\sqcup (B,*))$ induces a set map \[\Int(\cyl_\bullet(A,*)) \to
\Int(\cyl_\bullet(A,*) \sqcup_{(A,*)\sqcup (A,*)} ((B,*)\sqcup
(B,*))).\] By Proposition~\ref{int-ext-cyl}, we have
$\Int(\cyl_\bullet(A,*))=\Int(A,*) \p \{0\} \iso \Int(B,*)\p
\{0\}$. We deduce the inclusion of sets \[\Int(B,*)\p \{0\} \subset
\Int(\cyl_\bullet(A,*) \sqcup_{(A,*)\sqcup (A,*)} ((B,*)\sqcup
(B,*))).\] The proof is complete thanks to
Proposition~\ref{int-ext-cyl}.  \epf

\begin{cor} \label{ex4} Every map of $\Lambda(\I_\bullet) \backslash
  \Lambda_0(\{\varnothing \to (P(w),0)\mid w \in \Sigma^+\})$ is strongly
  proper, and in particular, bijective on states.
\end{cor}

\bpf Table~\ref{properness} summarizes the situation. \epf

\begin{table}
{\small
\begin{tabular}{|c|c|c|}
\hline&&\\
 $\I^\cts$ & $\varnothing \subset \{0\}$ & $\I^\cts\backslash \{\varnothing \subset \{0\}\}$\\
&&\\
 $\I_\bullet=\{\rho^*(f) \varoslash g\mid f\in \I^\cts$& $\{\varnothing \to (P(w),0)\mid w \in \Sigma^+\}$&$\I_\bullet\backslash\{\varnothing \to
(P(w),0)\mid w \in \Sigma^+\}$ \\ 
$ \hbox{ and }g\in \Mor(\csts_\bullet)\}$& proper (\ref{ex1}) & strongly proper (\ref{ex1})\\
&&\\
 $\Lambda_0$ & proper (\ref{ex2}) & strongly proper (\ref{ex2}) \\
&&\\
 $\Lambda_n$ with $n\geq 1$ & strongly proper (\ref{ex3}) &  strongly proper(\ref{ex3})  \\&&\\
\hline
\end{tabular}
}\newline\newline
\caption{Properness of the maps of $\Lambda(\I_\bullet)$}
\label{properness}
\end{table}

\bth \label{fib-carac} Let $(X,*)$ be a star-shaped transition
system. Then $(X,*)$ is fibrant if and only if its set of transitions
is closed under past-similarity. \eth

\bpf The ``only if'' part is Proposition~\ref{fib-carac-demi}. Let
$(X,*)$ be a star-shaped transition system such that the set of
transition is closed under past-similarity.  We have to prove that
$(X,*)$ is injective with respect to any map of $\Lambda(\I_\bullet)$.

Let us prove first that $(X,*)$ is injective with respect to the maps
of $\Lambda_0(\{\varnothing \to (P(w),0)\mid w \in \Sigma^+\})$.  Let
$i:\varnothing \to (P(w),0)$ with $w \in \Sigma^+$. By symmetry, it
suffices to prove the injectivity with respect to $i\star
\gamma^0$. Consider the diagram of solid arrows of $\csts_\bullet$
\[
\xymatrix@C=3em@R=5em {
  (P(w),0) \fr{\phi} \fd{i\star \gamma^0} & (X,*) \\
  \cyl_\bullet(P(w),0) \ar@{-->}[ru]^-{\ell} && }
\]
The map $i\star \gamma^0:(P(w),0) \to \cyl_\bullet(P(w),0)$ has the
retraction $\sigma:\cyl_\bullet(P(w),0) \to (P(w),0)$. Thus,
$\ell=\phi\sigma$ is a solution. So far, we have not used the fact
that the set of transitions of $(X,*)$ is closed under
past-similarity.

Let us prove now that $(X,*)$ is injective with respect to any map
of $\Lambda_0(\I_\bullet \backslash\{\varnothing \subset (P(w),0)\mid
w \in \Sigma^+\})$.  By symmetry, it suffices to prove the injectivity
with respect to
\[\{i\star \gamma^0 \mid i\in \I_\bullet\backslash\{\varnothing
\subset (P(w),0)\mid w \in \Sigma^+\}\}.\] The map $i\star \gamma^0$
is of the form $\cyl_\bullet(A,*) \sqcup_{(A,*)} (B,*) \to
\cyl_\bullet(B,*)$ with $i:(A,*)\to (B,*) \in \I_\bullet\backslash
\{\varnothing \subset (P(w),0)\mid w \in \Sigma^+\}$ where the map
from $(A,*)$ to $\cyl_\bullet(A,*)$ is $\gamma^0_{(A,*)}$. Consider
the diagram of solid arrows of $\csts_\bullet$
\[
\xymatrix@C=3em@R=5em {
  \cyl_\bullet(A,*) \sqcup_{(A,*)} (B,*) \fr{} \fd{i\star \gamma^0} & (X,*) \\
  \cyl_\bullet(B,*). \ar@{-->}[ru]^-{\ell} & }
\]
By adjunction, the lift $\ell$ in the diagram above exists if and only
the lift $\ell'$ in the commutative diagram of solid arrows of
$\csts_\bullet$
\[
\xymatrix@C=5em@R=5em {
  (A,*) \fr{\phi'} \fd{i} & \cocyl_\bullet(X,*) \fd{\pi^0}\\
  (B,*) \fr{\psi'} \ar@{-->}[ru]^-{\ell'} & (X,*)}
\]
exists.  By Table~\ref{properness}, the map $i$ is bijective on
states.  By Table~\ref{csts-star-shaped-cyl-path}, the map $\pi^0$ is
bijective on actions.  There is therefore one and exactly one way to
define $\ell'$ on states and on actions.  Let
$(\alpha,u_1,\dots,u_n,\beta)$ be a transition of $B$. Then the tuple
$(\psi'(\alpha),\psi'(u_1),\dots,\psi'(u_n),\psi'(\beta))$ is a
transition of $(X,*)$. By Table~\ref{csts-star-shaped-cyl-path}, we
obtain the equality of tuples
\[\ell'(\alpha,u_1,\dots,u_n,\beta)=((\psi'(\alpha),\alpha'),\psi'(u_1),\dots,\psi'(u_n),(\psi'(\beta),\beta'))\]
with $\alpha'\simeq_{past}\psi'(\alpha)$ and $\beta'
\simeq_{past}\psi'(\beta)$. Since the set of transitions of $(X,*)$ is
closed under past-similarity, the three tuples
\[(\psi'(\alpha),\psi'(u_1),\dots,\psi'(u_n),\beta')),
(\alpha',\psi'(u_1),\dots,\psi'(u_n),\psi'(\beta)),(\alpha',\psi'(u_1),\dots,\psi'(u_n),\beta')\]
are transitions of $X$. Thus, by
Table~\ref{csts-star-shaped-cyl-path}, the tuple
$\ell'(\alpha,u_1,\dots,u_n,\beta)$ is a transition of
$\cocyl_\bullet(X,*)$ and $\ell'$ is a well-defined map of
$\csts_\bullet$.

Let us prove now that $(X,*)$ is injective with respect to the maps
$f\star \gamma$ of $\Lambda_1(\I_\bullet\backslash \{(\varnothing \to
(P(w),0)) \mid w\in \Sigma^+\})$ and of $\Lambda_n(\I_\bullet)$ for
$n\geq 2$. Then we have \[f\in \Lambda_0(\I_\bullet\backslash
\{(\varnothing \to (P(w),0)) \mid w\in \Sigma^+\}) \cup
\Lambda_{n-1}(\I_\bullet).\] By Table~\ref{properness}, the map
$f:(A,*)\to (B,*)$ is bijective on states.  And by
Proposition~\ref{ex0}, it is bijective on actions. By adjunction, we
then have to prove that for any commutative diagram of solid arrows of
$\csts_\bullet$
\[
\xymatrix@C=4em@R=5em
{
(A,*) \fr{\phi} \fd{f} & \fd{\pi} \cocyl_\bullet(X,*) \\
(B,*) \ar@{-->}[ru]^-{\ell} \fr{\psi=(\psi_1,\psi_2)} & (X,*)\p (X,*),
}
\]
the lift $\ell$ exists. Since $f$ is bijective on actions, we have
$\psi_1(u)=\psi_2(u)$ for any action $u$ of $(B,*)$ because $\pi$ is
the diagonal on actions. Since $f$ is bijective on states,
$\psi(\alpha)$ is a state of $\cocyl_\bullet(X,*)$ for any state
$\alpha$ of $(B,*)$. Let $\ell(\alpha)=\psi(\alpha)$ for a state
$\alpha$ of $(B,*)$ and $\ell(u)=\psi_1(u)=\psi_2(u)$. We have to
prove that $\ell$ takes a transition of $(B,*)$ to a transition of
$\cocyl_\bullet(X,*)$. Let $(\alpha,u_1,\dots,u_n,\beta)$ be a
transition of $(B,*)$. Since the map $\psi_1:(B,*)\to (X,*)\p (X,*)
\to (X,*)$ is a map of star-shaped transition systems, the tuple
$(\psi_1(\alpha),\psi_1(u_1),\dots,\psi_1(u_n),\psi_1(\beta))$ is a
transition of $X$. Since the map $\psi_2:(B,*)\to (X,*)\p (X,*) \to
(X,*)$ is a map of star-shaped transition systems, the tuple
$(\psi_2(\alpha),\psi_2(u_1),\dots,\psi_2(u_n),\psi_2(\beta))$ is a
transition of $X$ as well.  Since $\psi(\alpha)$ is a state of
$\cocyl_\bullet(X,*)$, we have $\psi_1(\alpha)\simeq_{past}
\psi_2(\alpha)$ by Table~\ref{csts-star-shaped-cyl-path}. For the same
reason, we have $\psi_1(\beta)\simeq_{past} \psi_2(\beta)$.  Since the
set of transitions of $(X,*)$ is closed under past-similarity, we
deduce that the tuples
$(\psi_2(\alpha),\psi_1(u_1),\dots,\psi_1(u_n),\psi_1(\beta))$ and
$(\psi_1(\alpha),\psi_1(u_1),\dots,\psi_1(u_n),\psi_2(\beta))$ are
transitions of $(X,*)$. Therefore, the tuple
\[(\ell(\alpha),\ell(u_1),\dots,\ell(u_n),\ell(\beta)) =
(\psi(\alpha),\psi_1(u_1),\dots,\psi_1(u_n),\psi(\beta))\] is a
transition of $\cocyl_\bullet(X,*)$ by
Table~\ref{csts-star-shaped-cyl-path}.

It remains to prove that $(X,*)$ is injective with respect to the maps
of $\Lambda_1(\{(\varnothing \to (P(w),0)) \mid w\in \Sigma^+\})$) to
complete the proof. A map of
$\Lambda_1(\{(\varnothing \to (P(w),0)) \mid w\in \Sigma^+\})$ is of
the form $f \star \gamma$ where $f:(P(w),0)\to \cyl_\bullet(P(w),0)$
is the map $f=\gamma^\epsilon_{(P(w),0)}$ for $w\in \Sigma^n$ and
$n\geq 1$. Assume that $\epsilon=0$ without loss of generality. By
adjunction, we have then to prove that for any commutative diagram of
solid arrows of $\csts_\bullet$
\[
\xymatrix@C=5em@R=5em
{
(P(w),0) \fr{\phi} \fd{f=\gamma^0_{(P(w),0)}} & \fd{\pi} \cocyl_\bullet(X,*) \\
\cyl_\bullet(P(w),0) \ar@{-->}[ru]^-{\ell} \fr{\psi=(\psi_1,\psi_2)} & (X,*)\p (X,*),
}
\]
the lift $\ell$ exists. Since $f$ is bijective on actions, and since
$\pi$ is the diagonal map on actions, we have $\psi_1(u)=\psi_2(u)$
for any action $u$ of $(B,*)$. The only possible definition on actions
is $\ell(u)=\psi_1(u)=\psi_2(u)$. Let $\ell(\alpha)=\psi(\alpha)$.
For any $0\leq j\leq n$, we have $(j,0)\simeq_{past} (j,1)$ in
$\cyl_\bullet(P(w),0)$. Consequently, we have
$\psi_1(j,0)\simeq_{past} \psi_1(j,1)$ and $\psi_2(j,0)\simeq_{past}
\psi_2(j,1)$ in $(X,*)$. Since the diagram is commutative, the pair
$(\psi_1(j,0),\psi_2(j,0))$ is a state of $\cocyl_\bullet(X,*)$. By
Table~\ref{csts-star-shaped-cyl-path}, we deduce that $\psi_1(j,0)
\simeq_{past} \psi_2(j,0)$. By Proposition~\ref{pastsim-transitive},
past-similarity is transitive on $(X,*)$.  We deduce that $\psi_1(j,1)
\simeq_{past} \psi_2(j,1)$. We have proved that for any state $\alpha$
of $\cyl_\bullet(P(w),0)$, $\psi(\alpha)$ is a state of
$\cocyl_\bullet(X,*)$ (remember that $\pi$ is one-to-one on
states). It remains to prove that $\ell$ takes a transition of
$\cyl_\bullet(P(w),0)$ to a transition of $\cocyl_\bullet(X,*)$ to
complete the proof. Let $(\alpha,u_1,\dots,u_n,\beta)$ be a transition
of $\cyl_\bullet(P(w),0)$. Then both
$(\psi_1(\alpha),\psi_1(u_1),\dots,\psi_1(u_n),\psi_1(\beta))$ and
$(\psi_2(\alpha),\psi_2(u_2),\dots,\psi_2(u_n),\psi_2(\beta))$ are
transitions of $X$. Since the set of transitions of $(X,*)$ is closed
under past-similarity, the tuples
$(\psi_1(\alpha),\psi_1(u_1),\dots,\psi_1(u_n),\psi_2(\beta))$ and
$(\psi_1(\alpha),\psi_2(u_2),\dots,\psi_2(u_n),\psi_2(\beta))$ are
also transitions of $X$. The proof is complete using
Table~\ref{csts-star-shaped-cyl-path}.  \epf

\section{Characterization in the star-shaped case}
\label{desc-bullet}

\bd A star-shaped transition system is {\rm reduced} when two states are
past-similar if and only if they are equal. \ed

All reduced star-shaped transition systems are fibrant by
Theorem~\ref{fib-carac}.

\bp The full subcategory $\overline{\csts_\bullet}$ of reduced
star-shaped transition systems is a small-orthogonality class of
$\csts_\bullet$. \ep

\bpf Let $w=x_1\dots x_n \in \Sigma^n$ with $n\geq 1$. Let $(C(w),*)$
be the $\omega$-final lift of the map \[\omega(\cyl_\bullet(P(w),0))
\to (\St(\cyl_\bullet(P(w),0))/((n,1)=(n,2)),\{(x_1,1),\dots,(x_n,n)\}).\]
It can be depicted as follows: 
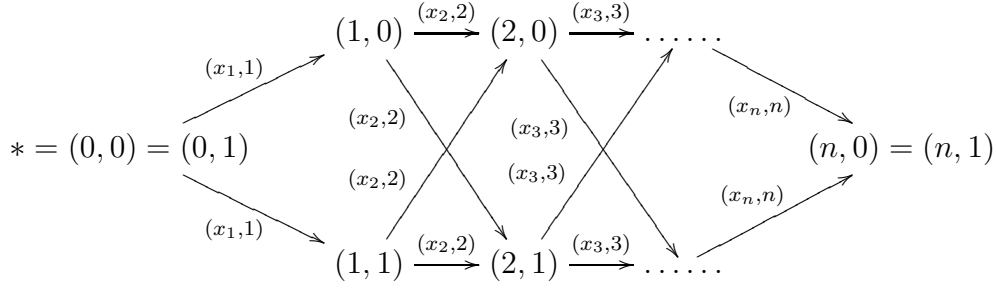
\[
\xymatrix
{
& (1,0) \fr{(x_2,2)}\ar@{->}[rdd]_-/-20pt/{(x_2,2)} & (2,0) \fr{(x_3,3)} \ar@{->}[rdd]_-/-20pt/{(x_3,3)}&\dots\dots   \ar@{->}[rd]_-{(x_n,n)} &\\
{*=(0,0)=(0,1)} \ar@{->}[ru]^-{(x_1,1)} \ar@{->}[rd]_-{(x_1,1)} &&&& (n,0) = (n,1)\\
& (1,1) \fr{(x_2,2)}\ar@{->}[ruu]^-/-20pt/{(x_2,2)} & (2,1) \fr{(x_3,3)}\ar@{->}[ruu]^-/-20pt/{(x_3,3)} &\dots\dots   \ar@{->}[ru]^-{(x_n,n)} & \\
}
\]
A star-shaped transition system is reduced if and only if it is
injective with respect to the maps $\cyl_\bullet(P(w),0) \to (C(w),*)$
with $w\in \Sigma^+$. Since these maps are onto on states and on
actions, they are epic.  Thus, injectivity is equivalent to
orthogonality in this case.  \epf

Let
$\underline{\mathcal{R}} = \{\cyl_\bullet(P(w),0) \to (C(w),*) \mid
w\in \Sigma^+\}$.
For any star-shaped transition system $(X,*)$, the canonical map
$(X,*) \to \mathbf{1}$ factors as a composite
$(X,*) \to R_\bullet^\perp(X,*) \to \mathbf{1}$ with the left-hand map
belonging to $\cell_{\csts_\bullet}(\underline{\mathcal{R}})$ and the
right-hand map belonging to
$\inj_{\csts_\bullet}(\underline{\mathcal{R}})$.  By
Theorem~\ref{colim-csts} and Corollary~\ref{onto-trans-csts}, every
map of $\cell_{\csts_\bullet}(\underline{\mathcal{R}})$ is onto on
states, on actions and on transitions. We deduce that every map of
$\cell_{\csts_\bullet}(\underline{\mathcal{R}})$ is epic in
$\csts_\bullet$. Thus, by \cite[Proposition~A.1]{biscsts1}, this
factorization is unique up to isomorphism. And for the same reason as
for CSA1, this construction provides the left adjoint to the inclusion
functor $\overline{\csts_\bullet} \subset \csts_\bullet$. By
\cite[Theorem~1.39]{MR95j:18001}, the category
$\overline{\csts_\bullet}$ is locally presentable.

\begin{nota} Let us denote by $\Psi_{(X,*)}:(X,*) \to R_\bullet^\perp(X,*)$ the unit
  map. 
\end{nota}

\begin{nota} Let $X$ be a weak transition system. Let $u$ and $v$ be
  two actions of $X$.  Denote by $u\simeq_{\CSA_1} v$ if
  $\mu(u)=\mu(v)$ and if there exist two states $\alpha$ and $\beta$
  of $X$ such that the triples $(\alpha,u,\beta)$ and
  $(\alpha,v,\beta)$ are transitions of $X$.  \end{nota}

\bp \label{onestep} Let $(X,*)$ be a star-shaped cubical transition
system. Let $\underline{r}(X)$ be the $\omega$-final lift of the map
$\omega(X) \longrightarrow
(\St(X)/\!\simeq_{past},\Ac(X)/\!\simeq_{\CSA_1})$ of $\set^{\{s\}\cup
  \Sigma}$. Let $*$ be the image by $X\to \underline{r}(X)$ of $*\in
X$. Then the pointed weak transition system $(\underline{r}(X),*)$ is
a star-shaped cubical transition system.  \ep

\bpf The weak transition system $X$ is cubical by hypothesis. Since
the set map $\Ac(X) \to \Ac(X)/\!\simeq_{\CSA_1}$ is onto, the weak
transition system $\underline{r}(X)$ is cubical by
Theorem~\ref{cubical-lift}.  The map $X\to \underline{r}(X)$ is onto
on states. Thus, every state of $\underline{r}(X)$ is reachable.  \epf

\bp \label{rperp-explicite} Let $(X,*)$ be a star-shaped transition
system. Let $X_0 = X$.  Suppose the weak transition system $X_\xi$
constructed for an ordinal $\xi \geq 0$. Let
$X_{\xi+1}=\underline{r}(X_\xi)$. For a limit ordinal $\xi$, let
$X_\xi = \liminj^\wts_{\beta<\xi} X_\beta$, the colimit being taken in
$\wts$. Then
\begin{enumerate}
\item For all ordinals $\xi$, the weak transition system $X_\xi$ is
  cubical and the pointed cubical transition system $(X_\xi,*)$ is
  star-shaped.
\item There exists an ordinal $\eta$ such that $X_\xi = X_\eta$ for all 
$\xi \geq \eta$.
\item There is the isomorphism $R_\bullet^\perp(X,*)=(X_\eta,*)$. 
\end{enumerate}
\ep 

\bpf The proof is in five steps. 

1) For a limit ordinal $\xi$, the weak transition system
$\liminj^\wts_{\zeta<\xi} X_\zeta$ is cubical if all $X_\zeta$ for
$\zeta < \xi$ are cubical since $\cts$ is a coreflective subcategory
of $\wts$. We have proved the first assertion using
Proposition~\ref{onestep}.

2) The second assertion holds for cardinality reasons. 

3) For any action $u$ and $v$ of $(X_\eta,*)$, we have
$u\simeq_{\CSA_1} v \Rightarrow u=v$ since $X_\eta=X_{\eta+1}$. This
means that $(X_\eta,*)$ satisfies CSA1.

4) Let $(\alpha,u_1,u_2,\beta)$, $(\alpha,u_1,\nu_1)$,
$(\alpha,u_1,\nu_2)$, $(\nu_1,u_2,\beta)$ and $(\nu_2,u_2,\beta)$ be
five transitions of $(X_\eta,*)$.  Thus, the triple
$((\alpha,\alpha),u_1,(\nu_1,\nu_2))$ is a transition of
$\cocyl_*(X_\eta,*)$ by Table~\ref{csts-cyl-path}.  Since
$\alpha\simeq_{past} \alpha$, we deduce that $(\alpha,\alpha)$ is a
reachable state of $\cocyl_*(X_\eta,*)$, and then that $(\nu_1,\nu_2)$
is a reachable state of $\cocyl_*(X_\eta,*)$, and therefore that
$\nu_1\simeq_{past} \nu_2$.  Thus, we obtain $\nu_1=\nu_2$ since
$X_\eta=\underline{r}(X_\eta)$.  Let
\begin{multline*}(\alpha,u_1,\dots,u_n,\beta),(\alpha,u_1,\dots,u_p,\nu_1),
(\alpha,u_1,\dots,u_p,\nu_2),(\nu_1,u_{p+1},\dots,u_n,\beta),\\
(\nu_2,u_{p+1},\dots,u_n,\beta)
\end{multline*} be five transitions of $(X_\eta,*)$ with $n\geq p+1$
and $p\geq 2$. Since $(X_\eta,*)$ is cubical, there exists a state
$\nu_3$ such that the tuples $(\alpha,u_1,\dots,u_{p-1},\nu_3)$ and
$(\nu_3,u_p,\dots,u_n,\beta)$ are two transitions of $(X_\eta,*)$.  By
the patching axiom, the triples $(\nu_3,u_p,\nu_1)$ and
$(\nu_3,u_p,\nu_2)$ are two transitions of $(X_\eta,*)$.  Since
$\nu_3\simeq_{past} \nu_3$, we deduce that $(\nu_3,\nu_3)$ is a
reachable state of $\cocyl_*(X_\eta,*)$, and then that $(\nu_1,\nu_2)$
is a reachable state of $\cocyl_*(X_\eta,*)$, and therefore that
$\nu_1\simeq_{past} \nu_2$.  Thus, we obtain $\nu_1=\nu_2$ since
$X_\eta=\underline{r}(X_\eta)$. Therefore, $(X_\eta,*)$ satisfies
CSA2.

5) We deduce that $(X_\eta,*)$ is a star-shaped Cattani-Sassone
transition system. By construction, the star-shaped transition system
$(X_\eta,*)$ is reduced. Let us prove by induction on the ordinal
$\xi$ that the map $(X,*)\to R_\bullet^\perp(X,*)$ factors as a
composite
$(X,*)\to (X_\xi,*) \stackrel{\phi_\xi}\to R_\bullet^\perp(X,*)$. The
case $\xi=0$ is trivial. If the map $(X,*)\to R_\bullet^\perp(X,*)$
factors as a composite $(X,*)\to (X_\xi,*) \to R_\bullet^\perp(X,*)$,
then for any pair of past-similar states $(\alpha,\beta)$ of $X_\xi$,
we have $\phi_\xi(\alpha) \simeq_{past} \phi_\xi(\beta)$ in
$R_\bullet^\perp(X,*)$. Thus, we obtain
$\phi_\xi(\alpha) = \phi_\xi(\beta)$ since $R_\bullet^\perp(X,*)$ is
reduced. And for any pair of actions $(u,v)$ of $X_\xi$ with
$u\simeq_{\CSA_1} v$, we have
$\phi_\xi(u) \simeq_{\CSA_1} \phi_\xi(v)$ in
$R_\bullet^\perp(X,*)$. Thus, we obtain $\phi_\xi(u) = \phi_\xi(v)$
since $R_\bullet^\perp(X,*)$ satisfies CSA1. We deduce that the map
$\omega(\phi_\xi):\omega(X_\xi) \to
\omega(\omega^*(R_\bullet^\perp(X,*)))$
of $\set^{\{s\}\cup \Sigma}$ factors uniquely as a composite
$\omega(X_\xi) \to \omega(X_{\xi+1})\to
\omega(\omega^*(R_\bullet^\perp(X,*)))$.
We obtain the factorization
$\phi_\xi : (X_\xi,*) \to (X_{\xi+1},*) \to R_\bullet^\perp(X,*)$. By
passing to the colimit, we then obtain the factorization
$X \longrightarrow X_\eta \stackrel{\phi_\eta} \longrightarrow
R_\bullet^\perp(X,*)$.
By the universal property of the adjunction, we deduce the isomorphism
$X_\eta \iso R_\bullet^\perp(X,*)$.  \epf

\bp \label{onto-by-step} Let $f:(X,*) \to (Z,*)$ be a map of
star-shaped cubical transition systems with $(Z,*)$ fibrant in
$\csts_\bullet$. Let $Y=\underline{r}(X)$. Then $(Y,*)$ is a
star-shaped cubical transition system. There exists a map
$g:(Y,*)\to (Z,*)$ of star-shaped cubical transition systems such
that for any state $\overline{\alpha}$ of $X$, we have
$gw(\overline{\alpha})\simeq_{past} f(\overline{\alpha})$ and for any
action $\overline{u}$ of $X$, we have
$gw(\overline{u}) = f(\overline{u})$ where $w:(X,*) \to (Y,*)$ is the
canonical map. \ep

\bpf Let $w:(X,*)\to (Y,*)$ be the canonical map. By
Proposition~\ref{onestep}, $(Y,*)$ is a star-shaped cubical transition
system.  By construction, the map $\omega(w)$ has a section $s$. Let
$g(\alpha) = fs(\alpha)$ for a state $\alpha$ of $Y$ and $g(u) =
fs(u)$ for an action $u$ of $Y$. Let $\overline{\alpha}\in \St(X)$.
Since $wsw(\overline{\alpha})=w(\overline{\alpha})$ in
$\St(X)/\!\simeq_{past}$, the pair of states
$(sw(\overline{\alpha}),\overline{\alpha})$ is in the transitive
closure of the binary relation $\simeq_{past}$ of $(X,*)$. Thus, the
pair of states $(fsw(\overline{\alpha}),f(\overline{\alpha}))$ is in
the transitive closure of the binary relation $\simeq_{past}$ of
$(Z,*)$. But $(Z,*)$ is fibrant by hypothesis. By
Proposition~\ref{pastsim-transitive}, we deduce that
$fsw(\overline{\alpha}) \simeq_{past} f(\overline{\alpha})$.  We have
$gw(\overline{\alpha})=fsw(\overline{\alpha})$ by definition of
$g$. We obtain $gw(\overline{\alpha}) \simeq_{past}
f(\overline{\alpha})$ for any state $\overline{\alpha}$ of $X$. Let
$\overline{u}$ be an action of $X$. By a similar argument, we prove
that the pair of actions $(gw(\overline{u}),f(\overline{u}))$ is in
the transitive closure of the binary relation $\simeq_{\CSA_1}$ of
$(Z,*)$. Since $Z$ satisfies CSA1, we deduce that
$gw(\overline{u})=f(\overline{u})$ for any action $\overline{u}$ of
$X$. It remains to prove that $g$ maps a transition of $(Y,*)$ to a
transition of $(Z,*)$. The weak transition system $Y$ is defined as
the $\omega$-final lift of the map $\omega(X) \to
(\St(X)/\!\simeq_{past},\Ac(X)/\!\simeq_{\CSA_1})$. That is to say, it
is equipped with the final structure. By \cite[Proposition~3.5]{hdts},
this final structure is obtained by considering the set $G_0$ of
transitions which are in the image of the map $(X,*)\to (Y,*)$, then
by applying the patching axiom on the transitions of $G_0$ to obtain a
set $G_1 \supseteq G_0$, and by transfinitely iterating the
process. The set of transitions $\bigcup_{\xi\geq 0} G_\xi$ is the
final structure~\footnote{\cite[Proposition~3.5]{hdts} also claims
  that the multiset axiom is automatically satisfied. This is due to
  the internal symmetry of the patching axiom and to the fact that
  $G_0$ satisfies the multiset axiom.}. We are going to prove by
transfinite induction on $\xi \geq 0$ that for any transition
$(\alpha,u_1,\dots,u_n,\beta)$ of $G_\xi$, the tuple
$(g(\alpha),g(u_1),\dots,g(u_n),g(\beta))$ is a transition of
$(Z,*)$. First of all, let $(\alpha,u_1,\dots,u_n,\beta) \in G_0$. By
definition of $G_0$, there exists a transition
$(\overline{\alpha},\overline{u_1},\dots,\overline{u_n},\overline{\beta})$
of $(X,*)$ such that
$(w(\overline{\alpha}),w(\overline{u_1}),\dots,w(\overline{u_n}),w(\overline{\beta}))
= (\alpha,u_1,\dots,u_n,\beta)$.  We have $w(\overline{\alpha}) =
\alpha = ws(\alpha)$ since $s$ is a section of $w$ on states. Thus,
the pair of states $(\overline{\alpha},s(\alpha))$ is in the
transitive closure of the binary relation $\simeq_{past}$ in
$(X,*)$. We obtain $f(\overline{\alpha}) \simeq_{past}
fs(\alpha)=g(\alpha)$ since $\simeq_{past}$ is transitive in $(Z,*)$
by Proposition~\ref{pastsim-transitive}. For the same reason, we
obtain $f(\overline{\beta}) \simeq_{past} fs(\beta) =g(\beta)$. The
tuple
$(f(\overline{\alpha}),f(\overline{u_1}),\dots,f(\overline{u_n}),f(\overline{\beta}))$
is a transition of $(Z,*)$ since $f$ is a map of transition
systems. By Theorem~\ref{fib-carac}, and since $(Z,*)$ is fibrant by
hypothesis, the tuple
$(g(\alpha),f(\overline{u_1}),\dots,f(\overline{u_n}),g(\beta))$ is
then a transition of $(Z,*)$. We have $w(\overline{u_i})=u_i=ws(u_i)$
for all $1\leq i \leq n$. Thus, the pair of actions
$(\overline{u_i},s(u_i))$ for any $1\leq i \leq n$ is in the
transitive closure of the binary relation $\simeq_{\CSA_1}$ in
$(X,*)$. We obtain $f(\overline{u_i}) = fs(u_i)=g(u_i)$ for all $1\leq
i \leq n$ since $Z$ satisfies CSA1. Therefore, the tuple
$(g(\alpha),g(u_1),\dots,g(u_n),g(\beta))$ is a transition of $(Z,*)$.
The step $\xi=0$ of the transfinite induction is proved. The case
$\xi$ limit ordinal is trivial. It remains to prove that if all
transitions of $G_\xi$ are mapped by $g$ to transitions of $(Z,*)$,
then the same fact holds for $G_{\xi +1}$. Consider the five tuples
\begin{multline*}
(\alpha,u_1, \dots, u_n, \beta), (\alpha,u_1, \dots, u_p, \nu_1),
(\nu_1, u_{p+1}, \dots, u_n, \beta),\\ (\alpha, u_1, \dots, u_{p+q},
\nu_2), (\nu_2, u_{p+q+1}, \dots, u_n, \beta) 
\end{multline*}
of $G_\xi$ with
$n\geq 3$, $p,q\geq 1$ and $p+q<n$. By definition, the tuple $(\nu_1,
u_{p+1}, \dots, u_{p+q}, \nu_2)$ belongs to $G_{\xi+1}$. By induction
hypothesis, the five tuples \begin{multline*}(g(\alpha),g(u_1), \dots, g(u_n),
g(\beta)), (g(\alpha),g(u_1), \dots, g(u_p), g(\nu_1)), (g(\nu_1), \\
g(u_{p+1}), \dots, g(u_n), g(\beta)), (g(\alpha), g(u_1), \dots,
g(u_{p+q}), g(\nu_2)),\\ (g(\nu_2), g(u_{p+q+1}), \dots, g(u_n),
g(\beta))\end{multline*} 
are transitions of $(Z,*)$.  By applying the patching axiom
in $(Z,*)$, we obtain that the tuple $(g(\nu_1), g(u_{p+1}), \dots,
g(u_{p+q}), g(\nu_2))$ is a transition of $(Z,*)$.  The proof is
complete.
\epf

\bth \label{loc-red} For any star-shaped transition system $(X,*)$, the map
$\Psi_{(X,*)}:(X,*) \to R_\bullet^\perp(X,*)$ is a weak equivalence of
$\csts_\bullet$. \eth

\bpf Let $f:(X,*) \to (Y,*)$ be a map of star-shaped transition
systems.  By Proposition~\ref{rperp-explicite} and
Proposition~\ref{onto-by-step}, there exists a map of star-shaped
transition system $g:R_\bullet^\perp(X,*) \to (Y,*)$ such that for any state
$\alpha$ of $(X,*)$, we have $g\Psi_{(X,*)}(\alpha) \simeq_{past}
f(\alpha)$ and for any action $u$ of $(X,*)$, we have
$g\Psi_{(X,*)}(u) = f(u)$.  By Theorem~\ref{car_homeq-2}, we deduce
that $g\Psi_{(X,*)}$ and $f$ are homotopy equivalent maps.  Thus, the
set map \[\pi_{\csts_\bullet}(R_\bullet^\perp(X,*),(Y,*)) \to
\pi_{\csts_\bullet}((X,*),(Y,*))\] induced by the precomposition with
$\Psi_{(X,*)}:(X,*) \to R_\bullet^\perp(X,*)$ is onto.  Let
$f,g:R_\bullet^\perp(X,*)\to (Z,*)$ be two maps of $\csts_\bullet$ such that
$f\Psi_{(X,*)}$ is homotopy equivalent to $g\Psi_{(X,*)}$. Then
$f\Psi_{(X,*)}$ and $g\Psi_{(X,*)}$ coincide on actions by
Theorem~\ref{car_homeq-2} and for any state $\alpha$ of $X$, the
states $f(\Psi_{(X,*)}(\alpha))$ and $g(\Psi_{(X,*)}(\alpha))$ are
past-similar.  For any state $\beta$ of $Y$, there exists a state
$\alpha$ of $X$ such that $\Psi_{(X,*)}(\alpha)=\beta$. Thus, for any
state $\beta$ of $Y$, the states $f(\beta)$ and $g(\beta)$ are
past-similar. Since $\Psi_{(X,*)}$ is onto on actions as well, $f$ and
$g$ coincide on actions.  We deduce that $f$ and $g$ are homotopy
equivalent by Theorem~\ref{car_homeq-2}. We deduce that the set map
\[\pi_{\csts_\bullet}(R_\bullet^\perp(X,*),(Z,*)) \to
\pi_{\csts_\bullet}((X,*),(Z,*))\]
induced by the precomposition with $\Psi_{(X,*)}$ is one-to-one.  \epf

\bth \label{carac-bullet} The left adjoint $R_\bullet^\perp:\csts_\bullet \to
\overline{\csts_\bullet}$ induces a Quillen equivalence between the
model category $\csts_\bullet$ and the category
$\overline{\csts_\bullet}$ equipped with the discrete model category
structure. \eth

\bpf For any reduced star-shaped transition system $(X,*)$, we have by
Proposition~\ref{pastsimilar-pathsimilar} and by
Table~\ref{csts-cyl-path} the equality
$\cocyl_\bullet(X,*)=(X,*)$. Let $\alpha$ be a state of $X$. Since
$(X,*)$ is star-shaped, there exists $w\in \Sigma^n$ with $n\geq 0$
and a map $f:(P(w),0) \to (X,*)$ such that $f(n)=\alpha$. By
functoriality, we obtain a map $\cyl_\bullet(f):\cyl_\bullet(P(w),0)
\to \cyl_\bullet(X,*)$. Thus, we have $(\alpha,0) \simeq_{past}
(\alpha,1)$ in $\cyl_\bullet(X,*)$. This implies that
$R_\bullet^\perp(\cyl(X,*))=(X,*)$ for any reduced star-shaped
transition system $(X,*)$. Using \cite[Theorem~3.1]{leftdet}, we
obtain a left determined Olschok model structure on
$\overline{\csts_\bullet}$ such that the cylinder and path functors
are the identity functor.  Therefore, we have the equalities
\begin{multline*}\pi_{\overline{\csts_\bullet}}((X,*),(Y,*)) =
  \pi^l_{\overline{\csts_\bullet}}((X,*),(Y,*))\\ =
  \pi^r_{\overline{\csts_\bullet}}((X,*),(Y,*)) =
  \overline{\csts_\bullet}((X,*),(Y,*))\end{multline*} for any reduced
Cattani-Sassone transition system $(X,*)$ and $(Y,*)$. This means that
the weak equivalences of $\overline{\csts_\bullet}$ are the
isomorphisms.  Thus, the left adjoint $R_\bullet^\perp:\csts_\bullet
\to \overline{\csts_\bullet}$ induces a homotopically surjective left
Quillen adjoint from $\csts_\bullet$ to $\overline{\csts_\bullet}$
equipped with the discrete model structure.  By Theorem~\ref{loc-red},
this left Quillen adjoint is a left Quillen equivalence.  \epf

The following corollary proves that a map of star-shaped transition
systems is a weak equivalence if and only if it becomes an isomorphism
after the identification of all past-similar states.

\begin{cor} \label{red-2} A map $f$ of $\csts_\bullet$ is a weak
  equivalence if and only if $R_\bullet^\perp(f)$ is an
  isomorphism. \end{cor}

\bpf Since all objects of $\csts_\bullet$ are cofibrant, a weak
equivalence $f$ is mapped to a weak equivalence $R_\bullet^\perp(f)$
of $\overline{\csts_\bullet}$, i.e. an isomorphism.  Conversely, if
$R_\bullet^\perp(f)$ is an isomorphism, then by Theorem~\ref{loc-red}
and the two-out-of-three property, $f$ is a weak equivalence.
\epf

\section{Causality and homotopy}
\label{futur}

This concluding section is written to interpret
Theorem~\ref{carac-bullet}.  All left determined model categories
constructed so far on higher dimensional transition systems, including
the ones of \cite{cubicalhdts} and \cite{csts} where
$R:\{0,1\}\to \{0\}$ is a generating cofibration, are Quillen
equivalent to discrete model categories. Similar left determined model
categories on flows \cite{model3} and multipointed $d$-spaces
\cite{mdtop} do not have such a behavior. This phenomenon could be
related to the absence of thin objects in higher dimension in the
formalism of higher dimensional transition systems which categorically
behave like labelled symmetric precubical sets
\cite[Theorem~11.6]{hdts}.

To be more specific in the sequel, we will be using the
$1$-dimensional paths $(P(w),0)$ with $w\in \Sigma^*$ (where
$\Sigma^*= \bigcup_{n\geq 0} \Sigma^n$ is the set of words over
$\Sigma$). The arguments developed here could be adapted to more
complicated notions of paths, in particular higher dimensional ones
like in \cite{zbMATH06567559}. Let \[\mathcal{P}=\{(P(w),0) \to
(P(ww'),0) \mid w,w'\in \Sigma^*\}\] be the set of extensions of paths.
The semantics of \cite{hdts} is used in this section. The reader does
not actually need to read the latter paper to understand the
sequel. Indeed, except for Proposition~\ref{debut} whose proof is just
sketched, the only facts to know are that:
\begin{enumerate}
\item All $\mathcal{P}$-cell complexes are realizations of process
  algebras.
\item All realizations of process algebras are colimits of cubes.
\end{enumerate}

\bd After \cite{0856.68067}, two star-shaped transition systems
$(X,*)$ and $(Y,*)$ are \emph{$\mathcal{P}$-bisimilar} if they are
related by a span of $\mathcal{P}$-injective maps $(X,*) \leftarrow
(Z,*) \rightarrow (Y,*)$. \ed

Since the class of $\mathcal{P}$-injective maps is closed under
pullback and composition, two star-shaped transition systems $(X,*)$
and $(Y,*)$ are $\mathcal{P}$-bisimilar if and and only if they are
related by a zig-zag of $\mathcal{P}$-injective maps. Note that a
$\mathcal{P}$-injective map between star-shaped transition systems is
always onto on states, on actions and on $1$-dimensional transitions.

\bp \label{debut} Any two weakly equivalent star-shaped transition
systems of the left determined Olschok model category $\csts_\bullet$
realizing process algebras are isomorphic, and hence
$\mathcal{P}$-bisimilar. There exist two $\mathcal{P}$-bisimilar
$\mathcal{P}$-cell complexes which are not weakly equivalent.  \ep

\bpf[Sketch of proof] A star-shaped transition systems realizing a
process algebra is reduced: the proof is by induction on the syntactic
description of the process algebra. Thus, if two of them are weakly
equivalent, they are isomorphic by Corollary~\ref{red-2}.  The two
$\mathcal{P}$-cell complexes $(P(u),0)\sqcup (P(u),0)$ and $(P(u),0)$
are $\mathcal{P}$-bisimilar since the unique map $(P(u),0)\sqcup
(P(u),0) \to (P(u),0)$ is $\mathcal{P}$-injective. They can be
depicted as follows with $\mu(u_1)=\mu(u_2)=u$:
\[
\xymatrix@C=1em@R=1em
{
&& \bullet   \\
(P(u),0)\sqcup (P(u),0) = & 0 \ar@{->}[ru]^-{u_1}\ar@{->}[rd]_-{u_2}&&&&& (P(u),0) = & 0 \fr{u_1} &\bullet  \\
&& \bullet \\
}
\]
However, they are reduced and not isomorphic. Thus, they are not
weakly equivalent by Corollary~\ref{red-2}. \epf

\bd A model structure on $\csts_\bullet$ is {\rm $\mathcal{P}$-causal}
(with respect to the semantics of \cite{hdts}) if any two star-shaped
transition systems realizing process algebras are weakly equivalent if
and only if they are $\mathcal{P}$-bisimilar.  \ed

\bth \label{big-pb} Consider a model structure of $\csts_\bullet$ such
that all $\mathcal{P}$-cell complexes are cofibrant and such that all
$\mathcal{P}$-injective maps are weak equivalences. Then there exist
two homotopy equivalent $\mathcal{P}$-cell complexes which are not
$\mathcal{P}$-bisimilar. In particular, this model structure is not
$\mathcal{P}$-causal. \eth

\bpf Consider the two $\mathcal{P}$-cell complexes $(M_0,0)=(P(uv),0)
\sqcup (P(u),0)$ and $(M_1,0)=(P(uv),0)$.  The star-shaped transition
systems $(M_0,0)$ and $(M_1,0)$ look as follows (with
$\mu(u_1)=\mu(u_2)=u$ and $\mu(v_1)=v$):
\[
\xymatrix@C=1em@R=1em
{
&& \bullet \fr{v_1} & \bullet  \\
(M_0,0) = & 0 \ar@{->}[ru]^-{u_1}\ar@{->}[rd]_-{u_2}&&&&& (M_1,0) = & 0 \fr{u_1} &\bullet \fr{v_1}&\bullet \\
&& \bullet \\
}
\]
They are not $\mathcal{P}$-bisimilar since the path
$0\stackrel{u_2}\to\bullet$ of $(M_0,0)$ cannot be extended. Let
$f:(M_0,0) \to (M_1,0)$ be the unique map defined on actions by the
mappings $u_1,u_2\mapsto u_1$ and $v_1 \mapsto v_1$. Let $g:(M_1,0)
\to (M_0,0)$ be the unique map defined on actions by the mappings
$u_1\mapsto u_1$ and $v_1 \mapsto v_1$. Consider a commutative diagram
of solid arrows of $\csts_\bullet$
\[
\xymatrix@C=4em@R=4em
{
(P(w),0) \fr{\phi} \fd{} & (M_0,0) \sqcup (M_0,0) \fd{} \\
(P(ww'),0) \fr{} \ar@{-->}[ru]^-{\ell} &(M_0,0)
}
\]
where $w,w'\in \Sigma^*$. The only possibilities for $w,ww'\in
\Sigma^*$ are $w,ww'\in \{\varnothing,u,uv\}$. Consequently, the map
$\phi$ factors as a composite
\[(P(w),0) \longrightarrow (M_0,0) \longrightarrow (M_0,0) \sqcup
(M_0,0).\] Thus, the lift $\ell$ exists. We deduce that the codiagonal
map $(M_0,0)\sqcup (M_0,0) \to (M_0,0)$ factors as a
composite \[\xymatrix@C=5em{(M_0,0)\sqcup (M_0,0) \fr{\iso}&
  (M_0,0)\sqcup (M_0,0) \fr{\in\inj_{\csts_\bullet}(\mathcal{P})}&
  (M_0,0)}\] where the left-hand map is a cofibration and the
right-hand map is $\mathcal{P}$-injective, i.e. by hypothesis a weak
equivalence. Consequently, $(M_0,0)\sqcup (M_0,0)$ is a good cylinder
of $(M_0,0)$ for this model structure. We deduce that the maps
$\id_{(M_0,0)},gf:(M_0,0)\rightrightarrows (M_0,0)$ are left homotopic
maps. By \cite[Proposition~7.4.8]{ref_model2}, and since $(M_0,0)$ is
cofibrant by hypothesis, $\id_{(M_0,0)}$ and $gf$ are right homotopic,
and then homotopic.  For the same reason, the maps $\id_{(M_1,0)}$ and
$fg$ are homotopic. Therefore, the star-shaped transition systems
$(M_0,0)$ and $(M_1,0)$ are homotopy equivalent in this model
structure.  \epf

Theorem~\ref{big-pb} and its proof tell us that localizing with
respect to the whole class of $\mathcal{P}$-injective maps is a very
bad idea. Theorem~\ref{big-pb} also tells us that the most obvious
candidate, the left determined model structure with respect to
$\mathcal{P}$, whose existence is a consequence of Vop\v{e}nka's
principle by \cite[Theorem~2.2]{rotho}, is not $\mathcal{P}$-causal
either. 

Theorem~\ref{big-pb} also holds by replacing the category of
star-shaped transition systems by any category of labelled precubical
sets of \cite{ccsprecub} or \cite{symcub}.  Indeed, the origin of the
problem is that the map $(X,*)\sqcup (X,*)\to (X,*)$ is
$\mathcal{P}$-injective for any star-shaped transition system $(X,*)$
not containing any cycle passing by the base state $*$, which means in
this case that a good cylinder of $(X,*)$ is $(X,*)\sqcup (X,*)$ in
such a model structure.

It turns out that there exist star-shaped transition systems which are
not colimits of cubes, e.g. the star-shaped transition systems of
Figure~\ref{not-past-similar} and Figure~\ref{past-similar}. It is
actually the main technical difference with any category of labelled
precubical sets of \cite{ccsprecub} or \cite{symcub}. To overcome the
problem arising from Theorem~\ref{big-pb}, the idea is to localize
with respect to a class of $\mathcal{P}$-injective maps between
star-shaped transition systems which are $\mathcal{P}$-optimized in
the sense that they use a minimal set of actions. For example, with
$(X,*)=(P(uv),*)$, the star-shaped transition system $(X,*)\sqcup
(X,*)$ looks as follows (with $\mu(u_1)=\mu(u_2)=u$ and
$\mu(v_1)=\mu(v_2)=v$):
\[
\xymatrix@C=1em@R=1em
{
&& \bullet \fr{v_1} & \bullet  \\
(X,*)\sqcup (X,*) = & {*} \ar@{->}[ru]^-{u_1}\ar@{->}[rd]_-{u_2} \\
&& \bullet \ar@{->}[r]_-{v_2}&\bullet\\
}
\]
and its $\mathcal{P}$-optimized version is:
\[
\xymatrix@C=1em@R=1em
{
&& \bullet \fr{v} & \bullet  \\
(Y,*) = & {*} \ar@{->}[ru]^-{u}\ar@{->}[rd]_-{u} \\
&& \bullet \ar@{->}[r]_-{v}&\bullet\\
}
\]
The star-shaped transition system of $(Y,*)$ is exactly the one of
Figure~\ref{not-past-similar}.  The star-shaped transition systems
$(X,*)\sqcup (X,*)$ and $(Y,*)$ are $\mathcal{P}$-bisimilar and
$(Y,*)$ uses as few actions as possible.  Note that $(X,*)\sqcup
(X,*)$ is never $\mathcal{P}$-optimized unless $X=\{*\}$.

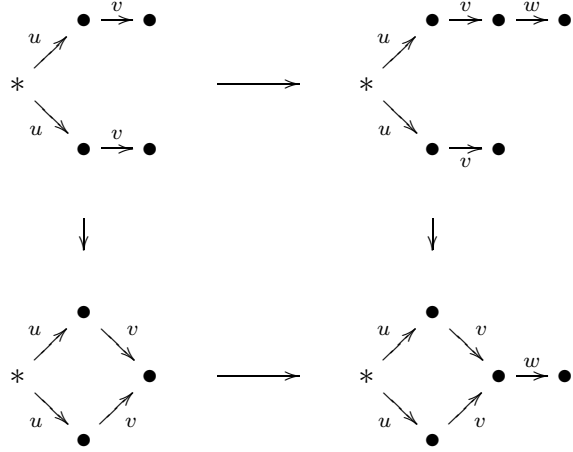
\begin{figure}
\[
\xymatrix@C=1em@R=1em
{
&\bullet \fr{v} & \bullet &&&& &\bullet   \fr{v} & \bullet \fr{w} & \bullet           \\
{*} \ar@{->}[ru]^-{u}\ar@{->}[rd]_-{u}&&&\fR{}&&&{*} \ar@{->}[ru]^-{u}\ar@{->}[rd]_-{u}\\
&\bullet \fr{v} & \bullet && &&   &\bullet   \ar@{->}[r]_-{v} & \bullet    \\
& \fd{}&&&&&& \fd{}\\
&&&&&&&\\
&\bullet \ar@{->}[rd]^-{v}&&&&&&\bullet \ar@{->}[rd]^-{v}\\
{*} \ar@{->}[ru]^-{u}\ar@{->}[rd]_-{u}&& \bullet& \fR{}&&&{*} \ar@{->}[ru]^-{u}\ar@{->}[rd]_-{u}&& \bullet \fr{w} &\bullet \\
&\bullet \ar@{->}[ru]_-{v}&&&&&&\bullet \ar@{->}[ru]_-{v}\\
}
\]
\caption{A pushout of a $\mathcal{P}$-injective map which is not $\mathcal{P}$-injective}
\label{left-proper-pb}
\end{figure}

The pushout diagram of Figure~\ref{left-proper-pb} highlights another
problem which seems to indicate that left properness could be an
obstacle. All star-shaped transition systems of
Figure~\ref{left-proper-pb} are $\mathcal{P}$-optimized. The left
vertical map is $\mathcal{P}$-injective. The right vertical map is not
$\mathcal{P}$-injective since the bottom path of the domain
$*\stackrel{u}\to \bullet \stackrel{v}\to \bullet$ cannot be extended.
It turns out that the right vertical map identifies two states having
the same past and not the same future.  There are two ways of
overcoming this situation: 1) noticing that the domain and the
codomain of the bottom horizontal arrow are not the optimizations of
star-shaped transition systems coming from a process algebra; indeed
the semantics of \cite{hdts} cannot create directed cycles; 2)
colocalizing with respect to the set of paths; then the domain and the
codomain of the bottom horizontal arrow are not cofibrant anymore.

After all these observations, we want to study localizations with
respect to $\mathcal{P}$-injective maps between
$\mathcal{P}$-optimized realizations of process algebras and also
localizations of the colocalization with respect to the set of paths
of $\csts_\bullet$. The colocalization is not left proper and entails
the introduction of a new underlying category.

\appendix

\section{Erratum}
\label{erratum}

\begin{figure}
\[
\xymatrix@C=3em@R=0.2em
{
\stackrel{\beta}\bullet & \fl{u} \stackrel{\alpha}\bullet & \stackrel{\gamma}\bullet \\
& \cap & \\
\stackrel{\beta}\bullet & \fl{u} \stackrel{\alpha}\bullet & \stackrel{\gamma}\bullet \\
 & \fr{u} \stackrel{\alpha'}\bullet & \stackrel{\gamma'}\bullet \\
& \cap & \\
\stackrel{\beta}\bullet & \fl{u} \stackrel{\alpha}\bullet \fr{u} &\stackrel{\gamma}\bullet
}
\]
\caption{Erratum}
\label{inclusions}
\end{figure}

\cite[Theorem~A.2]{biscsts1} is false. It is used in
\cite[Proposition~2.7]{biscsts1} and in \cite[Theorem~3.7]{biscsts1}.
\cite[Proposition~2.7]{biscsts1} is still true.
\cite[Theorem~3.7]{biscsts1} is not true. Indeed, 
the cofibration 
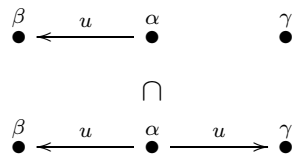
\[
\xymatrix@C=3em@R=0.2em
{
\stackrel{\beta}\bullet & \fl{u} \stackrel{\alpha}\bullet & \stackrel{\gamma}\bullet \\
& \cap & \\
\stackrel{\beta}\bullet & \fl{u} \stackrel{\alpha}\bullet \fr{u} & \stackrel{\gamma}\bullet
}
\]
is not a transfinite composition of pushouts of maps of 
\begin{multline*}
I^{\cts} = \{C:\varnothing \to \{0\}\} \\ \cup \{\de C_n[x_1,\dots,x_n]
\to C_n[x_1,\dots,x_n]\mid \hbox{$n\geq 1$ and $x_1,\dots,x_n \in
  \Sigma$}\} \\ \cup \{C_1[x] \to \dd{x}\mid x\in \Sigma\}.
\end{multline*}
However, it belongs to $\cell_\cts(I^{\cts} \cup \{R:\{0,1\} \to
\{0\}\})$ by the sequence of inclusions of Figure~\ref{inclusions}. The
correct statement of \cite[Theorem~3.7]{biscsts1} is obtained by
replacing the maps $C_1[x]\to \dd{x}$ for $x$ running over $\Sigma$ by
the maps $C_0 \sqcup C_0 \sqcup C_1[x] \to \dd{x}$ for $x$ running
over $\Sigma$ which are defined to be bijective on states (see
Notation~\ref{fix} of this paper).

\end{document}